\documentclass[11pt]{article}
\usepackage{amsfonts}
\begin{document}

\setcounter{secnumdepth}{5}
\renewcommand{\labelenumi}{(\roman{enumi}). }

\newcommand{\R}{I\hspace{-1.5 mm}R} 
\newcommand{\N}{I\hspace{-1.5 mm}N}
\newcommand{\1}{I\hspace{-1.5 mm}I}
\newtheorem{theorem}{Theorem}[section]
\newtheorem{lemma}[theorem]{Lemma}
\newtheorem{prop}[theorem]{Proposition}
\newtheorem{coro}[theorem]{Corollary}
\newtheorem{rem}[theorem]{Remark}
\newtheorem{ex}[theorem]{Example}
\newtheorem{claim}[theorem]{Claim}
\newtheorem{conj}[theorem]{Conjecture}

\newcounter{mycounter}[mycounter] \setcounter{mycounter}{1} 
\newcounter{ii}[mycounter] \setcounter{ii}{2} 
\newcounter{iii}[mycounter] \setcounter{iii}{3} 
\newcounter{iv}[mycounter] \setcounter{iv}{4}
\newcounter{v}[mycounter] \setcounter{v}{5}
\newcounter{vi}[mycounter] \setcounter{vi}{6}
\newcounter{vii}[mycounter] \setcounter{vii}{7}
\newcounter{viii}[mycounter] \setcounter{viii}{8}
\newcounter{ix}[mycounter] \setcounter{ix}{9}
\newcounter{x}[mycounter] \setcounter{x}{10}

\title{Existence of solutions to degenerate parabolic equations via the Monge-Kantorovich theory.}
\author{Martial Agueh}
\date{ School of Mathematics\\
Georgia Institute of Technology\\
Atlanta, GA 30332, USA\\
agueh@math.gatech.edu}
%\date{\today}
\maketitle

\begin{abstract}
We obtain solutions of the nonlinear degenerate parabolic equation
\[ \frac{\partial\,\rho}{\partial\,t} \,=\, \mbox{div}\,\Big\{ \rho\,\nabla c^\star\,\left[\,\nabla \left(F^\prime(\rho)+V\right)\,\right] \Big\} \]
as a steepest descent of an energy with respect to a convex cost functional. The method used here is variational. It requires less uniform convexity assumption than that imposed by Alt and Luckhaus in their pioneering work \cite{luckhaus:quasilinear}. In fact, their assumption may fail in our equation. This class of problems includes the Fokker-Planck equation, the porous-medium equation, the fast diffusion equation, and the parabolic p-Laplacian equation.
\end{abstract}

energy inequality.

\tableofcontents

\section{Introduction}
We consider a class of parabolic evolution equations, so-called doubly degenerate parabolic equations. These equations arise in many applications in physics and biology \cite{gardiner:handbook}, \cite{peletier:porous}, \cite{risken:fokker}, \cite{schuss:singular}. They are used to model a variety of physical problems: the evolution of a fluid in a certain domain: porous-medium equation \cite{otto:porous}, Fokker-Planck equation \cite{otto:variation}, etc. In this work, we focus on  these parabolic equations of the form

\begin{equation}
\label{intro1}
\left\{ \begin{array}{lcl}
\frac{\partial\,b(u)}{\partial\,t} = \mbox{div}\,\left(a\left(b(u), \nabla u\right)\right) &\mbox{on}& (0,\infty)\times \Omega \\ \\
u(t=0) = u_0 &\mbox{on}& \Omega\\ \\
a\left(b(u),\nabla u\right)\cdot\nu = 0 &\mbox{on}& (0,\infty)\times\partial\,\Omega, \end{array} \right.
\end{equation}
where 
\[ a\left(b(u),\nabla u\right) := f\left(b(u)\right)\nabla c^\star\left[\,\nabla(u+V)\,\right],\]
and $c^\star$ denotes the Legendre transform of a function $\,c:\R^d\rightarrow [0,\infty),\,$ that is, 
\[ c^\star(z)=\sup_{x\in\R^d} \{\langle x,z\rangle- c(x)\},\]
for $z\in\R^d$. Here, $\Omega$ is a bounded domain of $\R^d,\, \nu$ is the outward unit normal to $\partial \Omega$, $b:\R\rightarrow \R$ is a monotone nondecreasing function, $V:\bar{\Omega}\rightarrow \R$ is a potential, $c:\R^d \rightarrow [0,\infty)$ is a convex function, $f$ is a nonnegative real-valued function, and $u_0:\Omega\rightarrow \R$ is a measurable function. The unknown is $u:[0,\infty)\times\Omega \rightarrow \R, \quad u=u(t,x).$\\
\indent In a previous work, Alt and Luckhaus \cite{luckhaus:quasilinear} proved existence of weak solutions to (\ref{intro1}), when $V=0$, under the following ellipticity condition on $a(t,z):= f(t)\nabla c^\star(z)$,
\begin{equation}
\label{intro2}
\langle a(t,z_1)-a(t,z_2), z_1-z_2\rangle \geq \lambda\,|\,z_1-z_2\,|^p, 
\end{equation}
for some $\lambda>0\,$ and $\,p\geq 1$, and for all $z_1,\,z_2\in\R^d$. This amounts to imposing that $f$ is bounded below, and the cost function $c$ satisfies the ellipticity condition
\begin{equation}
\label{intro3}
\langle \nabla c^\star(z_1)-\nabla c^\star(z_2), z_1-z_2\rangle \geq \lambda\,|\,z_1-z_2\,|^p.
\end{equation}
Note that when $c(z)=\frac{|\,z\,|^q}{q}$ or equivalently $c^\star(z)=\frac{|\,z\,|^p}{p}$ where $p>1$ is the conjugate of $q>1$, $\frac{1}{p}+\frac{1}{q}=1$, condition (\ref{intro3}) reads as 
\begin{equation}
\label{intro3.1}
\langle |\,z_1\,|^{p-2}z_1-|\,z_2\,|^{p-2}z_2, z_1-z_2\rangle \geq \lambda |\,z_1-z_2\,|^p,
\end{equation}
which holds only if $p\geq 2$. In fact, when $1<p<2$, the reverse inequality in (\ref{intro3.1}) holds (see \cite{di}, pp. 13).  
In \cite{luckhaus:quasilinear}, the authors approximated (\ref{intro1}) by a time discretization, and they used a Galerkin type argument to solve the resulting elliptic problems. In the same paper, they proved uniqueness of solutions to (\ref{intro1}) when $V=0$, assuming that (\ref{intro2}) holds, and the distributional derivative $\frac{\partial b(u)}{\partial t}$ of a solution $u$ of (\ref{intro1}) is an integrable function. The last condition was removed by Otto in \cite{otto:contraction}, using the technique of ``doubling of variables'', introduced by Kru\v{z}kov \cite{kruzkov:first}, which consists of doubling the time variable of two solutions of (\ref{intro1}), and treating each solution as a constant with respect to the differential equation satisfied by the other solution.

In this work, we eliminate assumption (\ref{intro3}), and we impose instead, the following growth condition on the function $c$:
\begin{equation}
\label{intro5}
 \beta\,|z|^q \leq c(z) \leq \alpha\,\left(|z|^q+1\right),
\end{equation}
for $z\in\R^d$ and for some $\alpha,\,\beta >0$ and $q>1$. Notice that (\ref{intro5}) is much weaker than the ellipticity condition (\ref{intro3}) imposed by Alt and Luckhaus in \cite{luckhaus:quasilinear}. Typical examples are the functions $c(z)=\frac{|\,z\,|^q}{q}$ or equivalently $c^\star(z)=\frac{|\,z\,|^p}{p}$ with $\frac{1}{p}+\frac{1}{q}=1$ and $1<p<2$. Such functions satisfy (\ref{intro5}) but not (\ref{intro3}) or (\ref{intro3.1}) when $1<p<2$. As mentioned before, they actually satisfy the reverse inequality in (\ref{intro3}) or (\ref{intro3.1}) when $1<p<2$.

 We interpret (\ref{intro1}) as a dissipative system, and then, we introduce the internal energy density function $F:[0,\infty)\rightarrow \R$, satisfying $F^\prime = b^{-1}$. Setting $\rho:=b(u), \; \rho_0:=b(u_0)$, and $f(x)=\max(x,0)$, we rewrite (\ref{intro1}) as 
\begin{equation}
\label{intro8}
\left\{ \begin{array}{lcl}
\frac{\partial\,\rho}{\partial\,t} + \mbox{div}\,(\rho\,U_\rho)=0 &\mbox{on}& (0,\infty)\times \Omega \\ \\
\rho(t=0) = \rho_0 &\mbox{on}& \Omega\\ \\
\rho\, U_\rho\cdot \nu = 0 &\mbox{on}& (0,\infty)\times\partial\,\Omega . \end{array} \right.
\end{equation}
Here,
\[ U_\rho:= -\nabla c^\star\left[\,\nabla\left(F^\prime(\rho)+V\right)\,\right]\]
denotes the vector field describing the average velocity of a fluid evolving with the continuity equation in (\ref{intro8}), $\rho_0:\Omega\rightarrow [0,\infty)$ is the initial mass density of the fluid, and the unknown $\rho:[0,\infty)\times\Omega\rightarrow [0,\infty),\;\rho=\rho(t,x)$, is the mass density of the fluid at time $t$ and position $x$ of $\Omega$. The free energy associated with the fluid at time $t\in [0,\infty)$, is the sum of its internal energy and its potential energy,
\[ E(\rho(t)):=\int_\Omega \left[\,F\left(\rho(t,x)\right)+\rho(t,x)V(x)\,\right]\,\mbox{d}x.\]
Problem (\ref{intro8}) includes the
\begin{itemize}
\item Linear Fokker-Planck equation:
\[\frac{\partial\rho}{\partial t} = \Delta\rho +\mbox{div}(\rho\nabla V)\]
($c(z) = \frac{|\,z\,|^2}{2}$ and $F(x)=x\ln x$)
\item Porous-medium and Fast diffusion equations:
 \[\frac{\partial\rho}{\partial t} = \Delta\rho^m\]
($V=0$, $c(z) = \frac{|\,z\,|^2}{2}$, and $F(x)=\frac{x^m}{m-1}$ with $1\neq m\geq 1-\frac{1}{d}$).

\item Generalized heat equation:
\[\frac{\partial\rho}{\partial t} = \mbox{div}\,\left(|\,\nabla\rho^{\frac{1}{p-1}}\,|^{p-2}\nabla\rho^{\frac{1}{p-1}}\right)\]
($V=0$, $c(z):=\frac{|\,z\,|^q}{q}$ with $\frac{1}{p}+\frac{1}{q}$, and $F(x)=\frac{1}{p-1}x\ln x$ with $p>1$).

\item Parabolic $p$-Laplacian equation:
\[\frac{\partial\rho}{\partial t} = \mbox{div}\,\left(|\,\nabla\rho\,|^{p-2}\nabla\rho\right)\]
($V=0$, $c(z)=\frac{|\,z\,|^q}{q}$ with $\frac{1}{p}+\frac{1}{q}$ and $F(x)=\frac{x^m}{m(m-1)}$ with $m:=\frac{2p-3}{p-1}$ and $p\geq \frac{2d+1}{d+1}$).

\item Doubly degenerate diffusion equation (see \cite{otto:doubly}):
\begin{equation}
\label{doubly}
  \frac{\partial\rho}{\partial t} = \mbox{div}\,\left(|\,\nabla\rho^n\,|^{p-2}\nabla\rho^n\right)
\end{equation}
($V=0$, $c(z):=\frac{|\,z\,|^q}{q}$ with $\frac{1}{p}+\frac{1}{q}$, and $F(x)=\frac{nx^m}{m(m-1)}$ with $m:=n+\frac{p-2}{p-1}$ and $\frac{1}{p-1}\neq n\geq \frac{d-(p-1)}{d(p-1)}$).
\end{itemize}
The above restrictions on $m$, $n$ and $p$ are made so that $F$ satisfies the assumptions (HF\ref{f1}) and (HF\ref{f2}) below.

We are interested in the following questions: under what conditions does (\ref{intro8}) have solutions? Is the solution unique? What are the most relevant conditions on $c,F$ and $V$, which ensure that solutions converge asymptotically to an equilibrium?

In this work, we answered the first and the second questions. We proved existence and uniqueness of weak solutions to (\ref{intro8}), when the initial mass density $\rho_0$ is bounded below and above, that is, $\rho_0+\frac{1}{\rho_0}\in L^\infty(\Omega)$ (see Theorems \ref{theoexistbound} and \ref{theoexistbound1}). This restriction was made to simplify the proofs, and not to bury fundamental facts into technical computations. We include in Remark \ref{remextension}, a method which may be used to extend our existence result to the cases where $\frac{1}{\rho_0}$ fails to be bounded, and where $\rho_0$ belong to a wider class of probability densities $\rho_0\in L^p(\Omega),\,p\geq q$. In a coming paper, we establish large time asymptotic results for solutions of (\ref{intro8}).

Our approach in studying existence of solutions to (\ref{intro8}) was inspired by the works of Jordan-Kinderlehrer-Otto \cite{otto:variation} and Otto \cite{otto:doubly}. In \cite{otto:variation}, the authors observed that the Fokker-Planck equation can be interpreted  as the gradient flow of the entropy functional
\[ H(\rho):=\int_{\R^d} \left(\rho\ln \rho +\rho\,V\right)\,\mbox{d}x,\]
with respect to the Wasserstein metric $d_2.\,$  Recall that $\,d_2\,$ is a metric on the set of probability measures on $\R^d$, with finite second moments, defined by 
\[ d_2(\mu_0,\mu_1):= \left[\,\inf\Big\{ \int_{\R^d\times\R^d} \frac{|x-y|^2}{2}\,\mbox{d}\gamma (x,y):\quad \gamma \in \Gamma(\mu_0,\mu_1)\Big\}\,\right]^{1/2}, \]
where $\Gamma(\mu_0,\mu_1)$ denotes the set of probability measures on $\R^d\times\R^d$, having $\mu_0$ and $\mu_1$ as their marginals (see the definition below). This idea was generalized by Otto in \cite{otto:doubly} for doubly degenerate diffusion equations of the form (\ref{doubly}).

Let us outline here the proof of our existence theorem to (\ref{intro8}). For the sake of illustration, we assume that $V=0$. The proof consists of four main steps.\\
{\bf Step 1}.\hspace{1mm}We interpret (\ref{intro8}) as a ``steepest descent'' of the internal energy functional
\[ {\cal P}_a(\Omega)\ni \rho\mapsto E_i(\rho):=\int_\Omega F\left(\rho(x)\right)\,\mbox{d}x \]
against the Monge-Kantorovich work $W^h_c$, where $h>0$ is a time-step size, and ${\cal P}_a(\Omega)$ denotes the set of probability density functions $\rho:\Omega\rightarrow [0,\infty)$. In other words, given a mass density  $\rho^h_{k-1}$ of the fluid at time $t_{k-1}=(k-1)h$, we define the mass density $\rho^h_k$ at time $t_k=kh$, to be the unique minimizer of the variational problem
\begin{equation}
\label{intro10.0}
(P^h_k):\;\;\inf_{\rho\in{\cal P}_a(\Omega)}\Big\{ hW_c^h\left(\rho^h_{k-1},s\right)+ E_i(\rho)\,\Big\}
\end{equation}
(see Proposition \ref{propconvana}). So, at each time $t$, the system tends to decrease its internal energy $E_i(\rho)$, while trying to minimize the work to move from state $\rho(t)$ to state $\rho(t+h)$.\\
{\bf Step 2.}\hspace{1mm} We write the Euler-Lagrange equation of $(P^h_k)$, and then, deduce that
\begin{equation}
\label{intro10}
\frac{\rho^h_k -\rho^h_{k-1}}{h} = \mbox{div}\,\Big\{\rho^h_k\nabla c^\star\left[\,\nabla\left(F^\prime(\rho^h_k)\right)\right]\Big\} + A_k(h),
\end{equation}
weakly, for $k\in\N$ (Proposition \ref{propproperty}), where $A_k(h)$ tends to $0$, as $h$ goes to $0$. (\ref{intro10}) shows clearly why (\ref{intro10.0}) is a discretization of (\ref{intro8}). \\
{\bf Step 3.}\hspace{1mm} We define the approximate solution $\rho^h$ to (\ref{intro8}), as 
\[ \left\{\begin{array}{lcl}
\rho^h(t,x) &=& \rho^h_k(x)\;\;\mbox{if}\;\; t\in ((k-1)h,kh], \; k\in \N\\ \\
\rho^h(0,x) &=& \rho_0(x),
\end{array}\right. \]
and we deduce from (\ref{intro10}) that, $\rho^h$ satisfies
\begin{equation}
\label{intro11}
\left\{ \begin{array}{l}
\frac{\partial \rho^h}{\partial t} = \mbox{div}\,\Big\{\,\rho^h\,\nabla c^\star\left[\,\nabla \left(F^\prime(\rho^h)\right)\,\right]\,\Big\} +A(h) \quad \mbox{on}\quad (0,\infty)\times\Omega\\ \\
\rho^h(t=0) = \rho_0 \quad \mbox{on}\quad \Omega
\end{array}\right.
\end{equation}
in a weak sense (Proposition \ref{propinterpolation}), where $A(h)$ is shown to be $0\left(h^{\epsilon(q)}\right),\,\epsilon(q):=\min(1,q-1)$ (Proposition \ref{proplimsecond}).\\
{\bf Step 4.}\hspace{1mm} We let $h$ go to $0$ in (\ref{intro11}), and show that $(\rho^h)_h$ converges to a function $\rho$, which solves (\ref{intro8}) in a weak sense. Here, two convergence results are established: the weak convergence of $(\rho^h)_h$ to $\rho$ in $L^1\left((0,T)\times\Omega\right), \,0<T<\infty$, for a subsequence, which proves that $\left(\frac{\partial \rho^h}{\partial t}\right)_h$ converges weakly to $\frac{\partial \rho}{\partial t}$ in $\left[C^\infty_c(\R\times\R^d)\right]^\prime$, and the weak convergence of the nonlinear term $\left(\mbox{div}\{\rho^h\nabla c^\star\left[\nabla\left(F^\prime(\rho^h)\right)\right]\}\right)_h$ to $\mbox{div}\{\rho\nabla c^\star\left[\\\nabla\left(F^\prime(\rho)\right)\right]\}$ in $\left[C^\infty_c(\R\times\R^d)\right]'$, for a subsequence.\\
The first convergence follows from the second bound in (\ref{bba}) of Proposition \ref{propconvana}, a consequence of the maximum principle stated in Proposition \ref{propbound} (see Lemma \ref{lemstrongconv1}): {\it starting with a probability density function $\rho_0$ which is bounded above, that is, $\rho_0\leq N$ a.e., the probability density function $\rho^h_k$ - solution of $(P^h_k)$ - is bounded  above, as well, that is, $\,\rho^h_k\leq N$ a.e., for $k\in \N$}. As a consequence, $(\rho^h)_h$ is bounded in $L^\infty\left((0,\infty)\times\Omega\right)$, and then, converges to some $\rho$ in $L^1\left((0,T)\times\infty\right),\, 0<T<\infty$, for a subsequence.\\
The second convergence is one of the most difficult tasks in the proof of the existence theorem. Its proof requires elaborated intermediate results. Here, we see some technical differences with the works in \cite{luckhaus:quasilinear} and \cite{otto:doubly}. Indeed, due to the weaker condition (\ref{intro5}) imposed on $c$ here, compared to the stronger ellipticity condition (\ref{intro2}) or (\ref{intro3}) in \cite{luckhaus:quasilinear}, the method used in \cite{luckhaus:quasilinear} and \cite{otto:doubly} do not yield, here, strong convergence of the nonlinear term, as in \cite{luckhaus:quasilinear} and \cite{otto:doubly}. So, to prove the -- weak -- convergence of the nonlinear term, here, we proceed as follows:
\begin{enumerate}
\item First, we improve the previous convergence, by showing that, in fact, $(\rho^h)_h$ converges strongly to $\rho$, for a subsequence, in $L^1\left((0,T)\times\Omega\right)$ (Proposition \ref{propstrongconv}).
\item Then, we deduce that $\left(\mbox{div}\{\rho^h\nabla c^\star\left[\nabla\left(F^\prime(\rho^h)\right)\right]\}\right)_h$ converges weakly to \\
$\mbox{div}\{\rho\nabla c^\star\left[\,\nabla\left(F^\prime(\rho)\right)\,\right]\}$ in $\left[C^\infty_c(\R\times\R^d)\right]'$, for a subsequence (Theorem \ref{theoweakconv}).
\end{enumerate}
To prove (\roman{mycounter}), one needs to have a good control on the spatial derivative of $\rho^h$, for example, to show that $\Big\{\,\nabla\left(F^\prime(\rho^h)\right)\,\Big\}_h$ is bounded in $L^{q^\star}\left((0,T)\times\Omega\right),\,0<T<\infty$. The main ingredient used to establish this result is the following {\it Monge-Kantorovich type energy inequality}:
\begin{equation}
\label{introen}
E_i(\tilde{\rho_0}) - E_i(\tilde{\rho_1}) \geq \int_\Omega \langle \nabla\left(F^\prime(\tilde{\rho_1})\right), \tilde{S}(y)-y\rangle \tilde{\rho_1}(y)\,\mbox{d}y,
\end{equation}
 for $\tilde{\rho_0},\,\tilde{\rho_1}\in {\cal P}_a(\Omega)$. Here, $\tilde{S}$ denotes the $c$-optimal map that pushes $\tilde{\rho_1}$ forward to $\tilde{\rho_0}$ (see the definition in Proposition \ref{propintro1}). A more general statement of the energy inequality is given in Theorem \ref{theoenergyrel}. (\ref{introen}) can be seen as a consequence of the {\it displacement convexity} of the internal energy functional ${\cal P}_a(\Omega) \ni \rho \mapsto E_i(\rho)$, that is, the convexity of
\[ [0,1]\ni t\mapsto E_i(\tilde{\rho}_{1-t}),\]
where,
\begin{equation}
\label{rhot}
\tilde{\rho}_{1-t}:= \left((1-t)\,\mbox{id}\,+t\tilde{S}\right)_{\#}\tilde{\rho_1}
\end{equation}
is the shortest path joining $\tilde{\rho_1}$ and $\tilde{\rho_0}\,$ in ${\cal P}_a(\Omega)$. When $c(z)=\frac{|\,z\,|^2}{2}$, in which case $\tilde{S}$ is the gradient of a convex function, the interpolation in (\ref{rhot}) was introduced by McCann in \cite{mccann:convexity}.

Indeed, setting $\tilde{\rho_0}:=\rho^h_{k-1}$ and $\tilde{\rho_1}:=\rho^h_k\,$ in (\ref{introen}), and using the Euler-Lagrange equation of $(P^h_k)$, that is,
\begin{equation}
\label{introeul}
\frac{S^h_k-\mbox{id}}{h}=\nabla c^\star\left[\,\nabla\left(F^\prime(\rho^h_k)\right)\,\right], 
\end{equation}
where $S^h_k$ is the $c\left(\frac{\cdot}{h}\right)$-optimal map that pushes $\rho^h_k$ forward to $\rho^h_{k-1}$, we obtain that 
\begin{equation}
\label{intro11.0}
h\int_\Omega \langle \nabla\left(F^\prime(\rho^h_k)\right), \nabla c^\star\left[\,\nabla\left(F^\prime(\rho^h_k)\right)\,\right]\rangle\, \rho^h_k \leq E_i(\rho^h_{k-1}) - E_i(\rho^h_k).
\end{equation}
We integrate (\ref{intro11.0}) over $t\in [0,T]$, and we use Jensen's inequality, to deduce that 
\begin{equation}
\label{intro11.1}
\int_0^T\int_\Omega \langle \nabla\left(F^\prime(\rho^h)\right), \nabla c^\star\left[\nabla\left(F^\prime(\rho^h)\right)\,\right]\rangle \rho^h 
\leq E_i(\rho_0)-|\,\Omega\,|F\left(\frac{1}{|\,\Omega\,|}\right).
\end{equation}
We use condition (\ref{intro5}) -- precisely, $c(z)\leq \alpha\left(|\,z\,|^q+1\right)$ --, combined with (\ref{intro11.1}) and the fact that $(\rho^h)_h$ is bounded in $L^\infty\left((0,\infty)\times\Omega\right)$, to conclude that
\[\int_0^T\int_\Omega \rho^h\,\Big|\,\nabla\left(F^\prime(\rho^h)\right)\,\Big|^{q^\star} \leq \mbox{cst}\quad \mbox{(see Lemma \ref{lemstrongconv1}}).\]
Then, we use that $\left(\frac{1}{\rho^h}\right)_h$ is bounded in $L^\infty\left((0,\infty)\times\Omega\right)$ (see (\ref{bba})) -- a consequence of the minimum principle of Proposition \ref{propbound} -- to deduce that $\Big\{\nabla\left(F^\prime(\rho^h)\right)\Big\}_h$ is bounded in $L^{q^\star}\left((0,T)\times\Omega\right),\,0<T<\infty$. This yields (i).\\
To prove (ii), we first use (\ref{introeul}) and condition (\ref{intro5}) -- precisely $c(z)\geq \beta|\,z\,|^q$ -- to have that $\Big\{\nabla c^\star\left[\nabla\left(F^\prime(\rho^h)\right)\right]\Big\}_h$ is bounded in $L^q\left(\Omega\times(0,\infty)\right)$ (Lemma \ref{lemweakconv1}), from which we deduce that  $\Big\{\nabla c^\star\left[\nabla\left(F^\prime(\rho^h)\right)\right]\Big\}_h$ converges weakly to some $\sigma$ in $L^q\left(\Omega\times(0,T)\right)$, for a subsequence, and for all $0<T<\infty$. Next, we use (i) and the boundedness of  $\Big\{\nabla\left(F^\prime(\rho^h)\right)\Big\}_h$ in $L^{q^\star}\left(\Omega\times(0,\infty)\right)$, to obtain that $\Big\{\nabla\left(F^\prime(\rho^h)\right)\Big\}_h$ converges weakly to $\nabla\left(F'(\rho)\right)$ in  $L^{q^\star}\left(\Omega\times(0,T)\right)$ for a subsequence (Lemma \ref{lemweakconv1}). In the end, we extend the energy inequality (\ref{introen}) in time-space (Lemma \ref{lemweakconv2}), and we combine the new inequality with the strong convergence of $(\rho^h)_h$ to $\rho$, the weak convergence of $\Big\{\nabla\left(F^\prime(\rho^h)\right)\Big\}_h$ to $\nabla\left(F'(\rho)\right)$, and the weak convergence of $\Big\{\nabla c^\star\left[\nabla\left(F^\prime(\rho^h)\right)\right]\Big\}_h$ to $\sigma$, to establish that $\left(\mbox{div}\{\rho^h\nabla c^\star\left[\nabla\left(F^\prime(\rho^h)\right)\right]\}\right)_h$ converges weakly to $\mbox{div}(\rho\sigma)$ for a subsequence, and that $\mbox{div}(\rho\sigma)=\mbox{div}\{\rho\nabla c^\star\left[\nabla\left(F^\prime(\rho)\right)\right]\}$ (Theorem \ref{theoweakconv}). The convexity of $c^\star$ plays an important role in this proof.\\

\noindent{\bf \large Notations}

\begin{itemize}
\item $\Omega\subset\R^d,\, d\geq 1$ is open, bounded, convex and smooth, and $\Omega_T:=(0,T)\times\Omega$, for $0<T\leq\infty$.
\item $B_R(x)\subset \R^d$ denotes the open ball of radius $R$, centered at $x$, $B_R(x)^c:=\R^d \setminus B_R(x)$, and $p^\star$ denotes the conjugate index of $p>0$, that is, $\frac{1}{p}+\frac{1}{p^\star}=1$.
\item ${\cal P}_a(\Omega):=\Big\{\rho:\Omega \rightarrow [0,\infty)\;\mbox{measurable},\; \int_\Omega \rho(x)\,\mbox{d}x=1\Big\}$, and for $0<R<\infty,\;\; {\cal P}_a^{(R)}(\Omega):=\Big\{\rho\in{\cal P}_a(\Omega):\;\rho\leq R\;\mbox{a.e.}\,\Big\}$.
\item If $\varphi: \Omega\rightarrow \R$, then  $\|\,\varphi\,\|_{L^q(\Omega)}$ denotes the $L^q$-norm of $\varphi$, and $\mbox{spt}\,(\varphi)$ denotes the support of $\varphi$, that is, the closure of $\{x\in\Omega:\;\varphi(x)\not = 0\}$.
\item If $x=(x_1,\cdots, x_d)$ and $y=(y_1,\cdots, y_d)$ are vectors in $\R^d$, then $\left\langle x,y\right\rangle :=\sum_{i=1}^d x_i y_i$, and $\,|\,x\,|:=\sqrt{\langle x,x\rangle}$.
\item If $A\subset \R^D,\,D\geq 1$, is convex, and $G:A\rightarrow \R$ is convex, then $G^\star:\R^D\rightarrow\R$ denotes the Legendre transform of $G$, that is, 
\[G^\star(y) := \sup_{x\in \R^d} \{\langle x,y\rangle - \bar{G}(x) \},\;\mbox{where}\; \bar{G}(x):=\left\{\begin{array}{lr}
G(x) &\mbox{if}\;x\in A\\
+\infty &\mbox{otherwise.}
\end{array}\right.\]
By abuse of notations, we will identify $G$ and $\bar{G}$.
\item If $A\subset\R^d$ is Borel, then $|\,A\,|$ denotes the Lebesgue measure of $A$, and $\1_A(x):=\left\{ \begin{array}{ll}
1 & \mbox{if}\; x\in A\\
0 & \mbox{otherwise,}
\end{array} \right.$ denotes the characteristic function of $A$.
\end{itemize}
Throughout this manuscript, $M$ and $N$ are positive reals, a.e. refers to the d-dimensional Lebesgue measure, and 
\[c_h(z):=c\left(\frac{z}{h}\right),\, z\in\R^d,\quad\mbox{and}\quad E_i(\rho):=\int_\Omega F\left(\rho(x)\right)\,\mbox{d}x,\,\rho\in{\cal P}_a(\Omega).\]

\noindent{\bf \large Definitions}\\

\noindent{\bf Probability measures with marginals.} Let $\mu_0$ and $\mu_1$ be probability measures on $\R^d$. A Borel probability measure $\gamma$ on the product space $\R^d\times\R^d$ is said to have $\mu_0$ and $\mu_1$ as its marginals, if one of the following equivalent conditions holds:
\begin{enumerate}
\item for Borel $A\subset \R^d$,
\[ \gamma [A\times\R^d] = \mu_0[\,A\,] \;\;\mbox{and}\;\;  \gamma [\R^d\times A] = \mu_1[\,A\,]. \]
\item For $(\varphi, \psi)\in L^1_{\mu_0}(\R^d)\times L^1_{\mu_1}(\R^d)$, where $L^1_{\mu_i}(\R^d)$ denotes the space of $\mu_i$-integrable functions on $\R^d\; (i=1,2)$,
\[ \int_{\R^d\times\R^d} \left[\,\varphi(x)+\psi(y)\,\right]\,\mbox{d}\gamma(x,y) \,=\, \int_{\R^d} \varphi(x)\,\mbox{d}\mu_0(x) + \int_{\R^d} \psi(y)\,\mbox{d}\mu_1(y).\]
\end{enumerate}
We denote by $\Gamma(\mu_0,\mu_1)$, the set of all probability measures satisfying (i) or (ii). If $\mu_0$ and $\mu_1$ are absolutely continuous with respect to Lebesgue, and $\rho_0,\,\rho_1$ denote their respective density functions, we simply write $\Gamma(\rho_0,\rho_1)$.\\

\noindent{\bf Push-forward mapping.} Let $\mu_0$ and $\mu_1$ be probability measures on $\R^d$. A Borel map $T:\,\R^d\rightarrow \R^d$ is said to push $\mu_0$ forward to $\mu_1$, if 
\begin{enumerate}
\item $\;\mu_1[\,A\,] = \mu_0[\,T^{-1}(A)\,]$ for Borel $A\subset \R^d$, or equivalently
\item $\; \int_{\R^d} \varphi(y)\,\mbox{d}\mu_1(y) = \int_{\R^d} \varphi\left(T(x)\right)\,\mbox{d}\mu_0(x)$ for $\varphi \in L^1_{\mu_1}(\R^d)$.
\end{enumerate}
Whenever (i) or (ii) holds, we write that $\mu_1 = T_{\#}\mu_0$, and we say that $T$ pushes $\mu_0$ forward to $\mu_1$.\\ \\
The next proposition is due to Caffarelli \cite{caffarelli:allocation}, and Gangbo-McCann \cite{gangbo:optimal}. It asserts the existence and uniqueness of the minimizer for the Monge-Kantorovich problem.
\begin{prop}
\label{propintro1}
({\it Existence of optimal maps}).\\
Let $c:\R^d\rightarrow [0,\infty)$ be strictly convex, and $\rho_0, \rho_1\in {\cal P}_a(\Omega)$. Then,
\begin{enumerate}
\item there is a function $v:\overline{\Omega}\rightarrow \R$ such that, $T:\bar{\Omega}\rightarrow \bar{\Omega},\,T:= \mbox{id}-(\nabla c^\star)\circ\nabla u$ pushes $\rho_0$ forward to $\rho_1$, where $u(x) = \inf_{y\in \overline{\Omega}} \Big\{ c(x-y)-v(y)\Big\}$ for $x\in\overline{\Omega}$.
\item $T$ is the unique minimizer (a.e. with respect to $\rho_0$) of the Monge problem
\[ ({\cal M}):\quad \inf\,\Big\{ \int_\Omega c\left(x-T(x)\right)\,\rho_0(x)\,\mbox{d}x, \quad T_{\#}\rho_0 = \rho_1 \Big\}.\]
\item The joint measure $\gamma := (id\times T)_{\#}\rho_0$  uniquely solves the Kantorovich problem
\[ ({\cal K}): \quad \inf\,\Big\{ \int_{\R^d\times\R^d} c\left(x-y\right)\,\mbox{d}\gamma(x,y), \quad \gamma\in \Gamma(\rho_0,\rho_1) \Big\}. \] 
\item $T$ is one-to-one, that is, there exits a map $S:\bar{\Omega}\rightarrow \bar{\Omega}$ pushing $\rho_1$ forward to $\rho_0$, such that $T\left(S(y)\right)=y$ a.e. with respect to $\rho_1$, while $S\left(T(x)\right)=x$ a.e. with respect to $\rho_0$.\\
Moreover, $S=id\,+\,\nabla c^\star(-\nabla v)$, where $v(y) = \inf_{x\in \overline{\Omega}} \Big\{ c(x-y)-u(x)\Big\}$ for $y\in\overline{\Omega}$.\\
$v$ is called the $c$-transform of $u$, and it is denoted by $v := u^c$.
\end{enumerate}
\end{prop}  

We will refer to $T$ (respectively $S$) as the $c$-optimal map that pushes $\rho_0$ (respectively $\rho_1$) forward to $\rho_1$ (respectively $\rho_0$), and $\gamma$ will be called the $c$-optimal measure in $\Gamma(\rho_0,\rho_1)$.\\

\noindent{\bf Wasserstein metric.} Let $c:\R^d\rightarrow [0,\infty)$ be strictly convex, $h>0$, and $\rho_0,\,\rho_1\in {\cal P}_a(\Omega)$. We define 
\[ W_c^h\,(\rho_0,\rho_1) \,:=\,\inf\,\Big\{\int_{\R^d\times\R^d} c\left(\frac{x-y}{h}\right)\,\mbox{d}\gamma(x,y): \quad \gamma\in\Gamma(\rho_0,\rho_1)\Big\}.\]
If $c(z)=\frac{|z|^q}{q}$, we denote $W^h_c$ by $W^h_q$. When $c(z)=\frac{|z|^2}{2}$ and $h=1,\; d_2:=\sqrt{W^h_2}$ is called the Wasserstein metric.\\

We deduce from Proposition \ref{propintro1} that, there exist a unique probability measure $\gamma\in \Gamma(\rho_0,\rho_1)$, and a unique mapping $T$ that pushes $\rho_0$ forward to $\rho_1$, whose inverse $S$ pushes $\rho_1$ forward to $\rho_0$, such that
\begin{eqnarray*}
W^h_c(\rho_0,\rho_1) = \int_{\R^d\times\R^d} c\left(\frac{x-y}{h}\right)\,\mbox{d}\gamma(x,y) &=& \int_\Omega c\left(\frac{x-T(x)}{h}\right)\,\rho_0(x)\,\mbox{d}x\\
&=& \int_\Omega c\left(\frac{S(y)-y}{h}\right)\,\rho_1(y)\,\mbox{d}y.
\end{eqnarray*}

\noindent{\bf \large Assumptions}\\

\newcounter{line}
\begin{list}{(HC\arabic{line}) : }{\usecounter{line}}
\item $c:\R^d\rightarrow [\,0,\infty)$ is such that $\,0=c(0)<c(z)\,$ for $z\neq 0.$ \label{c2}\\
\item $\lim_{|x|\rightarrow \infty} \frac{c(x)}{|x|} = \infty,\,$ i.e. $c$ is coercive. \label{c3}\\
\item $\beta\,|\,z\,|^q \,\leq\, c(z) \,\leq\, \alpha\,(|\,z\,|^q + 1)\,$ for $\,z\in \R^d,\,$ where $\alpha,\,\beta >0\;\mbox{and}\; q>1$. \label{c4}\\
\end{list}
\begin{list}{(HF\arabic{line}) : }{\usecounter{line}}
\item Either $\lim_{x\rightarrow +\infty} \frac{F(x)}{x} = +\infty$, (i.e. $F:[0,\infty)\rightarrow \R\,$ has a super-linear growth at $+\infty)$, or $\lim_{x\rightarrow +\infty} \frac{F(x)}{x} = 0$ and $F^\prime(x)<0,\;\forall x>0$. \label{f1}\\
\item $(0,\infty)\ni x\mapsto x^d F(x^{-d})\,$ is convex. \label{f2}\\
\end{list}

\noindent We impose assumption (HF\ref{f1}) to ensure that the Legendre transform $F^\star$ of $F$ is finite on $(-\infty,0)$. In fact, (HF\ref{f1}), the strict convexity of $F$, and $F(0)=0$, imply that $F^\star(x)\in\R$ for $x<0$. Then, $F^\star\left(F^\prime(a)\right)$ is finite for $a>0$.\\
The following cost and energy density functions satisfy the above assumptions:
\begin{itemize}
\item  $\,c(z)=\sum_{i=1}^n A_i\,|\,z\,|^{q_i}$, where $n\in\N,\;q_i>1$, and $A_i>0$, (take  $q=\max_{\{i=1,\cdots, n\}}(q_i) = q_{i_0}$ for some $i_0\in \{1\cdots, n\},\; \beta = A_{i_0}$, and $\alpha =\sum_{i=1}^n A_i\,$).
\item $\,F(x)=x\ln x,\; F(x)=\frac{x^m}{m-1},\, m>1$ or $1-\frac{1}{d}\leq m<1$, and $F(x)=\sum_{i=1}^n A_i\,F_i(x)$, where $n\in\N$, $A_i>0$, and the $F_i$ are like the previous $F$.
\end{itemize}

\section{Calculus of Variations on ${\cal P}_a(\Omega)$}
We discretize (\ref{intro8}), and prove in section 2.1 that 
\begin{equation}
\label{convana2}
(P):\quad \inf\Big\{I(\rho):= hW^h_c(\rho_0,\rho)+E_i(\rho):\quad \rho\in{\cal P}_a(\Omega)\Big\}\end{equation}
admits a unique minimizer $\rho_1$.
The reason why we minimize such a functional $I(\rho)$ will be clear in section 2.2, where we find the Euler-Lagrange equation of $(P)$. In fact, we shall see that the Euler-Lagrange equation is nothing but the discretization of (\ref{intro8}). In section 2.3, we show that 
\begin{equation}
\label{convana1}
 E_i(\rho_0)-E_i(\rho_1) \geq \frac{d\,E(\rho_{1-t})}{dt} \Big|_{t=0},
\end{equation}
where $\rho_{1-t}$ (\ref{rhot}) denotes the probability density obtained by interpolating $\rho_0$ and $\rho_1$ along the``geodesic'' joining them in ${\cal P}_a(\Omega)$. We refer to (\ref{convana1}) as the (internal) energy inequality. We shall see later on, that (\ref{convana1}) is an essential ingredient in the proof of the convergence of the approximate sequence $(\rho^h)_h$ (see the definition in section 2.4) to solutions of (\ref{intro8}).

\subsection{Existence of solutions to a minimization problem $(P)$}
Throughout this section, $h>0$, and $\rho_0\in{\cal P}_a(\Omega)$ is such that $\rho_0\leq M$ a.e. We show that
\begin{equation}
\label{pr}
(P_R):\;\;\inf\Big\{I(\rho):= hW^h_c(\rho_0,\rho)+E_i(\rho):\;\;\rho\in {\cal P}_a^{(R)}(\Omega)\Big\}
\end{equation}
admits a unique minimizer $\rho_{1R}$, for $R\geq M$ (Proposition \ref{propconvanabis}), and that $\rho_{1R}\in {\cal P}_a^{(M)}(\Omega)$ for $R>2M$, that is, $0\leq \rho_{1R}\leq M$ a.e. (Proposition \ref{propbound}). We deduce that $(P)$ (\ref{convana2}) has a unique minimizer $\rho_1$, which satisfies $0\leq\rho_1\leq M$ a.e. (Proposition \ref{propconvana}).
\begin{prop}
\label{propconvanabis}
Let $R\geq M$, and assume that $F:[0,\infty)\rightarrow \R$ and $c:\R^d\rightarrow [0,\infty)$ are strictly convex, and $c$ satisfies (HC\ref{c2}). Then $(P_R)$ has a unique minimizer $\rho_{1R}$, which satisfies
\begin{equation}
\label{convana3pr}
|\,\Omega\,|\,F\left(\frac{1}{|\,\Omega\,|}\right)\leq  E_i(\rho_{1R}) \leq E_i(\rho_0).
\end{equation} 
\end{prop}
{\tt Proof:} Let $I_{inf}$ denote the infimum of $I(\rho)$ over $\rho\in {\cal P}_a^{(R)}(\Omega)$. Since $\rho_0\in {\cal P}_a^{(R)}(\Omega)$,\\
$E_i(\rho_0) <\infty$ and $c(0)=0$, we have that $I_{inf}\leq \frac{1}{h} E_i(\rho_0)$. Moreover, because of Jensen's inequality and the fact that $c\geq 0$ and $\rho\in {\cal P}_a(\Omega)$, we have that $I_{inf}\geq \frac{|\,\Omega\,|}{h}F\left(\frac{1}{|\,\Omega\,|}\right)$. We deduce that  $I_{inf}$ is finite. Now, let $\left(\rho^{(n)}\right)_n$ be a minimizing sequence for $(P_R)$. We have that $\left(\rho^{(n)}\right)_n$ is bounded in $L^\infty(\Omega)$. As a consequence, $\left(\rho^{(n)}\right)_n$ converges weakly-$\star$ to a function $\rho_{1R}$ in $L^\infty(\Omega)$, and then, weakly in $L^1(\Omega)$, for a subsequence, since $\Omega$ is bounded. Clearly, $\rho_{1R}\in {\cal P}_a^{(R)}(\Omega)$. Furthermore, because of Proposition 5.3.1 \cite{agueh:thesis}, we have that ${\cal P}_a(\Omega)\ni \rho\mapsto I(\rho)$ is weakly lower semi-continuous on $L^1(\Omega)$, as the sum of weakly lower semi-continuous functions. Therefore,   
\[I(\rho_{1R}) \leq \liminf_{n\rightarrow \infty}\,I(\rho^{(n)}) = I_{inf}\leq I(\rho_{1R}),\]
which shows that $\rho_{1R}$ is a minimizer of $(P_R)$. The uniqueness of $\rho_{1R}$ follows from the convexity of ${\cal P}_a(\Omega)\ni \rho\mapsto W^h_c(\rho_0,\rho)$ and the strict-convexity of ${\cal P}_a(\Omega)\ni \rho\mapsto \int_\Omega F(\rho)\,\mbox{d}x$ (see Proposition 5.3.1 \cite{agueh:thesis}).\\
Next, we observe that $I(\rho_{1R})\leq I(\rho_0)$, and since $W^h_c(\rho_0,\rho_0)=0$ and $W^h_c(\rho_0,\rho_{1R})\geq 0$ (because of (HC\ref{c2})), we deduce that $E_i(\rho_{1R})\leq E_i(\rho_0)$. We use Jensen's inequality and the fact that $\rho_{1R}\in {\cal P}_a^{(R)}(\Omega)$, to conclude that $|\,\Omega\,|\,F\left(\frac{1}{|\,\Omega\,|}\right)\leq E_i(\rho_{1R})$ \hfill $\Box$ 
\begin{prop}
\label{propbound}
({\it Maximum/minimum principle})\\
Let $R>2M$, and $\rho_0$ be such that $N\leq \rho_0\leq M$ a.e. Assume that $F:[0,\infty)\rightarrow \R$ and $c:\R^d\rightarrow [0,\infty)$ are strictly convex, and $c$ satisfies (HC\ref{c2}). Then, the minimizer $\rho_{1R}$ for $(P_R)$ (\ref{pr}) satisfies $N\leq \rho_{1R}\leq M$ a.e. Therefore, $\rho_{1R}$ does not depend on $R$. 
\end{prop}
{\tt Proof:} The proof we present here is similar to that in \cite{otto:dynamics}, where $c(z)=\frac{|\,z\,|^2}{2}$ and $F(x)=x\ln(x)$. Since the proof of ``$\,\rho_{1R}\geq N$ a.e.'' is analogue to that of ``$\,\rho_{1R}\leq M$ a.e.'',\hspace{1.5mm} we only prove that $\,\rho_{1R}\leq M$ a.e. Suppose by contradiction that $E := \{y\in \Omega:\; \rho_{1R}(y) > M \}$ has a positive Lebesgue measure. The idea is to come up with $\rho_{1R}^{(\epsilon)}\in {\cal P}_a^{(R)}(\Omega)$, such that $I(\rho_{1R}) > I(\rho_{1R}^{(\epsilon)})$. This contradicts the fact that $\rho_{1R}$ is the minimizer of $I$ over ${\cal P}_a^{(R)}(\Omega)$.\\
Let $\gamma_R$ be the $c_h$-optimal measure in $\Gamma(\rho_0,\rho_{1R})$. We have that
\begin{equation}
\label{gammaR}
\gamma_R (E^c\times E) > 0,
\end{equation}
 where $E^c := \R^d\setminus E$; otherwise
\begin{eqnarray*}
M|\,E\,| < \int_E \rho_{1R}(y)\,\mbox{d}y &=&  \gamma_R (\R^d\times E) = \gamma_R (E\times E)\leq \gamma_R (E\times \R^d)\\
& = & \int_E \rho_0(x)\,\mbox{d}x \leq M\,|\,E\,|, 
\end{eqnarray*}
which yields a contradiction. Consider the measure $\nu := \gamma_R\, \1_{E^c\times E}\,$ defined by
\[\int_{\R^d\times\R^d} \xi(x,y)\,\mbox{d}\nu (x,y) = \int_{E^c\times E}\xi(x,y)\,\mbox{d}\gamma_R (x,y),\]
for $\xi\in C_0(\R^d\times\R^d)$, or equivalently
\[ \nu(F) = \gamma_R \left[\,F\cap (E^c\times E)\,\right],\] 
for Borel sets $F\subset \R^d\times\R^d$. Denote by $\nu_0$ and $\nu_1$ its marginals, that is,
\[\int_{\R^d\times\R^d} \left[\,\varphi(x)+\psi(y)\,\right]\,\mbox{d}\nu(x,y) = \int_{\R^d}\varphi (x)\,\mbox{d}\nu_0(x) + \int_{\R^d}\psi(y)\,\mbox{d}\nu_1(y),\]  
for $\varphi,\,\psi\in C_0(\R^d)$. Since $\nu <<\gamma_R$ and $\gamma_R\in \Gamma(\rho_0,\rho_{1R})$, we have that $\nu_0 << \rho_0(x)\,\mbox{d}x$ and $\nu_1 << \rho_{1R}(y)\,\mbox{d}y$. As a consequence, $\nu_0$ and $\nu_1$  are absolutely continuous with respect to Lebesgue. Denote by $v_0$ and $v_1$ their respective density functions. We have that,
\begin{enumerate}
\item  $0\leq v_0\leq M\;\;{\mbox a.e.},\;\;$ and $\,0\leq v_1\leq R\;\;{\mbox a.e.},\;$ and
\item $v_0 = 0\,$ a.e. on $\,E,\,$ and $\,v_1 = 0\,$ a.e. on $\,E^c$.
\end{enumerate}
For $\epsilon\in(0,1)$, we define $\rho_{1R}^{(\epsilon)}:= \rho_{1R}+\epsilon(v_0-v_1)$, and the probability measure $\gamma_R^{(\epsilon)}$ by
 \begin{eqnarray*} 
\int_{\R^d\times\R^d} \xi(x,y)\,\mbox{d}\gamma_R^{(\epsilon)}(x,y) &:=& \int_{\R^d\times\R^d} \xi(x,y)\,\mbox{d}\gamma_R(x,y)\\
 \\
& &\; + \;\epsilon\int_{E^c\times E} \left[\xi(x,x) - \xi(x,y)\right]\,\mbox{d}\gamma_R(x,y),
\end{eqnarray*}
for $\xi\in C_0(\R^d\times\R^d)$. Because of (\roman{mycounter}), (\roman{ii}) and the fact that $2M<R$,  we have that $0\leq \rho_{1R}^{(\epsilon)}\leq R$, and
\[ \int_\Omega \rho_{1R}^{(\epsilon)}(y)\,\mbox{d}y = 1 + \epsilon\left[ \gamma_R(E^c\times E) - \gamma_R(E^c\times E) \right] = 1.\]
Hence, $\rho_{1R}^{(\epsilon)}\in {\cal P}_a^{(R)}(\Omega)$. Moreover, since $\gamma_R\in \Gamma(\rho_0,\rho_{1R})$, and $\nu$ has marginals $\nu_0=v_0(x)\,\mbox{d}x$ and $\nu_1=v_1(y)\mbox{d}y$, we have that $\gamma_R^{(\epsilon)}\in \Gamma(\rho_0,\rho_{1R}^{(\epsilon)})$. Now, we show that $I(\rho_{1R}^{(\epsilon)}) < I(\rho_{1R})$, for $\epsilon$ small enough. Indeed, 
\begin{equation}
\label{convana9}
I(\rho_{1R}^{(\epsilon)}) - I(\rho_{1R}) = h\left[ W^h_c (\rho_0,\rho_{1R}^{(\epsilon)}) - W^h_c (\rho_0,\rho_{1R})\right] + \int_\Omega \left[ F(\rho_{1R}^{(\epsilon)}) - F(\rho_{1R})\right].
\end{equation}
Because $\gamma_R^{(\epsilon)}\in \Gamma(\rho_0,\rho_{1R}^{(\epsilon)})$ and $c(0)=0$, we have that
 \begin{eqnarray}
\label{convana10} 
W^h_c (\rho_0,\rho_{1R}^{(\epsilon)}) - W^h_c (\rho_0,\rho_{1R}) &\leq& \int_{\R^d\times\R^d}c\left(\frac{x-y}{h}\right)\,\mbox{d}\gamma_R^{(\epsilon)}(x,y)\nonumber\\
& & \;\;-\int_{\R^d\times\R^d} c\left(\frac{x-y}{h}\right)\,\mbox{d}\gamma_R(x,y)\nonumber \\ 
& = & -\,\epsilon\int_{E^c\times E} c\left(\frac{x-y}{h}\right)\,\mbox{d}\gamma_R(x,y).
\end{eqnarray} 
On the other hand, according to (\roman{mycounter}) and (\roman{ii}), we have, for $\epsilon$ small enough, that
\begin{equation}
\label{convana10.1er}
\rho_{1R}^{(\epsilon)} = \rho_{1R}-\epsilon v_1 \geq M - \epsilon v_1 > 0\quad \mbox{on}\quad E,
\end{equation}
and 
\begin{equation}
\label{convana10.2e}
\rho_{1R}^{(\epsilon)} = \rho_{1R}+\epsilon v_0 \geq \epsilon v_0 >0 \quad \mbox{on} \quad E^c\cap [\,v_0>0\,].
\end{equation}
We combine (\roman{mycounter}), (\roman{ii}), (\ref{convana10.1er}), (\ref{convana10.2e}), and the fact that $F\in C^1\left((0,\infty)\right)$ is convex, and $\nu=\gamma_R\1_{E^c\times E}$ has marginals $\nu_0=v_0(x)\,\mbox{d}x$ and $\nu_1=v_1(y)\,\mbox{d}y$, to obtain that
 \begin{eqnarray*}
\int_\Omega \left[ F(\rho_{1R}^{(\epsilon)}) - F(\rho_{1R})\right] &=& \int_{E^c} \left[\,F\left(\rho_{1R}+\epsilon v_0\right)-F(\rho_{1R})\,\right] + \int_E \left[\,F\left(\rho_{1R}-\epsilon v_1\right)-F(\rho_{1R})\,\right] \\ 
&\leq& \epsilon\left[ \int_{E^c\cap [\,v_0>0\,]} F^\prime(\rho_{1R}+\epsilon v_0)v_0 - \int_E F^\prime (\rho_{1R}-\epsilon v_1)v_1 \right]\\
&\leq& \epsilon\left[ \int_{E^c} F^\prime(M +\epsilon v_0)v_0 - \int_E F^\prime (M -\epsilon v_1)v_1 \right]\\
&=& \,\epsilon\,\left[\int_{E^c\times E} \left( F^\prime(M+\epsilon v_0(x)) - F^\prime(M -\epsilon v_1(y))\right)\,\mbox{d}\gamma_R(x,y)\,\right].
\end{eqnarray*}
And since $F\in C^2\left((0,\infty)\right)$, (\roman{mycounter}) and the above estimate give that
\begin{equation}
\label{convana11}
\int_\Omega \left[\,F(\rho_{1R}^{(\epsilon)} - F(\rho_{1R})\,\right] = 0(\epsilon^2).
\end{equation}
Combining (\ref{convana9}), (\ref{convana10}) and (\ref{convana11}), we conclude that, for $\epsilon$ small enough,
\[I(\rho_{1R}^{(\epsilon)}) - I(\rho_{1R}) \leq -\epsilon h\int_{E^c\times E} c\left(\frac{x-y}{h}\right)\,\mbox{d}\gamma_R(x,y) < 0,\]
where, the last inequality holds because of (HC\ref{c2}) and (\ref{gammaR}). \hfill$\Box$
\begin{prop}
\label{propconvana}
({\it Existence and uniqueness for $(P)$})\\
Assume that $N\leq \rho_0\leq M$ a.e., and $F:[0,\infty)\rightarrow \R$ and $c:\R^d\rightarrow [0,\infty)$ are strictly convex, and $c$ satisfies (HC\ref{c2}). Then, $\rho_1:=\rho_{1R}$ (defined in Proposition \ref{propbound}) is the unique minimizer for $(P)$ (\ref{convana2}). Therefore, 
\begin{equation}
\label{bba}
N\leq \rho_1\leq M \;{\mbox a.e.},
\end{equation}
and
\begin{equation}
\label{convana3}|\,\Omega\,|\,F\left(\frac{1}{|\,\Omega\,|}\right)\leq  E_i(\rho_1) \leq E_i(\rho_0).
\end{equation} 
\end{prop}
{\tt Proof:} Let $\rho\in{\cal P}_a(\Omega)$, and $\left(\rho^{(R)}\right)_{R>2M}$ be a sequence in ${\cal P}_a^{(R)}(\Omega)$ converging to $\rho$, such that,
\begin{equation}
\label{convana10.3e}
\int_\Omega F\left(\rho^{(R)}\right) \leq \int_\Omega F(\rho),
\end{equation}
as in Proposition 1.4.1 \cite{agueh:thesis}. Since $\rho_1$ is the minimizer for $(P_R)$ (Proposition \ref{propbound}), we have, using (\ref{convana10.3e}), that
\begin{equation}
\label{number2}
hW^h_c(\rho_0,\rho_1)+\int_\Omega F(\rho_1) \leq hW^h_c(\rho_0,\rho^{(R)})+\int_\Omega F(\rho).
\end{equation}
And, since $\left(\rho^{(R)}\right)_R$ converges to $\rho$ in $L^1(\Omega)$, Proposition 5.3.2 \cite{agueh:thesis} gives that 
\begin{equation}
\label{number3}
\lim_{R\uparrow \infty} W^h_c \left(\rho_0,\rho^{(R)}\right) = W^h_c (\rho_0,\rho).
\end{equation}
We let $R$ go to $\infty$ in (\ref{number2}), and we use (\ref{number3}), to conclude that $\rho_1$ is a minimizer for $(P)$. The uniqueness of the minimizer follows from the strict convexity of ${\cal P}_a(\Omega)\ni \rho\mapsto I(\rho)$ as in Proposition \ref{propconvanabis}, and (\ref{bba}) and (\ref{convana3}) are direct consequences of Proposition \ref{propbound} and (\ref{convana3pr}) \hfill$\Box$\\
In the remaining of this section, we state two propositions needed to establish the convergence of the approximate solution $\rho^h$ (\ref{inter}) of (\ref{intro8}), as $h$ goes to $0$. The first proposition stated below, shows that the interpolant densities $\rho_{1-t},\,t\in [0,1]$ (\ref{rhot}) between two probability densities $\rho_0$ and $\rho_1$, which are bounded above, are also bounded above.

\begin{prop}
\label{propintbound}
Let $\rho_0,\,\rho_1\in {\cal P}_a(\Omega)$ be such that $\rho_0,\,\rho_1\leq M\;\mbox{a.e.}$, and assume that $c:\R^d\rightarrow [0,\infty)$ is strictly convex, of class $C^1$, and satisfies $c(0)=0$ and (HC\ref{c4}). Denote by $S$ the $c$-optimal map that pushes $\rho_1$ forward to $\rho_0$, and define the interpolant map 
\[S_t:= (1-t)\,\mbox{id} + tS,\]
 for $t\in[\,0,1\,]$. Then, for $\xi\in C_c(\R^d),\;\xi\geq 0$,
\begin{equation}
\label{intbound1}
\int_\Omega \xi\,\left(S_t(y)\right)\,\rho_1(y)\,\mbox{d}y \,\leq\, M\,\int_{\R^d}\xi\,(x)\,\mbox{d}x. 
\end{equation}
\end{prop}
{\tt Proof:} The proof will be done in two steps. In step 1, we prove (\ref{intbound1}) for sufficiently regular cost functions. Here, we use the fact that $\nabla S$ is diagonalizable with positive eigenvalues when $c, c^\star \in C^2(\R^d)$, and $A\mapsto \left(\mbox{det}\right)^{1/d}$ is concave on the set of $d\times d$ diagonalizable matrices with positive eigenvalues. In step 2, we approximate a general cost function $c$ by regular cost functions $c_k$, and we obtain (\ref{intbound1}) in the limit as $k$ goes to $\infty$.\\
{\bf{Step 1}.}\hspace{1mm} {\sl $c$ is strictly convex, and $c,\,c^\star\in C^2(\R^d)$.}\\
Proposition \ref{propappendix1} gives that $\mu_{1-t}:=(S_t)_{\#}\rho_1$ is absolutely continuous with respect to Lebesgue, for $t\in[0,1]$. Let $\rho_{1-t}$ denote the density function of $\mu_{1-t}$. Then, (\ref{intbound1}) reads as
\[ \int_{\Omega} \xi\,(x)\,\rho_{1-t}\,(x)\,\mbox{d}x \,\leq\, M\,\int_{\R^d} \xi\,(x)\,\mbox{d}x.\]
Thus, it suffices to show that $\rho_{1-t}\leq M$.\\
Recall that, because of Proposition \ref{propappendix1}, there exists a set $K\subset \Omega$ of full measure for $\mu_1:=\rho_1(y)\,\mbox{d}y$, such that $S_t$ is injective on $K$, and for $y\in K$ and $t\in [0,1]$, $\nabla S(y)$ is diagonalizable with positive eigenvalues, and
\begin{equation}
\label{intbound2}
0\not = \rho_1(y) = \rho_{1-t}\left(S_t(y)\right) \,\mbox{det}\,\left[\,\nabla S_t(y)\,\right],
\end{equation}
where $\nabla S_t(y) = (1-t)\mbox{id} + t\nabla S(y)$. Since $\rho_0,\rho_1\leq M$ a.e., and $S_{\#}\rho_1=\rho_0$, we can choose $K$, such that $\rho_1(y),\rho_0\left(S(y)\right) \leq M$ for $y\in K$. We use $t=1$ in (\ref{intbound2}), and the fact that $\rho_0\left(S(y)\right)\leq M$, to deduce that 
\begin{equation}
\label{intbound3}
\mbox{det}\,\left[\,\nabla S(y)\,\right] \geq \frac{\rho_1(y)}{M}.
\end{equation}
Because $A\mapsto \left(\mbox{det}\,A\right)^{1/d}$ is concave on the set of $d\times d$ diagonalizable matrices with positives eigenvalues, we have that
\begin{equation}
\label{intbound4}
\left[\,\mbox{det}\,\nabla S_t(y)\,\right]^{1/d} \geq (1-t) + t\left(\,\mbox{det}\,\left[\nabla S(y)\right]\,\right)^{1/d}.
\end{equation}
We use (\ref{intbound3}), (\ref{intbound4}), and the fact that $\rho_1(y)\leq M$, to obtain that 
\begin{equation}
\label{intbound5}
\mbox{det}\,\left[\,\nabla S_t(y)\,\right] \geq \frac{\rho_1(y)}{M}.
\end{equation}
We combining (\ref{intbound2}), (\ref{intbound5}), and we use the fact that $S_t$ is injective on K, to deduce that $\rho_{1-t} \leq M $ on $S_t(K)$. But, since $\mu_1(K^c)=0$, and $\mu_{1-t}=(S_t)_{\#}\mu_1$, we have that $\mu_{1-t}\left[\,\left(S_t(K)\right)^c\,\right]=0$, and then $\rho_{1-t}=0$ on $\left[\,S_t(K)\,\right]^c$. We  conclude that $\rho_{1-t}\leq M$.\\
{\bf{Step 2}.}\hspace{1mm} {\sl $c$ satisfies the assumptions of the proposition. }\\
Let $(c_k)_k$ be a sequence of strictly convex cost functions satisfying
\begin{equation}
\label{exibc12}
\left\{ \begin{array}{l}
c_k, \, c_k^\star \in C^2(\R^d),\\
c_k \rightarrow c\;\;\mbox{locally in}\;\;C^1(\R^d),\;\mbox{as}\; k\rightarrow \infty,\\
0 = c_k(0)<c_k(z)\;\;\mbox{for}\;\;\,z\neq 0 \end{array}\right.
\end{equation}
(see \cite{agueh:thesis}, Proposition 1.3.1). Denote by $S_k$ the $c_k$-optimal map that pushes $\rho_1$ forward to $\rho_0$, and set
\[ S_k^{(t)}:= (1-t)\,\mbox{id}+tS_k,\]
for $t\in[0,1]$. Lemma 2.2.2 \cite{agueh:thesis} gives that $\left(S_k^{(t)}\right)_k$ converges to $S_t$ a.e. on $[\rho_1\not =0]$ for a subsequence, and because of step 1, we have that 
\begin{equation}
\label{intbound8}
\int_\Omega \xi\,\left(S_k^{(t)}(y)\right)\rho_1(y)\,\mbox{d}y \leq M\int_{\R^d} \xi(x)\,\mbox{d}x.
\end{equation}
We let $k$ go to $\infty$ in (\ref{intbound8}), and we use that $0\leq \xi\in C_c(\R^d)$, and Fatou's lemma, to conclude (\ref{intbound1}) \hfill$\Box$\\

Next, we state a proposition needed in the next section, to prove the strong convergence of the approximate solutions $(\rho^h)_h$ of (\ref{intro8}) in $L^1\left((0,T)\times\Omega\right)$, for $0<T<\infty$.
\begin{prop}
\label{propinterpolation1}
Let $f:\R\rightarrow \R$ and $g:\R\rightarrow \R$ be strictly convex, of class $C^1(\R)$, such that $\lim_{t\rightarrow \infty} \frac{g(t)}{t}=\infty$. Given $M,\,\delta>0$, and $p,\,q>1$, define\\
$A_{M,\delta}:=\Big\{(u_1,u_2)\in L^q(\Omega)^2:\;\|\,u_j\,\|_{L^q(\Omega)}\leq M,\;\|\,g^\prime(u_j)\,\|_{W^{1,p}(\Omega)}\leq M,\;\;\mbox{and}\\ \protect{\qquad\qquad} \;\;\int_\Omega\left[\,f^\prime(u_2)-f^\prime(u_1)\,\right]\left[\,u_2-u_1\,\right]\leq \delta,\,\;(j=1,2)\Big\},$\\
and set
\[\Lambda_M(\delta):=\sup_{(u_1,u_2)\in A_{M,\delta}}\|\,u_2-u_1\,\|_{L^1(\Omega)}.\]
Then
\[\lim_{\delta\downarrow 0} \Lambda_M(\delta)=0.\]
\end{prop}

\noindent{\tt Proof:} Suppose by contradiction that there exist $\kappa>0$ and $(u_j^\delta)_{\delta\downarrow 0},\,(j=1,2)$, such that $(u_1^\delta,u_2^\delta)\in A_{M,\delta}$, and
\begin{equation}
\label{interp1}
 \|\,u_2^\delta-u_1^\delta\,\|_{L^1(\Omega)} > \kappa.
\end{equation}
By the Sobolev embedding theorem, $\left(g^\prime(u_j^\delta)\right)_\delta$ converges strongly in $L^p(\Omega)$, and then, a.e., for a (non-relabeled) subsequence. Since $g\in C^1(\R)$ is strictly convex and has a super-linear growth at $\infty$, we have that $(g^\prime)^{-1}$ is continuous. We deduce that
\begin{enumerate}
\item $(u^\delta_j)_\delta$ converges to some function $u_j\,$ a.e., for $j=1,2$.
\end{enumerate}
We use (\roman{mycounter}), $\|\,u_j^\delta\,\|_{L^q(\Omega)}\leq M$, and the fact that $q>1$ to conclude that $(u_j^\delta)_\delta$ converges strongly to $u_j$ in $L^1(\Omega)$. And since $\|\,u^\delta_1-u^\delta_2\,\|_{L^1(\Omega)}>\kappa$, we obtain that
\begin{equation}
\label{interp2}
\|\,u_2-u_1\,\|_{L^1(\Omega)} > \kappa.
\end{equation}
Now, we use (\roman{mycounter}), the convexity of $f$, and the fact that $\int_\Omega\left[\,f^\prime(u^\delta_2)-f^\prime(u^\delta_1)\,\right]\left[\,u^\delta_2-u^\delta_1\,\right]\leq \delta$, to have that
\[ 0\leq \int_\Omega\left[f^\prime(u_2)-f^\prime(u_1)\right]\left[u_2-u_1\right] \leq \liminf_{\delta\downarrow 0} \int_\Omega\left[f^\prime(u_2^\delta)-f^\prime(u_1^\delta)\right]\,\left[u_2^\delta-u_1^\delta\right] \leq 0.\]
This implies that
\begin{equation}
\label{interp3}
\left[f^\prime\left(u_2(x)\right)-f^\prime\left(u_1(x)\right)\right]\left[u_2(x)-u_1(x)\right] = 0 \;\;\mbox{for a.e.}\; x\in\Omega.
\end{equation}
Since $f\in C^1(\R)$ is strictly convex, we have that $f^\prime$ is one-to-one, and then, (\ref{interp3}) implies that $u_1(x)=u_2(x)$ for a.e. $x\in \Omega$. This yields a contradiction to (\ref{interp2}) \hfill$\Box$

\subsection{Properties of the minimizer for $(P)$}
We establish the Euler-Lagrange equation for $(P)$ (\ref{convana2}), and we derive some properties of the minimizer for this problem. The next proposition is the first step towards showing that $(P)$ is a discretization of (\ref{intro8}), or in other words, (\ref{intro8}) is the steepest descent of the internal energy functional $E_i$, with respect to the Monge-Kantorovich work $W^h_c $.

\begin{prop} 
\label{propproperty}
Let $\rho_0\in {\cal P}_a(\Omega)$ be such that $N\leq \rho_0\leq M$ a.e. Assume that $F:[0,\infty)\rightarrow \R$ is strictly convex, and satisfies $F\in C^2\left((0,\infty)\right)$, and $c:\R^d\rightarrow [0,\infty)$ is strictly convex, of class $C^1$, and satisfies (HC\ref{c2}) - (HC\ref{c3}). If $\rho_1$ denotes the minimizer for $(P)$, then the followings hold:
\begin{equation}
\label{exibcN1}
\int_{\Omega \times \Omega} \langle\, \nabla c\left(\frac{x-y}{h}\right), \psi(y)\, \rangle \,\mbox{d}\gamma (x,y) + \int_\Omega P(\rho_1(y))\, div\,\psi(y)\,\mbox{d}y = 0,
\end{equation}
for $\psi \in C^\infty _c (\Omega, \R^d)$; here $P(x) := P_F(x):= x\,F^\prime(x) - F(x)$ for $x\in(0,\infty)$, and $\gamma$ is the $c_h$-optimal measure in $\Gamma(\rho_0,\rho_1)$. Moreover,
\begin{enumerate}
\item $P(\rho_1) \in W^{1,\infty} (\Omega)$.
\item If $S$ is the $c_h$-optimal map that pushes $\rho_1$ forward to $\rho_0$, then
\begin{equation}
\label{exibc3}
\frac{S(y)-y}{h} \,=\, \nabla c^\star \left[ \nabla (F^\prime (\rho_1(y)))\right], 
\end{equation}
for a.e. $y\in \Omega$, and
\begin{eqnarray}
\label{exibc4}
\Big|\,\int_\Omega \frac{\rho_1(y)-\rho_0(y)}{h}\, \varphi(y) \,\mbox{d}y &+& \int_\Omega \rho_1(y)\langle\, \nabla c^\star \left[\,\nabla \left(F^\prime(\rho_1(y))\right)\,\right], \nabla \varphi(y)\, \rangle \,\mbox{d}y\;\Big|\nonumber \\ 
&\leq&  \frac{1}{2h}\,\sup_{x\in \bar{\Omega}} \Big|\,D^2\varphi(x)\,\Big|\; \int_{\Omega\times\Omega} |\,x-y\,|^2\,\mbox{d}\gamma (x,y),  
\end{eqnarray}
for $\varphi \in C^2(\bar{\Omega})$.
\end{enumerate}
\end{prop}
{\tt Proof:} Since $c\in C^1(\R^d)$ is strictly convex and satisfies (HC\ref{c3}), we have that $c^\star\in C^1(\R^d)$, and $(\nabla c)^{-1} = \nabla c^\star$. Following \cite{otto:variation}, we consider the flow map $(\phi _{\epsilon})_{\epsilon \in \R}$  in $C^{\infty}(\Omega,\Omega)$, defined by
\begin{equation}
\label{property1}
 \left\{ \begin{array}{ll}
  \frac{\partial\phi_\epsilon}{\partial\epsilon} = \psi \circ \phi _{\epsilon}\\
  \phi _0 = \mbox{id}.
  \end{array} \right. 
\end{equation}
where $\psi\in C^\infty_c(\Omega,\R^d)$. We have that $det\,(\nabla\phi_\epsilon) \neq 0$,  and 
\begin{equation}
\label{exibc6}
\frac{\partial (det\,\nabla\phi_\epsilon)}{\partial \epsilon}\Big |_{\epsilon =0} = \mbox{div}\,\psi.
\end{equation}
We define on $\Omega$, the probability measure $\mu _\epsilon := (\phi_\epsilon)_{\#}\rho_1$. Since $\phi_\epsilon$ is a $C^1$-diffeomorphism, then $\mu_\epsilon$ is absolutely continuous with respect to Lebesgue. Let $\rho_\epsilon$ denote its density function. Clearly, $\rho_\epsilon\in{\cal P}_a(\Omega)$, and 
 \begin{equation}
 \label{exibc7}
 (\rho_\epsilon \circ \phi_\epsilon ) \, det\, (\nabla\phi_\epsilon) = \rho_1 \quad \mbox{a.e}.
  \end{equation}
Next, we define on $\Omega\times\Omega$, the probability measure $\gamma_\epsilon := (id\,\times \phi_\epsilon)_{\#}\gamma$, i.e. 
\[\int_{\Omega\times\Omega}\xi(x,y)\,\mbox{d}\gamma_\epsilon (x,y) = \int_{\Omega\times\Omega} \xi\left(x,\phi_\epsilon(y)\right)\,\mbox{d}\gamma(x,y), \qquad \forall\, \xi \in C(\Omega\times\Omega). \]
We have that $\gamma_\epsilon \in \Gamma (\rho_0,\rho_\epsilon)$, and then, the mean-value theorem gives that
\begin{eqnarray*}
\lefteqn{\frac{W^h_c (\rho_0,\rho_\epsilon) - W^h_c (\rho_0,\rho_1)}{\epsilon} } \\
& & \leq  \int \frac{1}{\epsilon}\left[\,c_h(x-\phi_\epsilon(y))-c_h(x-y)\,\right]\,\mbox{d}\gamma(x,y) \\ \\ 
& & = - \int \langle\; \nabla c_h\left[\,x-y + \theta (y-\phi _\epsilon (y))\,\right], \frac{\phi _\epsilon - \phi _0}{\epsilon}(y)\;\rangle \,\mbox{d}\gamma(x,y),
\end{eqnarray*}
where $\theta \in [0,1]$. Because of (\ref{property1}), we have that $|\,\frac{\phi_\epsilon -\phi_0}{\epsilon}\,|\leq \|\,\psi\,\|_{L^\infty}$, for $\epsilon >0$. Then, we use that $c\in C^1(\R^d)$, the Lebesgue dominated convergence theorem, and (\ref{property1}), to obtain that
\begin{equation}
\label{exibc8}
\limsup _{\epsilon \downarrow 0} \frac{W^h_c (\rho_0,\rho_\epsilon) - W^h_c (\rho_0,\rho_1)}{\epsilon} \,\leq\, - \int \langle\, \nabla c_h(x-y), \psi (y)\,\rangle \,\mbox{d}\gamma(x,y).
\end{equation}
 On the other hand, because of (\ref{exibc7}), we have that
\[ \int_\Omega F\left(\rho_\epsilon(x)\right)\,\mbox{d}x = \int _{\Omega} F(\rho_\epsilon \circ \phi_\epsilon(y))\,\mbox{det}\, \nabla\phi_\epsilon (y)\,\mbox{d}y = \int _{\Omega} F\left(\frac{\rho_1(y)}{\,\mbox{det}\,\nabla\phi _\epsilon (y)}\right)\,\mbox{det}\,\nabla\phi _\epsilon (y)\,\mbox{d}y.\]
And since, $F\in C^1\left((0,\infty)\right)$, we deduce by the mean-value theorem that
\begin{eqnarray}
\label{property2}
\lefteqn{\int _{\Omega} \frac{F(\rho_\epsilon (x))-F(\rho_1(x))}{\epsilon}\,\mbox{d}x} \nonumber\\
& & = \frac{1}{\epsilon} \int _{\Omega} \left[\,\left(F\left(\frac{\rho_1}{det\,\nabla\phi _\epsilon }\right)-F(\rho_1)\right)\, det\,\nabla\phi _\epsilon  \,+\,F(\rho_1)( det\,\nabla\phi _\epsilon - 1)\,\right] \nonumber\\
& & = \int _{\Omega} \left[\,- F^{\prime}\left(\rho_1+\theta\left( \frac{\rho_1}{det\, \nabla\phi _\epsilon }-\rho_1\right)\right)\,\rho_1 \frac{det\,\nabla\phi_\epsilon -1}{\epsilon}\,\right] \nonumber \\ 
& & \quad +\,\int_{\Omega} \left[F(\rho_1)\frac{det\,\nabla\phi_\epsilon -1}{\epsilon} \right],
\end{eqnarray} 
where $\theta\in[0,1]$. We combine (\ref{property1}), (\ref{exibc6}) and (\ref{property2}) to have that
\begin{equation}
\label{exibc9}
\lim _{\epsilon\downarrow 0}\int _{\Omega} \frac{F(\rho_\epsilon (y))-F(\rho_1(y))}{\epsilon}\,\mbox{d}y = -\int_\Omega P(\rho_1(y))\,\mbox{d}iv\,\psi(y)\,\mbox{d}y.
\end{equation}
We use (\ref{exibc8}) and (\ref{exibc9}), to conclude that 
\begin{equation}
\label{property4}
\int_{\Omega\times\Omega} \langle\,\nabla c_{h}(x-y),\psi (y)\,\rangle\,\mbox{d}\gamma(x,y) +\frac{1}{h}\,\int _{\Omega} \,P(\rho_1(y))\,\mbox{div}\,\psi(y)\,\mbox{d}y \leq 0.
\end{equation}
Since $\nabla c_h(z)=\frac{1}{h}\,\nabla c\left(\frac{z}{h}\right)$, and $\psi$ is arbitrarily chosen in $C^\infty_c(\Omega,\R^d)$, (\ref{property4}) implies (\ref{exibcN1}).\\
(i). By (\ref{bba}), $N\leq \rho_1\leq M$ a.e., and since $F\in C^1\left((0,\infty)\right)$, we have that $P(\rho_1) \in L^\infty (\Omega)$. Now, let $\varphi \in C^\infty _c (\Omega)$, and for an arbitrary  $i\in \N$, define $\psi=(\psi_j)_{j=1,\cdots,d}\in C^\infty _c(\Omega, \R^d)$ by $\psi_j:= \delta_{ij}\,\varphi$, where $\delta_{ij}$ denotes the Kronecker symbol. Because of (\ref{exibcN1}), we have that
\begin{eqnarray*}
\Big|\,\int_\Omega P(\rho_1(y))\frac{\partial \varphi}{\partial z_i}(y)\,\Big| &=& \Big|\int_{\Omega\times\Omega} \frac{\partial c}{\partial z_i}\,\left(\frac{x-y}{h}\right)\varphi (y)\,\mbox{d}\gamma (x,y)\,\Big| \\
&\leq& \sup_{x,y\in\Omega} \Big|\,\frac{\partial c}{\partial z_i}\left(\frac{x-y}{h}\right) \,\Big|\;\int_\Omega |\,\varphi(y)\,|\rho_1(y)\,\mbox{d}y\\ \\ 
&\leq&  M \;\|\,\varphi\,\|_{L^1(\Omega)}\,\sup_{x,y\in\Omega} \Big|\,\frac{\partial c}{\partial z_i} \left(\frac{x-y}{h}\right)\,\Big|.
\end{eqnarray*} 
And since $c\in C^1(\R^d)$, we deduce (i).\\
(ii). Because $P(\rho_1)\in W^{1,\infty}(\Omega)$, we can integrate by parts in (\ref{exibcN1}). We use that $\gamma\in \Gamma(\rho_0,\rho_1)$ and $S_{\#}\rho_1=\rho_0$, to obtain that
\begin{eqnarray*}
\int_\Omega \langle\,\nabla c\left(\frac{S(y)-y}{h}\right), \psi(y)\,\rangle \,\rho_1(y)\,\mbox{d}y &=& \int_\Omega \langle\,\nabla \left[\,P(\rho_1(y))\,\right], \psi(y)\,\rangle \,\mbox{d}y \\ 
&=& \int_\Omega \rho_1(y)\,\langle\,\nabla \left[F^\prime (\rho_1(y))\right], \psi(y)\,\rangle \,\mbox{d}y,
\end{eqnarray*} 
for $\psi\in C^\infty_c(\Omega,\R^d)$. And since $\psi$ is arbitrarily chosen, we deduce that
\begin{equation}
\label{exibc10}
\nabla c\left(\frac{S(y)-y}{h}\right)\,\rho_1(y) = \nabla \left[\,F^\prime(\rho_1(y))\,\right]\,\rho_1(y),
\end{equation}
for a.e. $y\in \Omega$. We combine (\ref{exibc10}), and the fact that $(\nabla c)^{-1}=\nabla c^\star$ and $\rho_1\not = 0$ a.e., to conclude (\ref{exibc3}).\\
Next, we consider $\varphi\in C^2(\bar{\Omega})$, we take the scalar product of both sides of (\ref{exibc3}) with $\rho_1(y)\nabla \varphi(y)$, and we use that $\gamma=\left(\mbox{id}\times S\right)_{\#}\rho_1$, to obtain that
\begin{equation}
\label{exibc11}
\frac{1}{h}\int_{\Omega \times \Omega} \langle\,y-x,\nabla \varphi(y)\,\rangle\, d\gamma(x,y) = - \int_\Omega \langle\,\nabla c^\star\left[\,\nabla\left(F^\prime\left(\rho_1(y)\right)+V(y)\right)\,\right], \nabla\,\varphi(y)\,\rangle\,\rho_1(y)\,\mbox{d}y. 
\end{equation}
Now, we express $\frac{1}{h}\int_{\Omega \times \Omega} \langle\,y-x,\nabla \varphi(y)\,\rangle\, d\gamma(x,y)$ in terms of $\int_\Omega \frac{\rho_1(y)-\rho_0(y)}{h}\,\varphi(y)\,\mbox{d}y$. Since $\gamma \in \Gamma(\rho_0,\rho_1)$, we have that
\[\int _{\Omega}\frac {\rho_1(y)-\rho_0(y)}{h}\,\varphi(y)\,\mbox{d}y = \frac{1}{h} \int_{\Omega \times \Omega}\,\left[\,\varphi (y)-\varphi (x)\,\right]\,\mbox{d}\gamma(x,y).\] 
Combining the above equality with the first order Taylor expansion of $\varphi$ around $y$, we obtain that
\begin{eqnarray}
\label{property5}
\lefteqn{\Big|\, \frac{1}{h} \int _{\Omega \times \Omega} \langle\,y-x,\nabla \varphi (y)\,\rangle d\gamma(x,y)-\frac{1}{h} \int_\Omega \left(\rho_1(y)- \rho_0(y)\right)\varphi(y)\;\Big|} \nonumber\\
& & \leq  \frac{1}{2h} \sup _{x\in \bar{\Omega}} |\,D^{2}\varphi (x)\,|\; \int_{\Omega \times \Omega} |\, x-y\, |^2\,\mbox{d}\gamma(x,y).
\end{eqnarray}
We substitute (\ref{exibc11}) into (\ref{property5}) to conclude  (\ref{exibc4})\hfill$\Box$

\subsection{Energy inequality}
We establish an inequality relating the internal energy $E_i(\rho_0)$ and $E_i(\rho_1)$ of two probability density functions $\rho_0$ and $\rho_1$. This inequality will be called {\it energy inequality} and will be used later on, to improve compactness properties of the approximate sequence $\rho^h$ (see the definition in section 2.4), to solutions of (\ref{intro8}). First, we prove this inequality for smooth cost functions $c$ and $c^\star$, whose Legendre transform $c^\star$ are $C^2$. Instead of using the density function $F$, we consider a more general function $G,$ which satisfies some assumptions to be specified later on. The (internal) energy inequality reads as
\begin{equation}
\label{convana12}
\int_\Omega G(\rho_0(y))\,\mbox{d}y - \int_\Omega G(\rho_1(y))\,\mbox{d}y \geq -\int_\Omega P_G(\rho_1(y))\,\mbox{div}\,(S(y)-y)\,\mbox{d}y, 
\end{equation}
where $S$ is the $c$-optimal map that pushes $\rho_1$ forward to $\rho_0$, and $P_G(x):=xG^\prime(x)-G(x)$. For smooth cost functions $c$, this inequality is simply a consequence of the displacement convexity of ${\cal P}_a(\Omega)\ni \rho\mapsto \int_\Omega G(\rho(x))\,\mbox{d}x$, that is, the convexity of $[0,1]\ni t\mapsto \int_\Omega G(\rho_{1-t}(x))\,\mbox{d}x$, where $\rho_{1-t}$ is the probability density obtained by interpolating $\rho_0$ and $\rho_1$ along the ``geodesic'' joining them in ${\cal P}_a(\Omega)$ (see Proposition \ref{propappendix1}). To prove (\ref{convana12}), we rather follow a more direct procedure, using the following result of  Cordero and Otto (Proposition \ref{propappendix1}): {\sl if $\rho_0,\,\rho_1\in{\cal P}_a(\Omega),\; c,\,c^\star \in C^2(\R^d)$, and $S$ is the $c$-optimal map that pushes $\rho_1$ forward to $\rho_0$, then $\nabla S(y)$ is diagonalizable with positive eigenvalues for $\mu_1:=\rho_1(y)\,\mbox{d}y$ - a.e. $y\in \Omega$. More!
over, the pointwise Jacobian $\mbox{det}\,\nabla S$, satisfies
\begin{equation}
\label{energy1}
0\not = \rho_1(y) = \mbox{det}\,\nabla S(y)\,\rho_0(S(y)),
\end{equation}
for $\mu_1$ - a.e. $y\in \Omega.\,$}

\begin{prop}
\label{propenergy}
({\it Energy inequality for regular cost functions})\\
Let $\rho_0,\,\rho_1\in{\cal P}_a(\Omega)$ be density functions of two Borel probability measures $\mu_0$ and $\mu_1$ on $\R^d$, respectively.  Let $\bar{c}:\R^d\rightarrow [0,\infty)$ be strictly convex, such that $\bar{c},\,\bar{c}^\star\in C^2(\R^d)$. Let $G:[0,\infty)\rightarrow \R$ be differentiable on $(0,\infty)$, such that $G(0)=0$, and $(0,\infty)\ni x\mapsto x^dG(x^{-d})$ be convex and nonincreasing. Then, the internal energy inequality (\ref{convana12}) holds.
In addition, if $P_G(\rho_1)\in W^{1,\infty}(\Omega)$ and $\rho_1>0$ a.e., then
\begin{equation}
\label{energy0bis}
\int_\Omega G(\rho_0(y))\,\mbox{d}y - \int_\Omega G(\rho_1(y))\,\mbox{d}y \geq \int_\Omega \langle\nabla [\,G^\prime\left(\rho_1(y)\right)\,], S(y)-y\rangle \rho_1(y)\,\mbox{d}y. 
\end{equation}
\end{prop}
{\tt Proof:} Set 
\[ A(x):= x^dG(x^{-d}),\;x\in(0,\infty).\]
We observe that
\begin{equation}
\label{energy0bis2}
A^\prime(x)=-dx^{d-1}P_G(x^{-d}).
\end{equation}
Since $A$ is nonincreasing, we have that $P_G\geq 0$, and then
\begin{enumerate}
\item $(0,\infty) \ni x\mapsto \frac{G(x)}{x}$ is nondecreasing.
\end{enumerate}
\noindent Proposition \ref{propappendix1} gives that $\nabla S(y)$ is diagonalizable with positive eigenvalues, and that (\ref{energy1}) holds for $\mu_1$ - a.e. $y\in \Omega$. So, $\rho_0(S(y))\not = 0$ for $\mu_1$ - a.e. $y\in \Omega$. We use that $G(0)=0,\; S_{\#}\rho_1 = \rho_0$, and (\ref{energy1}), to deduce that
\begin{eqnarray}
\label{energy2}
\int_\Omega G(\rho_0(x))\,\mbox{d}x &=& \int_{[\,\rho_0\not = 0\,]} \frac{G\left(\rho_0(x)\right)}{\rho_0(x)}\rho_0(x)\,\mbox{d}x = \int_\Omega \frac{G(\rho_0(S(y)))}{\rho_0(S(y))}\rho_1(y)\,\mbox{d}y \nonumber\\
&=& \int_\Omega G\left(\frac{\rho_1(y)}{det\,\nabla S(y)}\right)\,\mbox{d}et\,\nabla S(y)\,\mbox{d}y.
\end{eqnarray} 
Comparing the geometric mean $\left(\mbox{det}\,\nabla S(y)\right)^{1/d}$ to the arithmetic mean $\frac{\mbox{tr}\,\nabla S(y)}{d}$, we have that
\[ \frac{\rho_1(y)}{\mbox{det}\,\nabla S(y)} \geq \rho_1(y)\left(\frac{d}{\mbox{tr}\,\nabla S(y)}\right)^d.\]
Then, we deduce from (\roman{mycounter}) and the above inequality, that
\begin{equation}
\label{energy3}
 G\left(\frac{\rho_1(y)}{det\,\nabla S(y)}\right)\,\mbox{d}et\,\nabla S(y) \geq \Lambda^d\,G\left(\frac{\rho_1(y)}{\Lambda^d}\right) = \rho_1(y)A\left(\frac{\Lambda}{\rho_1(y)^{1/d}}\right),
\end{equation}
where, 
\[\Lambda := \frac{tr\,\nabla S(y)}{d}.\]
Now, we use (\ref{energy0bis2}) and the convexity of $A$, to obtain that
\begin{eqnarray}
\label{convana14}
\rho_1(y)A\left(\frac{\Lambda}{\rho_1(y)^{1/d}}\right) &\geq& \rho_1(y)\left[\,A\left(\frac{1}{\rho_1(y)^{1/d}}\right)+A^\prime\left(\frac{1}{\rho_1(y)^{1/d}}\right)\left(\frac{\Lambda-1}{\rho_1(y)^{1/d}}\right)\,\right] \nonumber\\
&=& \rho_1(y)\left[\,\frac{G_1\left(\rho_1(y)\right)}{\rho_1(y)} - d(\Lambda-1)\frac{P_G\left(\rho_1(y)\right)}{\rho_1(y)}\,\right] \nonumber\\
&=& G_1\left(\rho_1(y)\right) - P_G\left(\rho_1(y)\right)\,tr\,(\nabla S(y)-\mbox{id}).
\end{eqnarray}
Combining (\ref{energy2}) - (\ref{convana14}), we conclude that
\begin{eqnarray*}
\int_{\Omega}G\left(\rho_0(y)\right)\,\mbox{d}y - \int_\Omega G\left(\rho_1(y)\right)\,\mbox{d}y &\geq& -\int_\Omega P_G\left(\rho_1(y)\right)\,tr\,(\nabla S(y)-\mbox{id})\,\mbox{d}y\\
&=&-\int_\Omega P_G\left(\rho_1(y)\right)\,\mbox{div}\,(S(y)-y)\,\mbox{d}y.
\end{eqnarray*}
Next, assume that $P_G(\rho_1)\in W^{1,\infty}(\Omega)$ and $\rho_1>0$ a.e. Since $P_G\geq 0$, we can approximate $P_G(\rho_1)$ by nonnegative functions in $C^\infty_c(\R^d)$. We use Proposition \ref{propappendix1} - (\roman{iv}), to obtain that
\begin{eqnarray}
\label{energy3bis}
-\int_\Omega P_G(\rho_1(y))\,\mbox{d}iv\,(S(y)-y)\,\mbox{d}y &\geq& \int_\Omega \langle\, \nabla [P_G(\rho_1(y))],S(y)-y\,\rangle\,\mbox{d}y\\ 
&=& \int_\Omega \langle\,\nabla [G^\prime(\rho_1(y))],S(y)-y\,\rangle\,\rho_1(y)\,\mbox{d}y.\nonumber
\end{eqnarray}
We combine (\ref{convana12}) and (\ref{energy3bis}) to conclude (\ref{energy0bis}) \hfill$\Box$\\

The next theorem extends the energy inequality (\ref{energy0bis}) to general cost functions $c$.

\begin{theorem}
\label{theoenergyrel}
({\it Energy inequality for general cost functions}).\\
Let $\rho_0,\,\rho_1\in{\cal P}_a(\Omega)$ be such that $\rho_1>0$ a.e., and $c:\R^d\rightarrow [0,\infty)$ be strictly convex, of class $C^1$ and satisfy $c(0)=0$ and (HC\ref{c4}). Let $G:[0,\infty)\rightarrow \R$ be differentiable on $(0,\infty)$, such that $G(0)=0,\; (0,\infty)\ni x\mapsto x^d\,G(x^{-d})$ be convex and nonincreasing, $\nabla\left(G^\prime(\rho_1)\right)\in L^\infty(\Omega)$, and $P_G(\rho_1)\in W^{1,\infty}(\Omega)$. Denote by $S$, the $c$-optimal map that pushes $\rho_1$ forward to $\rho_0$. Then,
\begin{equation}
\label{en1}
\int_\Omega G(\rho_0(y))\,\mbox{d}y - \int_\Omega G(\rho_1(y))\,\mbox{d}y \geq \int_\Omega \langle\nabla [\,G^\prime\left(\rho_1(y)\right)\,], S(y)-y\rangle \rho_1(y)\,\mbox{d}y. 
\end{equation}
\end{theorem}
{\tt Proof:} Let $(c_k)_k$ be a sequence of regular cost functions satisfying (\ref{exibc12}). By Proposition \ref{propenergy}, we have that
\begin{equation}
\label{energyrel9}
\int_\Omega G(\rho_0(y))\,\mbox{d}y - \int_\Omega G(\rho_1(y))\,\mbox{d}y \geq \int_\Omega \langle\,\nabla\left(G^\prime\left(\rho_1(y)\right)\right) ,S_k(y)-y\,\rangle \rho_1(y)\,\mbox{d}y,
\end{equation}
for all $k\in\N$, where $S_k$ denotes the $c_k$-optimal map that pushes $\rho_1$ forward to $\rho_0$. We let $k$ go to $\infty$ in (\ref{energyrel9}), and we use that $\nabla\left(G^\prime(\rho_1)\right)\in L^\infty(\Omega)$, and $(S_k)_k$ converges to $S$ in $L^2_{\rho_1}(\Omega, \R^d)$ (\cite{agueh:thesis}, Lemma 2.2.2), to conclude (\ref{en1}); here $L^2_{\rho_1}(\Omega, \R^d)$ denotes the set of measurable functions $\varphi:\Omega\rightarrow \R^d$ whose square are summable with respect to the measure $\mu_1:=\rho_1(y)\,\mbox{d}y$, i.e. $\int_\Omega |\,\varphi(y)\,|^2\rho_1(y)\,\mbox{d}y <\infty$. \hfill$\Box$

\subsection{Approximate solutions to the parabolic equation}
Throughout this section, we assume that $\rho_0+\frac{1}{\rho_0}\in L^\infty(\Omega)$. For fixed $h>0$ and $i\in\N$, we denote by $\rho^h_i$ the minimizer of
\begin{equation}
\label{pih} 
(P^h_i):\quad \inf\,\Big\{hW^h_c (\rho^h_{i-1},\rho)+E_i(\rho):\quad \rho\in{\cal P}_a(\Omega)\Big\},
\end{equation}
where $\rho^h_0:=\rho_0$ (see Proposition \ref{propconvana}). We define the approximate solution $\rho^h$ to (\ref{intro8}), as
\begin{equation}
\label{inter}
\rho^h(t,x) := \left\{ \begin{array}{lcl}
\rho_0(x) & \mbox{if} & t=0 \\ \\
\rho^h_i(x) & \mbox{if} & t\in(t_{i-1},t_i],
\end{array}
\right.
\end{equation}
where $t_i=ih,\;i\in\N$. The next proposition shows that
\[ \frac{\partial \rho^h}{\partial t} = \mbox{div}\,\Big\{\rho^h\,\nabla c^\star\left[\,\nabla\left(F^\prime(\rho^h)\right)\,\right]\Big\} + \Lambda(h),\]
in a weak sense. We show in the next section, that 
\[ \|\,\Lambda(h)\,\|_{(W^{2,\infty}(\Omega))^\star} = 0\left(h^{\epsilon(q)}\right),\]
where $\epsilon(q):=\min(1,q-1)$.

\begin{prop}
\label{propinterpolation}
Assume that $F:[0,\infty)\rightarrow \R$ is strictly convex and satisfies $F\in C^2\left((0,\infty)\right)$, and $c:\R^d\rightarrow [0,\infty)$ is strictly convex, of class $C^1$, and satisfies (HC\ref{c2}) - (HC\ref{c3}). Then 
\begin{eqnarray}
\label{exibc25}
\Big|\int_0^T\int _{\Omega} (\rho_0-\rho^h)\,\partial^h_t \xi\,\mbox{d}x\,\mbox{d}t &+& \int_0^T\int_{\Omega} \langle\,\rho^h \nabla c^\star\left[\,\nabla\left(F^{\prime}(\rho^h)\right)\,\right],\nabla \xi \,\rangle \,\mbox{d}x\,\mbox{d}t \Big| \\ 
&\leq& \frac{1}{2}\,\sup_{[0,T]\times\bar{\Omega}}\,\Big|D^2 \xi(t,x)\Big|\; \sum _{i=1}^{T/h} \int _{\Omega \times \Omega}| x-y|^2\, d\gamma_i^h (x,y),\nonumber
\end{eqnarray}
where, $\xi:\R\times\Omega\rightarrow \R$ is such that $\xi(t,.)\in C^2(\bar{\Omega})$ for $t\in\R$, and $\mbox{spt}\,\xi(.,x)\subset [-T,T]$ for $x\in\Omega$, and for some $T>0$. Here,
\[ \partial ^h_t\xi(t,x) := \frac{\xi(t+h,x) - \xi(t,x)}{h},\]
and $\gamma^h_i$ is the $c_h$-optimal measure in $\Gamma(\rho^h_{i-1}, \rho^h_i)$.
\end{prop}
{\tt Proof:} Without loss of generality, we assume that $\frac{T}{h}\in\N$. Because of (\ref{exibc4}), we have that
\[\Big|\,\int_{\Omega} A^h_i(t,x)\,\mbox{d}x\,\Big|\leq B^h_i,\]
for $t\in (0,T)$, where
\[A^h_i(t,x) := \frac{\rho_i^h(x) - \rho^h_{i-1}(x)}{h} \;\xi (t,x) +\left\langle\,\rho_i^h(x)\nabla c^{\star}\left[\,\nabla \left(F^{\prime}\left(\rho_i^h(x)\right)\right)\,\right],\nabla \xi (t,x)\,\right\rangle,\]
and 
\[B^h_i \,:=\, \frac{1}{2h}\,\sup_{[0,T]\times\bar{\Omega}}\,\Big|\,D^2\xi(t,x)\,\Big|\; \int_{\Omega \times \Omega }|\,x-y\,|^2\, d\gamma_i^h (x,y).\]
We integrate the above inequality over $t\in(0,T)$, to obtain that
\begin{equation}
\label{exibc26}
\Big|\,\sum _{i=1}^{T/h}\, \int _{t_{i-1}}^{t_i}\,\mbox{d}t\int _{\Omega} A^h_i(t,x)\,\mbox{d}x\,\Big|\; \leq\, h \sum _{i=1}^{T/h} B^h_i.
\end{equation}
The right hand side of (\ref{exibc26}) gives that
\begin{equation}
\label{exibc27}
h \sum _{i=1}^{T/h} B^h_i\,=\, \frac{1}{2}\sup _{[0,T]\times\bar{\Omega}}\,\Big|\,D^2\xi(t,x)\,\Big|\;\sum _{i=1}^{T/h} \int_{\Omega \times \Omega}|\,x-y\,|^2 d\gamma_i^h(x,y),
\end{equation}
while, on the left hand side, we have that
\begin{eqnarray}
\label{exibc28}
\sum_{i=1}^{T/h}\,\int _{t_{i-1}}^{t_i}\int_\Omega A^h_i(t,x)\,\mbox{d}x\,\mbox{d}t &=& \sum_{i=1}^{T/h} \int_{t_{i-1}}^{t_i} \int_\Omega \frac{\rho_i^h(x) - \rho^h_{i-1}(x)}{h}\,\xi (t,x)\,\mbox{d}x\,\mbox{d}t \\
& & +\; \int_0^T\int_\Omega\left\langle\,\rho^h \nabla c^{\star}\left[\nabla \left(F^{\prime}(\rho^h)\right)\right],\nabla \xi\,\right\rangle\,\mbox{d}x\,\mbox{d}t.\nonumber
\end{eqnarray}
By a direct computation, the first term on the right hand side of (\ref{exibc28}) gives that
\begin{eqnarray*}
\sum_{i=1}^{T/h} \int _{t_{i-1}}^{t_i} \int_\Omega \frac{\rho_i^h(x) - \rho^h_{i-1}(x)}{h}\,\xi (t,x)\,\mbox{d}x\,\mbox{d}t &=& \frac{1}{h}\int_0^T \int_{\Omega} \rho^h (t,x)\,\xi (t,x)\,\mbox{d}x\,dt \\
& & - \frac{1}{h}\,\sum_{i=2}^{T/h}\int _{t_{i-1}}^{t_i}\int _{\Omega} \rho^h(\tau-h,x)\,\xi (\tau,x)\,\mbox{d}x\,\mbox{d}\tau\\
& & - \frac{1}{h}\int _0^h\int _{\Omega}\rho_0(x)\,\xi (t,x)\,\mbox{d}x\,\mbox{d}t.
\end{eqnarray*}
We use the substitution $\tau =t+h$ in the above expression, to obtain that 
\begin{eqnarray*}
\lefteqn{\sum_{i=1}^{T/h} \int _{t_{i-1}}^{t_i} \int_\Omega \frac{\rho_i^h(x) - \rho^h_{i-1}(x)}{h}\,\xi (t,x)\,\mbox{d}x\,\mbox{d}t}\\
& & = \frac{1}{h}\,\int_0^T\int_{\Omega} \rho^h(t,x)\,\xi (t,x)\,\mbox{d}x\,\mbox{d}t\,-\, \frac{1}{h}\int_0^{T-h}\, \rho^h(t,x)\,\xi(t+h,x)\,\mbox{d}x\,\mbox{d}t \\
& &\quad -\,\frac{1}{h}\, \int_0^h \int_{\Omega} \rho_0(x)\,\xi (t,x)\,\mbox{d}x\,\mbox{d}t \\
& & = -\,\int_0^T \int_{\Omega} \rho^h(t,x)\,\partial ^h_t\,\xi (t,x)\,\mbox{d}x\,\mbox{d}t \,+\, \frac{1}{h}\,\int_{T-h}^T \rho^h(t,x)\,\xi (t+h,x)\\
& & \quad -\,\frac{1}{h}\,\int_0^h\int _{\Omega} \rho_0(x)\,\xi (t,x)\, dt\,\mbox{d}x.
\end{eqnarray*}
Noting that
\[-\frac{1}{h} \int_0^h \int _{\Omega}\rho_0(x)\,\xi (t,x)\,\mbox{d}t\,\mbox{d}x =\int_0^T\int_\Omega \rho_0(x)\,\partial^h_t\xi (t,x)\,\mbox{d}x\,\mbox{d}t, \] 
and $\xi(t+h)=0$ for $t\in(T-h,T)$, we deduce that
\begin{eqnarray}
\label{exibc29}
\lefteqn{\sum_{i=1}^{T/h} \int _{t_{i-1}}^{t_i} \int_\Omega \frac{\rho_i(x) - \rho^h_{i-1}(x)}{h}\,\xi (x,t)\,\mbox{d}x\,\mbox{d}t} \\ 
& &= \int_0^T\int_{\Omega} \left(\rho_0(x)-\rho^h(t,x)\right)\,\partial ^{h}_{t}\xi(t,x)\,\mbox{d}x\,\mbox{d}t. \nonumber
\end{eqnarray}
We combine (\ref{exibc26})  - (\ref{exibc29}) to conclude (\ref{exibc25}) \hfill $\Box$

\section{Existence and uniqueness of solutions}
Below, we study the limit of (\ref{exibc25}), as $h$ goes to $0$. The first three sections deal with the limits of the three terms of inequality (\ref{exibc25}), and the last section proves the existence theorem to (\ref{intro8}), when $\rho_0$ is bounded below and above. Here and after, $\rho^h$ is defined as in (\ref{inter}).

\subsection{Second moments of the optimal measures}
We show that
\begin{equation}
\label{exisv1}
\sum_{i=1}^{T/h}\int_{\Omega\times\Omega}|x-y|^2\,\mbox{d}\gamma^h_i(x,y) = 0(h^{\epsilon(q)}),
\end{equation}
where $\epsilon(q):=\min(1,q-1),\,\gamma^h_i$ denotes the $c_h$-optimal measure in $\Gamma(\rho^h_{i-1},\rho^h_i)$, and $\rho^h_i$ is the unique minimizer of (\ref{pih}). The first step toward proving (\ref{exisv1}) is the next lemma, which states that $\sum_{i=1}^\infty W^h_c(\rho^h_{i-1},\rho^h_i)$ is bounded, uniformly in $h$.

\begin{lemma}
\label{lemlimsecond}
Assume that $F:[0,\infty)\rightarrow \R$ and $c:\R^d\rightarrow [0,\infty)$ are strictly convex, and $c$ satisfies (HC\ref{c2}). Then
\begin{equation}
\label{exisv1bis}
\sum_{i=1}^\infty h\,W^h_c(\rho^h_{i-1},\rho^h_i) \leq E_i(\rho_0) - |\,\Omega\,|F\left(\frac{1}{|\,\Omega\,|}\right).
\end{equation}
\end{lemma}
{\tt Proof:} Let $T>0$, be such that $\frac{T}{h}\in\N$. Since $c(0)=0$, Proposition \ref{propconvana} gives that
\[ h\,W^h_c(\rho^h_{i-1},\rho^h_i) \leq E_i(\rho^h_{i-1})-E_i(\rho^h_i),\]
for $i\in\N$. We sum both sides of the above inequality over $i$, to obtain that
\[\sum_{i=1}^{T/h} h\,W^h_c(\rho^h_{i-1},\rho^h_{i}) \leq E_i(\rho_0) - \int_\Omega F(\rho^h_{T/h}(x))\,\mbox{d}x.\]
We apply Jensen's inequality to the integral term above, and we let $T$ go to $\infty$, to conclude (\ref{exisv1bis})\hfill $\Box$

\begin{prop}
\label{proplimsecond}
Assume that $F:[0,\infty)\rightarrow \R$ and $c:\R^d\rightarrow [0,\infty)$ are strictly convex, and $c$ satisfies $c(0)=0$ and $c(z)\geq\beta\,|\,z\,|^q$, for some $\beta>0$ and $q>1$. Then, for $T>0$ and $h\in(0,1)$ such that $\frac{T}{h}\in\N$,
\begin{equation}
\label{exibc30}
\sum_{i=1}^{T/h} \int_{\Omega\times\Omega}|\,x-y\,|^2\,\mbox{d}\gamma_i^h(x,y) \leq  M(\Omega,T,F,\rho_0,q,\beta)\,h^{\epsilon(q)},
\end{equation}
where $\epsilon(q):=\min(1,q-1)$.
\end{prop}
{\tt Proof.} Since $c(z)\geq \beta\,|z|^q$, we have that
\begin{equation}
\label{exibc31}
\int_{\Omega\times\Omega}|\,x-y\,|^q\,\mbox{d}\gamma^h_i(x,y) \,\leq\, \frac{h^q}{\beta}\;W^h_c(\rho^h_{i-1},\rho^h_i),
\end{equation}
for $i\in\N$. We distinguish two cases, based on the values of $q$.\\
{\underline{\em Case 1:}}\hspace{1mm} $1 < q\leq 2$.\\
Because of (\ref{exibc31}), we have, for $i\in \N$, that
\begin{eqnarray*}
\int_{\Omega\times\Omega}|\,x-y\,|^2\,\mbox{d}\gamma^h_i(x,y) &\leq& \sup_{x,y\in\bar{\Omega}}\,|\,x-y\,|^{(2-q)}\,\int_{\Omega\times\Omega}|\,x-y\,|^q\,\mbox{d}\gamma^h_i(x,y)\\
&\leq & \frac{(\mbox{diam}\,\Omega)^{(2-q)}}{\beta}\,h^q\,W^h_c(\rho^h_{i-1},\rho^h_i),
\end{eqnarray*}
where, $\mbox{diam}\,\Omega$ denotes the diameter of $\Omega$. We sum both sides of the above inequality over $i$, and we use (\ref{exisv1bis}), to conclude that
\[ \sum_{i=1}^{T/h} \int_{\Omega\times\Omega}|\,x-y\,|^2\,\mbox{d}\gamma_i^h(x,y)\, \leq M(\Omega, F, \rho_0, q, \beta)\,h^{q-1}.\]

\noindent{\underline{\em Case 2:}}\hspace{1mm} $q > 2$.\\
Because of Jensen's inequality and (\ref{exibc31}), we have that
\[\int_{\Omega\times\Omega}|\,x-y\,|^2\,\mbox{d}\gamma^h_i(x,y) \leq  \left(\int_{\Omega\times\Omega}|\,x-y\,|^q\,\mbox{d}\gamma^h_i(x,y)\right)^{2/q} \leq \frac{h^2}{\beta^{2/q}}\,\left[\,W^h_c(\rho^h_{i-1},\rho^h_i)\,\right]^{2/q}.\]
We sum both sides of the above inequality over $i$, and we use H\"older's inequality on the right hand side term, to obtain that
\begin{eqnarray}
\label{exisv2}
\sum_{i=1}^{T/h}\int_{\Omega\times\Omega}|\,x-y\,|^2\,\mbox{d}\gamma_i^h(x,y) &\leq& \frac{h^2}{\beta^{2/q}}\,\left(\frac{T}{h}\right)^{1-\frac{2}{q}}\,\left[\,\sum_{i=1}^{T/h} W^h_c(\rho^h_{i-1},\rho^h_i)\,\right]^{2/q} \nonumber\\
& & =T^{1-\frac{2}{q}}\frac{h}{\beta^{2/q}}\,\left[\,\sum_{i=1}^{T/h} h\,W^h_c(\rho^h_{i-1},\rho^h_i)\,\right]^{2/q}.
\end{eqnarray}
We combine (\ref{exisv1bis}) and (\ref{exisv2}), to conclude that
\[\sum_{i=1}^{T/h}\int_{\Omega\times\Omega}|\,x-y\,|^2\,\mbox{d}\gamma_i^h(x,y)  \leq M(\Omega,T, F, \rho_0, q,\beta)\,h  \qquad\qquad\Box \]

\subsection{Strong convergence of the approximate solutions}
We prove that $(\rho^h)_h$ is compact in $L^1(\Omega_T)$, for $0<T<\infty$. The main ingredient in the proof is the energy inequality (\ref{en1}). It allows us to obtain a uniform bound in $h$, of the $L^{q^\star}$- norm of $\nabla \left(F^\prime(\rho^h)\right),\,$ which leads to the compactness of $(\rho^h)_h$ in  $L^1(\Omega_T)$. We first show that $(\rho^h)_h$ converges weakly in $L^1(\Omega_T)$ for a subsequence. We introduce the following constant needed in the next lemma:
\[\overline{M}(\Omega, T, F, \rho_0, q, \alpha) := M(\alpha, q)\left(E_i(\rho_0)-|\,\Omega\,|F\left(\frac{1}{|\,\Omega\,|}\right)+\alpha T |\,\Omega\,|\,\|\,\rho_0\,\|_{L^\infty(\Omega)}\right),\]
where $M(\alpha, q)$ is a constant which depends only on $\alpha$ and $q$.

\begin{lemma}
\label{lemstrongconv1}
Assume that  $c:\R^d\rightarrow [0,\infty)$ is strictly convex, of class $C^1$ and satisfies (HC\ref{c2}), and $F:[0,\infty)\rightarrow \R$ is strictly convex, of class $C^2\left((0,\infty)\right)$. If $\rho_0\in{\cal P}_a(\Omega)\cap L^\infty(\Omega)$, then, 
\begin{equation}
\label{strongconv2.1}
\|\,\rho^h\,\|_{L^\infty\left((0,\infty);L^\infty(\Omega)\right)} \leq \|\,\rho_0\,\|_{L^\infty(\Omega)}.
\end{equation}
Therefore, there exists $\rho:[0,\infty)\times\Omega \rightarrow \R$ and a subsequence of $(\rho^h)_{h\downarrow 0}$, which converges to $\rho$, weakly in $L^1(\Omega_T)$, for $\,0<T<\infty$.\\
In addition, if $\frac{1}{\rho_0}\in L^\infty(\Omega)$, $c$ satisfies (HC\ref{c4}), and $F$ satisfies $F(0)=0$ and (HF\ref{f2}), then
\begin{equation}
\label{strongconv2.3}
\int_{\Omega_T} \rho^h\Big|\,\nabla\left(F^\prime(\rho^h)\right)\,\Big|^{q^\star} \leq \overline{M}(\Omega, T, F, \rho_0, q, \alpha).\end{equation}
\end{lemma}
{\tt Proof:} Because of the second bound in (\ref{bba}) - a consequence of the maximum principle of Proposition \ref{propbound} -, we have that $\rho^h_i\leq \|\,\rho_0\,\|_{L^\infty(\Omega)}$ for $i\in\N$, which reads as $\|\,\rho^h(t)\,\|_{L^\infty(\Omega)} \leq \|\,\rho_0\,\|_{L^\infty(\Omega)}$ for $t\in [0,\infty)$. We take the supremum of the previous inequality over $t\in (0,\infty)$, to deduce (\ref{strongconv2.1}).\\
Due to (\ref{strongconv2.1}), we have that $(\rho^h)_h$ is precompact in $L^1(\Omega_T)$, for $0<T<\infty$. We use the standard diagonal argument, to conclude that  $(\rho^h)_{h\downarrow 0}$ converges weakly to some function $\rho:[0,\infty)\times\Omega \rightarrow \R$ in $L^1(\Omega_T)$, for a subsequence.\\
Because of Proposition \ref{propproperty}, (\ref{bba}), and the fact that $\nabla\left(P(\rho^h_i)\right)=\rho^h_i\nabla\left(F^\prime(\rho^h_i)\right)$, we have that $P(\rho^h_i)\in W^{1,\infty}(\Omega)$ and $\nabla\left(F^\prime(\rho^h_i)\right)\in L^\infty(\Omega)$, for $i\in\N.$ Then, we choose $G:=F$ in the energy inequality (\ref{en1}), and we use (\ref{exibc3}), to obtain that
\[ h\int_\Omega \left\langle\,\nabla\left(F^\prime(\rho^h_i)\right), \nabla c^\star \left[\,\nabla\left(F^\prime(\rho^h_i)\right)\,\right]\,\right\rangle \rho^h_i \leq \int_\Omega F(\rho^h_{i-1}) - \int_\Omega F(\rho^h_i).\]
We sum both sides of the subsequent inequality over $i$, and we use Jensen's inequality, to deduce that
\begin{equation}
\label{strongconv2.4}
\int_{\Omega_T} \left\langle\,\nabla\left(F^\prime(\rho^h)\right), \nabla c^\star \left[\,\nabla\left(F^\prime(\rho^h)\right)\,\right]\,\right\rangle \rho^h \leq \int_\Omega F(\rho_0) - |\,\Omega\,|\,F\left(\frac{1}{|\,\Omega\,|}\right).
\end{equation} 
Because of (\ref{appendn3}) of Proposition \ref{lemappend3}, and the fact that $c(z)\leq \alpha\,(|\,z\,|^q+1),$ we have that
\[\left\langle\, z, \nabla c^\star (z)\,\right\rangle \geq c^\star(z) \geq  M(\alpha, q)\,|z\,|^{q^\star}\,-\,\alpha,\]
and then, (\ref{strongconv2.4}) implies that
\begin{equation}
\label{strongconv2.5}
M(\alpha, q)\,\int_{\Omega_T}\rho^h\,\Big|\,\nabla\left(F^\prime(\rho^h)\right)\,\Big|^{q^\star} \leq \int_\Omega F(\rho_0) - |\,\Omega\,|\,F\left(\frac{1}{|\,\Omega\,|}\right) + \alpha\int_{\Omega_T}\rho^h.
\end{equation}
We combine (\ref{strongconv2.1}) and (\ref{strongconv2.5}), to obtain that
\[M(\alpha, q)\,\int_{\Omega_T}\rho^h\,\Big|\,\nabla\left(F^\prime(\rho^h)\right)\,\Big|^{q^\star} \leq \int_\Omega F(\rho_0) - |\,\Omega\,|\,F\left(\frac{1}{|\,\Omega\,|}\right) + \alpha\,T|\,\Omega\,|\,\|\,\rho_0\,\|_{L^\infty(\Omega)}.\]
We divide both sides of the above inequality by $M(\alpha, q)$ to conclude (\ref{strongconv2.3})\hfill$\Box$\\

\begin{lemma}
\label{lemstrongconv4}
{\it (Space-compactness)}\\
Assume that  $c:\R^d\rightarrow [0,\infty)$ is strictly convex, of class $C^1$ and satisfies $c(0)=0$ and (HC\ref{c4}), and $F:[0,\infty)\rightarrow \R$ is strictly convex, of class $C^2\left((0,\infty)\right)$, and satisfies $F(0)=0$ and (HF\ref{f2}). If $\rho_0\in{\cal P}_a(\Omega)$ is such that $\rho_0+\frac{1}{\rho_0}\in L^\infty(\Omega)$, then, for $\eta\not = 0$ and $0<T<\infty$,
\begin{equation}
\label{spacecomp1}
\int_{\Omega^{(\eta)}_T} \Big|\,\rho^h(t,x+\eta e) - \rho^h(t,x)\,\Big| \leq M(\Omega, T, F, \rho_0, \alpha, q)\,|\eta|,
\end{equation}
where $e$ is a unit vector of $\R^d,\;\Omega^{(\eta)}:=\{\,x\in\Omega:\;{\mbox dist}\,(x,\partial\Omega)>|\,\eta\,|\,\}$, and  $\Omega_T^{(\eta)}:=(0,T)\times\Omega^{(\eta)}$.
\end{lemma}
{\tt Proof:} Since $\rho_0+\frac{1}{\rho_0}\in L^\infty(\Omega)$, (\ref{bba}) implies that $(\rho^h)_h$ is bounded below and above. Then, we use that $F\in C^2\left((0,\infty)\right)$, to obtain that
\begin{eqnarray}
\label{spacecomp2}
\Big\|\,\nabla \rho^h\,\Big\|^{q^\star}_{L^{q^\star}(\Omega_T)} &=& \int_{\Omega_T} \frac{1}{\rho^h\left[F^{\prime\prime}(\rho^h)\,\right]^{q^\star}}\,\rho^h\,\Big|\,\nabla\left(F^\prime(\rho^h)\right)\,\Big|^{q^\star}\nonumber\\
&\leq& M(\Omega,\rho_0,F)\int_{\Omega_T} \rho^h\,\Big|\,\nabla\left(F^\prime(\rho^h)\right)\,\Big|^{q^\star}.
\end{eqnarray}
We combine (\ref{strongconv2.3}) and (\ref{spacecomp2}), to conclude that $\left(\nabla \rho^h\right)_h$ is bounded in $L^{q^\star}(\Omega_T)$. As a consequence, we have that $(\rho^h)_h$ is bounded in $W^{1,q^\star}(\Omega_T)$. Approximating $\rho^h$ by $C^\infty(\Omega_T)$-functions, and using the mean-value theorem, and the fact that $(\nabla \rho^h)_h$ is bounded in $L^{q^\star}(\Omega_T)$, we have that
\begin{equation}
\label{spacecomp4}
\int_{\Omega_T^{(\eta)}} |\,\rho^h(t,x+\eta e)-\rho^h(t,x)\,|^{q^\star} \leq M(\Omega, T, F, \rho_0, \alpha, q)\, |\,\eta\,|^{q^\star}.
\end{equation}
We combine (\ref{spacecomp4}) and H\"older's inequality, to conclude that
\begin{eqnarray*}
\int_{\Omega_T^{(\eta)}} \Big|\,\rho^h(t,x+\eta e)-\rho^h(t,x)\,\Big| &\leq& |\,\Omega_T\,|^{1/q}\,\left(\int_{\Omega_T^{(\eta)}} \Big|\,\rho^h(t,x+\eta e)-\rho^h(t,x)\,\Big|^{q^\star}\right)^{1/q^\star}\\
&\leq& M(\Omega, T, F, \rho_0, \alpha, q)\,|\,\eta\,| \qquad\qquad \Box
\end{eqnarray*}

Now, we focus on the time-compactness of $(\rho^h)_h$ on $\Omega_T,\,0<T<\infty$. The following constant will be needed in the next lemma:
\begin{eqnarray*}
\lefteqn{ \overline{\overline{M}}(\Omega, T, F, \rho_0, q, \alpha, \beta)}\\
& & :=\frac{\|\,\rho_0\|_{L^\infty(\Omega)}^{1/q^\star}}{\Big\|\,\frac{1}{\rho_0}\,\Big\|_{L^\infty(\Omega)}^{1/q^\star}}\,M(q, \alpha, \beta)\left(E_i(\rho_0)-|\,\Omega\,|F\left(\frac{1}{|\,\Omega\,|}\right)+\alpha T |\,\Omega\,|\,\|\,\rho_0\,\|_{L^\infty(\Omega)}\right),
\end{eqnarray*}
where $M(q,\alpha,\beta)$ is a constant which only depends on $q,\,\alpha$ and $\beta$.

\begin{lemma}
\label{lemstrongconv5}
Assume that  $c:\R^d\rightarrow [0,\infty)$ is strictly convex, of class $C^1$ and satisfies $c(0)=0$ and (HC\ref{c4}), and $F:[0,\infty)\rightarrow \R$ is strictly convex, of class $C^2\left((0,\infty)\right)$, and satisfies $F(0)=0$ and (HF\ref{f2}). If $\rho_0\in{\cal P}_a(\Omega)$ is such that $\rho_0+\frac{1}{\rho_0}\in L^\infty(\Omega)$, then, for $\tau>0$ and $0<T<\infty$,
\begin{eqnarray}
\label{timecomp1}
\int_{\Omega_T} \left[\,F^\prime\left(\rho^h(t+\tau,x)\right) - F^\prime\left(\rho^h(t,x)\right)\,\right]\,\left[\,\rho^h(t+\tau,x) - \rho^h(t,x)\,\right] \nonumber\\ 
\leq \overline{\overline{M}}(\Omega, T, F, \rho_0, q, \alpha, \beta)\,\tau.
\end{eqnarray}
\end{lemma}
{\tt Proof:} Without loss of generality, we assume that $\frac{T}{h}\in\N$ and $\tau=Nh$, for some $N\in\N$. For simplicity, we set
\[L(h,\tau) := \int_{\Omega_T} \left[\,F^\prime\left(\rho^h(t+\tau,x)\right) - F^\prime\left(\rho^h(t,x)\right)\,\right]\,\left[\,\rho^h(t+\tau,x) - \rho^h(t,x)\,\right],\]  
and
\[ J(i,h,N):=\int_\Omega \left[\,F^\prime\left(\rho^h_{i+N}(x)\right) - F^\prime\left(\rho^h_i(x)\right)\,\right]\,\left[\,\rho^h_{i+N}(x) - \rho^h_i(x)\,\right].\]
It is straightforward to check that
\begin{equation}
\label{strongcon5.1}
L(h,\tau) = \sum_{i=1}^{T/h} h\,J(i,h,N).
\end{equation}
Since $(W^h_c)^{1/q}$ does not satisfy the triangle inequality, we introduce the $q$-Wasserstein metric $d^h_q:=(W^h_q)^{1/q}$, defined by  
\begin{equation}
\label{strongconv5.2}
d^h_q(\rho^h_{i}, \rho^h_{i+N}):= \left(\int_\Omega \Big|\,\frac{y-S^h_q(y)}{h}\,\Big|^q \rho^h_{i+N}(y)\,\mbox{d}y\right)^{1/q},
\end{equation}
where $S^h_q$ denotes the $|\,\frac{\cdot}{h}\,|^q$-optimal map that pushes $\rho^h_{i+N}$ forward to $\rho^h_{i}$. Then, setting $\varphi^h_{i,N} := F^\prime\left(\rho^h_{i+N}\right) - F^\prime\left(\rho^h_{i}\right)$, we obtain that
\[ J(i,h,N) = \int_\Omega\left[\;\varphi^h_{i,N}(y)-\varphi^h_{i,N}\left(S^h_q(y)\right)\,\right]\,\rho^h_{i+N}(y)\,\mbox{d}y.\]
Since $\rho_0+\frac{1}{\rho_0}\in L^\infty(\Omega),\;F\in C^2\left((0,\infty)\right)$, and $\rho^h_i\nabla\left(F^\prime(\rho^h_i)\right) = \nabla\left(P(\rho^h_i)\right) \in L^\infty(\Omega)$ (see Proposition \ref{propproperty} - (\roman{mycounter})), (\ref{bba}) gives that $\varphi^h_{i,N}\in W^{1,\infty}(\Omega)$. So, approximating $\varphi^h_{i,N}$ by $C^\infty(\Omega)$-functions, and using that $\left(S^h_q\right)_{\#}\rho^h_{i+N}=\rho^h_i$ and the mean-value theorem, we rewrite $J(i,h,N)$ as follows:
\[J(i,h,N) = \int_\Omega\int_0^1 \left\langle\,\nabla \varphi^h_{i,N}\left((1-t)y + tS^h_q(y)\right), y-S^h_q(y)\,\right\rangle\,\rho^h_{i+N}(y)\,\mbox{d}t\,\mbox{d}y.\]
We combine H\"older's inequality and (\ref{strongconv5.2}), to deduce that
\begin{equation}
\label{strongconv5.3}
J(i,h,N)\leq h\,d^h_q\left(\rho^h_{i},\rho^h_{i+N}\right)\left[\,\int_\Omega\int_0^1 \Big|\,\nabla \varphi^h_{i+N}\left((1-t)y+tS^h_q(y)\right)\,\Big|^{q^\star}\rho^h_{i+N}(y)\,\mbox{d}t\,\mbox{d}y\,\right]^{1/q^\star} 
\end{equation}
But, observe that $\rho^h_{i},\,\rho^h_{i+N}\leq \|\,\rho_0\,\|_{L^\infty(\Omega)}$ because of Proposition \ref{propbound}, and $|\,\nabla\varphi^h_{i,N}\,|^{q^\star}\in L^\infty(\Omega)$. So, we approximate $|\,\nabla\varphi^h_{i,N}\,|^{q^\star}$ by nonnegative functions in $C^\infty_c(\R^d)$, and we use (\ref{intbound1}) in Proposition \ref{propintbound}, to deduce that
\begin{equation}
\label{strongconv5.4}
\int_\Omega \Big|\,\nabla \varphi^h_{i,N}\left((1-t)y+t\rho^h_q(y)\right)\,\Big|^{q^\star}\rho^h_{i+N}(y)\,\mbox{d}y \leq \|\,\rho_0\,\|_{L^\infty(\Omega)}\,\int_{\R^d} \Big|\,\nabla \varphi^h_{i,N}(y)\,\Big|^{q^\star}\,\mbox{d}y.
\end{equation} 
We combine (\ref{strongcon5.1}), (\ref{strongconv5.3}) and (\ref{strongconv5.4}), to have that
\[ L(h,\tau)\leq \|\,\rho_0\,\|_{L^\infty(\Omega)}^{1/q^\star}\,h^2\sum_{i=1}^{T/h} d^h_q\left(\rho^h_{i}, \rho^h_{i+N}\right)\,\|\,\nabla\varphi^h_{i,N}\,\|_{L^{q^\star}(\Omega)}.\]
And since $d^h_q$ is a metric, the triangle inequality gives that
\[L(h,\tau)\leq  \|\,\rho_0\,\|_{L^\infty(\Omega)}^{\frac{1}{q^\star}}\,h^2\sum_{k=1}^N\sum_{i=1}^{T/h} \|\,\nabla\varphi^h_{i,N}\,\|_{L^{q^\star}(\Omega)}\,d^h_q\left(\rho^h_{i+k-1}, \rho^h_{i+k}\right).\]
Then, we apply H\"older's inequality to the interior sum, to deduce that
\begin{equation}
\label{strongconv5.6}
L(h,\tau)\leq \|\,\rho_0\,\|_{L^\infty(\Omega)}^{\frac{1}{q^\star}}\,h^{2-\frac{1}{q^\star}}\left(\sum_{i=1}^{T/h} h\,\|\,\nabla \varphi^h_{i,N}\,\|^{q^\star}_{L^{q^\star}(\Omega)}\right)^{1/q^\star} \sum_{k=1}^N\left[\,\sum_{i=1}^{T/h} d^h_q\left(\rho^h_{i+k-1}, \rho^h_{i+k}\right)^q\,\right]^{1/q}
\end{equation}
Because of (\ref{bba}) and (\ref{strongconv2.3}), $\left(h^{1/q^\star}\,\Big\|\,\nabla\left(F^\prime(\rho^h_i)\right)\,\Big\|_{L^{q^\star}(\Omega)}\right)_{i=1,\cdots,\frac{T}{h}}$ and \\
 $\left(h^{1/q^\star}\,\Big\|\,\nabla\left(F^\prime(\rho^h_{i+N})\right)\,\Big\|_{L^{q^\star}(\Omega)}\right)_{i=1, \cdots, \frac{T}{h}} \,$ belong to $l_{q^\star}(\Omega)$. Then, we combine H\"older's inequality, Minkowski's inequality, (\ref{bba}) and (\ref{strongconv2.3}), to have that
\begin{eqnarray}
\label{strongconv5.7}
\lefteqn{\left(\sum_{i=1}^{T/h} h\,\Big\|\nabla\varphi^h_{i,N}\,\Big\|^{q^\star}_{L^{q^\star}(\Omega)}\right)^{1/q^\star}}\nonumber\\
& &\leq \left(\sum_{i=1}^{T/h}\left( h^{1/q^\star}\,\Big\|\,\nabla \left(F^\prime(\rho^h_{i+N})\right)\,\Big\|_{L^{q^\star}(\Omega)} + h^{1/q^\star}\, \Big\|\,\nabla \left(F^\prime(\rho^h_{i})\right)\,\Big\|_{L^{q^\star}(\Omega)}\right)^{q^\star}\right)^{1/q^\star} \nonumber \\
& & \leq \left(\sum_{i=1}^{T/h} h\,\Big\|\,\nabla \left(F^\prime(\rho^h_{i+N})\right)\,\Big\|^{q^\star}_{L^{q^\star}(\Omega)}\right)^{1/q^\star} + \left(\sum_{i=1}^{T/h} h\,\Big\|\,\nabla \left(F^\prime(\rho^h_{i})\right)\,\Big\|^{q^\star}_{L^{q^\star}(\Omega)}\right)^{1/q^\star}\nonumber\\
& & \leq \frac{1}{\Big\|\frac{1}{\rho_0}\,\Big\|_{L^\infty(\Omega)}^{1/q^\star}}\,\left[\,\overline{M}\,(\Omega, T, F, \rho_0, q, \alpha)\,\right]^{1/q^\star}.
\end{eqnarray}
On the other hand, since $c(z)\geq \beta\,|\,z\,|^q$, we have that $\left(d^h_q\right)^q \leq \frac{1}{\beta}W^h_c$, and then,
\[\sum_{k=1}^N\left[\,\sum_{i=1}^{T/h} d^h_q\left(\rho^h_{i+k-1}, \rho^h_{i+k}\right)^q\,\right]^{1/q} \leq \frac{1}{(\beta h)^{1/q}}\, \sum_{k=1}^N \left[\,\sum_{i=1}^{T/h} h\,W^h_c\left(\rho^h_{i+k-1}, \rho^h_{i+k}\right)\,\right]^{1/q}.\]
We use (\ref{exisv1bis}) and the above inequality, to deduce that
\begin{equation}
\label{strongconv5.8}
\sum_{k=1}^N\left[\,\sum_{i=1}^{T/h} d^h_q\left(\rho^h_{i+k-1}, \rho^h_{i+k}\right)^q\,\right]^{1/q} \leq \frac{1}{\beta}\,\left[\, E_i(\rho_0) - |\,\Omega\,|F\left(\frac{1}{|\,\Omega\,|}\right)\,\right]^{1/q}\,Nh^{-1/q}.
\end{equation}
We combine (\ref{strongconv5.6}) - (\ref{strongconv5.8}), and we use that $\tau=Nh$, to conclude that
\[ L(h,\tau)\leq \overline{\overline{M}}(\Omega, T, F, \rho_0, q, \alpha, \beta)\,\tau.\qquad\qquad \Box\]

\begin{lemma}
\label{lemstrongconv6} {\it (Time-compactness)}\\
Assume that the assumptions of Lemma \ref{lemstrongconv4} hold. If $\rho_0\in {\cal P}_a(\Omega)$ is such that $\rho_0+\frac{1}{\rho_0}\in L^\infty(\Omega)$, then, for $0<T<\infty$, and small $\tau>0,$ 
\[\int_{\Omega_T}\Big|\,\rho^h(t+\tau,x)-\rho^h(t,x)\,\Big| \leq M(R,\Omega, T, F, \rho_0, \alpha, q, \beta)\sqrt{\tau} + T\Lambda(\sqrt{\tau}), \]
where $\Lambda$ is such that $\,\lim_{\tau \downarrow 0} \Lambda (\sqrt{\tau})=0$.
\end{lemma}
{\tt Proof:} Let $R>0$, and for fixed $h,\,T$ and $\tau$, define
\[ \begin{array}{l}
 E_R:=\Big\{ t\in(0,T):\quad \Delta_{h,\tau}(t):=\Big\|\,\rho^h(t)\,\Big\|_{L^q(\Omega)} + \Big\|\,\rho^h(t+\tau)\,\Big\|_{L^q(\Omega)} \\ \\
\qquad\qquad + \,\Big\|\,F^\prime\left(\rho^h(t)\right)\,\Big\|_{W^{1,q^\star}(\Omega)} +\,\Big\|\,F^\prime\left(\rho^h(t+\tau)\right)\,\Big\|_{W^{1,q^\star}(\Omega)}\\ \\
\qquad\qquad + \frac{1}{\tau}\int_\Omega \left[\,F^\prime\left(\rho^h (t+\tau)\right) - F^\prime\left(\rho^h(t)\right)\,\right]\left[\,\rho^h (t+\tau) - \rho^h (t)\,\right]> R\Big\}.\end{array}\]
Because of (\ref{bba}), (\ref{strongconv2.3}), (\ref{timecomp1}), and the fact that $F\in C^2\left((0,\infty)\right)$, we have that $(0,T)\ni t\mapsto\Delta_{h,\tau}(t)$ belongs to $L^1\left((0,T)\right)$. Hence
\begin{equation}
\label{timecomp2}
|\,E_R\,|\leq \frac{M(\Omega, T, F, \rho_0, q,\alpha, \beta)}{R}.
\end{equation}
We combine (\ref{strongconv2.1}) and (\ref{timecomp2}), to have that
\begin{equation}
\label{timecomp3}
\int_{E_R}\int_\Omega \Big|\,\rho^h(t+\tau,x)-\rho^h(t,x)\,\Big| \leq 2\,\|\,\rho_0\,\|_{L^\infty(\Omega)}\,|\,\Omega\,|\,|\,E_R\,|\leq \frac{M(\Omega, T, F, \rho_0, q, \alpha,\beta)}{R}.
\end{equation}
On the other hand, if $t\in E_R^c:=(0,T)\setminus E_R$, setting $\rho^h(t):=u_1$ and $\rho^h(t+\tau):=u_2$, we clearly have that $\|\,u_i\,\|_{L^q(\Omega)} \leq R, \; \|\,F^\prime(u_i)\,\|_{W^{1,q^\star}(\Omega)} \leq R\,$ for $i=1,2$, and $\int_\Omega \left[\,F^\prime(u_2) - F^\prime(u_1)\,\right]\left[\,u_2 - u_1\,\right] \leq R\,\tau$. Then, Proposition \ref{propinterpolation1} gives that
\begin{equation}
\label{timecomp4}
\int_{E_R^c}\int_\Omega |\,\rho^h(t+\tau,x)-\rho^h(t,x)\,| \leq \int_{E_R^c}\Lambda(R\,\tau) \leq T\Lambda(R\,\tau),
\end{equation}
where $\,\Lambda(R\,\tau):=\Lambda_R(R\,\tau)\,$ is defined as in Proposition \ref{propinterpolation1}. We combine (\ref{timecomp3}) - (\ref{timecomp4}), and we choose $\,R=\frac{1}{\sqrt{\tau}},\,$ to conclude the proof \hfill$\Box$\\

Having proved the space-compactness and time-compactness of $(\rho^h)_h$, we are now ready to show that $(\rho^h)_h$ converges strongly to $\rho$ in $L^1(\Omega_T),\,0<T<\infty,$ for a subsequence, where $\rho$ is defined as in Lemma \ref{lemstrongconv1}.

\begin{prop}
\label{propstrongconv}
Assume that  $c:\R^d\rightarrow [0,\infty)$ is strictly convex, of class $C^1$ and satisfies $c(0)=0$ and (HC\ref{c4}), and $F:[0,\infty)\rightarrow \R$ is strictly convex, of class $C^2\left((0,\infty)\right)$, and satisfies $F(0)=0$ and (HF\ref{f2}). If $\rho_0\in{\cal P}_a(\Omega)$ is such that $\rho_0+\frac{1}{\rho_0}\in L^\infty(\Omega)$, then, for $0<T<\infty$, there is a subsequence of $(\rho^h)_{h\downarrow 0}$ which converges strongly to $\rho$ in $L^r(\Omega_T)$, for $1\leq r<\infty$, where $\rho$ is defined as in Lemma \ref{lemstrongconv1}.
\end{prop} 
{\tt Proof:} Fix $\delta>O$, and define $\Omega_T^{(\delta)}$ as in Lemma \ref{lemstrongconv4}. Because of (\ref{strongconv2.1}), we have that $(\rho^h)_h$ is bounded in $L^1\left(\Omega_T^{(\delta)}\right)$. Furthermore, for $\epsilon>0$, and small $\tau>0$ and $\eta\in(0,\delta)$, we have that $\Omega_T^{(\delta)}\subset \Omega_T^{(\eta)}\subset \Omega_T$, and then, Lemma \ref{lemstrongconv4} and Lemma \ref{lemstrongconv6} give that
\[ \int_{\Omega_T^{(\delta)}} |\,\rho^h(t,x+\eta e)-\rho^h(t,x)\,| < \epsilon, \;\mbox{and}\; \int_{\Omega_T^{(\delta)}} |\,\rho^h(t+\tau,x)-\rho^h(t,x)\,| < \epsilon,\] 
uniformly in $h$. We deduce that, $(\rho^h)_h$ is precompact in $L^1\left(\Omega_T^{(\delta)}\right)$ (see \cite{adams:sobolev}, Theorem 2.21). We observe that $\lim_{\delta\rightarrow 0}|\,\Omega\setminus \Omega^{(\delta)}\,|=0$, and then, we use the diagonal argument, to obtain that  $(\rho^h)_h$ converges strongly to $\rho$ in $L^1(\Omega_T)$, for a subsequence. And since  $(\rho^h)_h$ is bounded in $L^\infty(\Omega_T)$ (see (\ref{strongconv2.1})), we conclude that it converges to $\rho$ in $L^r(\Omega_T)$, for $1\leq r <\infty$ (up to a subsequence) \hfill$\Box$

\subsection{Weak convergence of the nonlinear terms}
We use the energy inequality (\ref{en1}) to show that $\left(\mbox{div}\{\rho^h\nabla c^\star\left[\nabla\left(F^\prime(\rho^h)\right)\right]\}\right)_{h}$ converges weakly to $\mbox{div}\{\rho\nabla c^\star\left[\nabla\left(F^\prime(\rho)\right)\right]\}$ in $\Omega_T$, for a subsequence. Throughout this section, $(\rho^h)_h$  denotes the (non-relabeled) subsequence of $(\rho^h)_h$ which converges to $\rho$ in $L^r(\Omega_T)$, for $1\leq r<\infty$, as in Proposition \ref{propstrongconv}, and
\[ \sigma^h:=\nabla c^\star\left[\,\nabla\left(F^\prime(\rho^h)\right)\,\right].\]
The next lemma shows that $(\sigma^h)_h$ is bounded in $L^q(\Omega_\infty)$, and $\left(\nabla \left(F^\prime(\rho^h)\right)\right)_h$ converges weakly to $\nabla\left(F^\prime(\rho)\right)$ in $L^{q^\star}(\Omega_T)$ for a subsequence..

\begin{lemma}
\label{lemweakconv1}
Assume that $c:\R^d\rightarrow [0,\infty)$ is strictly convex, of class $C^1$ and satisfies $c(0)=0$ and $c(z)\geq \beta\,|\,z\,|^q$ for some $\beta>0$ and $q>1$, and $F:[0,\infty)\rightarrow \R$ is strictly convex, of class $C^2\left((0,\infty)\right)$. If $\rho_0\in{\cal P}_a(\Omega)$ is such that $\rho_0+\frac{1}{\rho_0}\in L^\infty(\Omega)$, then
\begin{equation}
\label{weakconv1}
\|\,\sigma^h\,\|^q_{L^q(\Omega_\infty)} \leq \frac{1}{\beta\,\Big\|\,\frac{1}{\rho_0}\,\Big\|_{L^\infty(\Omega)}}\,\left[\,E_i(\rho_0) - |\,\Omega\,| F\left(\frac{1}{|\,\Omega\,|}\right)\,\right].
\end{equation}
\begin{enumerate}
\item Therefore, there is a subsequence of $\left(\sigma^h\right)_{h\downarrow 0}$, which converges weakly to a function $\sigma$ in $L^q(\Omega_T)$, for $0<T<\infty$.\\

In addition, if $c$ satisfies (HC\ref{c4}), and $F$ satisfies $F(0)=0$ and (HF\ref{f2}), then\\

\item there is a subsequence of $\Big\{\nabla\left(F^\prime(\rho^h)\right)\Big\}_{h\downarrow 0}$, which converges weakly to $\nabla\left(F^\prime(\rho)\right)$, in $L^{q^\star}\left((0,T)\times\Omega\right)$, for $0<T<\infty$.
\end{enumerate}
\end{lemma}
{\tt Proof:} By (\ref{exibc3}), we have that
\begin{equation}
\label{weakconv1.1er}
\frac{S^h_i(y)-y}{h} \,=\, \nabla c^\star\,\left[\,\nabla \left(F^\prime(\rho^h_i(y))\right)\,\right],
\end{equation}
for $i\in\N$, where $S^h_i$ denotes the $c_h$-optimal map that pushes $\rho^h_i$ forward to $\rho^h_{i-1}$. We use (\ref{bba}) and (\ref{weakconv1.1er}), to deduce that
\begin{eqnarray*}
\|\,\sigma^h\,\|^q_{L^q(\Omega_\infty)} &=& \sum_{i=1}^\infty\,h\,\int_\Omega \Big|\,\nabla c^\star\,\left[\,\nabla \left(F^\prime(\rho^h_{i}(y))\right)\,\right]\,\Big|^q\,\mbox{d}y = \sum_{i=1}^\infty h\,\int_\Omega \Big|\,\frac{S^h_i(y)-y}{h}\,\Big|^q\,\mbox{d}y\\
&\leq& \frac{1}{\Big\|\,\frac{1}{\rho_0}\,\Big\|_{L^\infty(\Omega)}}\,\sum_{i=1}^\infty h\,\int_\Omega \Big|\,\frac{S^h_i(y)-y}{h}\,\Big|^q\,\rho^h_{i}(y)\,\mbox{d}y.
\end{eqnarray*}
Since $c(z)\geq \beta\,|\,z\,|^q$, we obtain that
\begin{equation}
\label{weakconv1er}
\|\,\sigma_h\,\|^q_{L^q(\Omega_\infty)} \leq \frac{1}{\beta\,\Big\|\,\frac{1}{\rho_0}\,\Big\|_{L^\infty(\Omega)}}\,\sum_{i=1}^\infty\,h\,W^h_c(\rho^h_{i-1}, \rho^h_i).
\end{equation}
We combine (\ref{exisv1bis}) and (\ref{weakconv1er}), to conclude (\ref{weakconv1}). (\roman{mycounter}) is a direct consequence of (\ref{weakconv1}).\\
Now, fix $0<T<\infty$. By Proposition \ref{propstrongconv}, $(\rho^h)_h$ converges strongly to $\rho$, in $L^1(\Omega_T)$, and by (\ref{strongconv2.1}) and the fact that $F^\prime$ is continuous on $(0,\infty),\;\left(F^\prime(\rho^h)\right)_h$ is bounded in $L^\infty(\Omega_T)$. We deduce that $\left(F^\prime(\rho^h)\right)_h$ converges weakly to $F^\prime(\rho)$ in $L^{q^\star}(\Omega_T)$. And, since $\Big\{\,\nabla\left(F^\prime(\rho^h)\right)\,\Big\}_h$ is bounded in $L^{q^\star}(\Omega_T)$ (because of (\ref{bba}) and (\ref{strongconv2.3})), we conclude (\roman{ii}) \hfill$\Box$\\

\noindent The next lemma extends the energy inequality (\ref{en1}) to the time-space domain $(0,\infty)\times\Omega$.

\begin{lemma}
\label{lemweakconv2}
{\it (Energy inequality in time-space)}\\
Assume that $c:\R^d\rightarrow [0,\infty)$ is strictly convex, of class $C^1$, and satisfies $c(0)=0$ and (HC\ref{c4}), and $F:[0,\infty)\rightarrow \R$ is strictly convex, of class $C^2\left((0,\infty)\right)$, and satisfies $F(0)=0,\;F\in C^2\left((0,\infty)\right)$ and (HF\ref{f2}). If $\rho_0\in{\cal P}_a(\Omega)$ is such that $\rho_0+\frac{1}{\rho_0}\in L^\infty(\Omega)$, and $t\mapsto u(t)$ is a nonnegative function in $C^2_c(\R)$, then
\begin{eqnarray*}
\lefteqn{\int_0^\infty\int_\Omega \left\langle\, \rho^h \nabla \left(F^\prime(\rho^h)\right), \nabla c^\star\,\left[\,\nabla \left(F^\prime(\rho^h)\right)\,\right]\,\right\rangle\, u(t)} \\ 
& &\leq \frac{1}{h}\,\int_0^h\int_\Omega F\left(\rho_0(x)\right)u(t)\,+\, \int_0^\infty\int_\Omega F(\rho^h)\,\partial ^h_t\,u(t),\\ 
\end{eqnarray*}
where 
\[\quad \partial^h_t\,u(t) := \frac{u(t+h)\,-\,u(t)}{h}.\]
\end{lemma}
{\tt Proof:} Let $T$ be such that $\frac{T}{h}\in\N$, and assume that $\mbox{spt}\,u\subset [-T,T]$. We choose $G:=F$ in the energy inequality (\ref{en1}), and we use (\ref{weakconv1.1er}), to obtain that
\begin{eqnarray*} 
\lefteqn{\int_\Omega \frac{F\left(\rho^h_i(y)\right) - F\left(\rho^h_{i-1}(y)\right)}{h}\,\mbox{d}y }\\
& & \leq -\int_\Omega \left\langle\,\nabla \left[\,F^\prime\left(\rho^h_i(y)\right)\,\right], \nabla c^\star\,\left[\,\nabla \left(F^\prime(\rho^h_i(y))\right)\,\right]\,\right\rangle \, \rho^h_i(y)\,\mbox{d}y,
\end{eqnarray*}
for all $i\in\N$. Since $u\geq 0$, we deduce that
\begin{eqnarray}
\label{exibc39}
\lefteqn{ \sum_{i=1}^{T/h}\,\int_{t_{i-1}}^{t_i}\int_\Omega\frac{F\left(\rho^h_i(y)\right) - F\left(\rho^h_{i-1}(y)\right)}{h}\,u(t) }\\
& & \leq -\int_{\Omega_T} \rho^h\,\left\langle\,\nabla \left(F^\prime(\rho^h)\right), \nabla c^\star\,\left[\,\nabla \left(F^\prime(\rho^h)\right)\,\right]\,\right\rangle\,u(t).\nonumber
\end{eqnarray}
By direct computations, the left hand side of the above inequality gives that
\begin{eqnarray*}
\lefteqn{\sum_{i=1}^{T/h}\,\int_{t_{i-1}}^{t_i}\int_\Omega\frac{F\left(\rho^h_i(y)\right) - F\left(\rho^h_{i-1}(y)\right)}{h}\,u(t)}\\
& & = \frac{1}{h}\,\int_{\Omega_T} F\left(\rho^h(t,y)\right)u(t) \,-\, \frac{1}{h}\, \int_{\Omega_h}F\left(\rho_0(y)\right)u(t)\\
& & \quad -\,\frac{1}{h}\,\int_h^T\int_\Omega F\left(\rho^h(t-h)\right)u(t).
\end{eqnarray*}
We use the substitution $\tau=t-h$ in the last integral, and the fact that $u(t+h)=0$ for $t\in(T-h,T)$, to have that
\begin{eqnarray}
\label{weakconv2}
\lefteqn{\sum_{i=1}^{T/h}\,\int_{t_{i-1}}^{t_i}\int_\Omega\frac{F\left(\rho^h_i(y)\right) - F\left(\rho^h_{i-1}(y)\right)}{h}\,u(t)}\\
& &= -\int_{\Omega_T}F\left(\rho^h(t,y)\right)\,\partial^h_t u(t) - \frac{1}{h}\,\int_{\Omega_h}F\left(\rho_0(y)\right)\,u(t).\nonumber
\end{eqnarray}
We combine (\ref{exibc39}) and (\ref{weakconv2}), and we let $T$ go to $\infty$, to complete the proof  \hfill$\Box$ 

\begin{theorem}
\label{theoweakconv}
Assume that $c:\R^d\rightarrow [0,\infty)$ is strictly convex, of class $C^1$, and satisfies $c(0)=0$ and (HC\ref{c4}), and $F:[0,\infty)\rightarrow \R$ is strictly convex, of class $C^2\left((0,\infty)\right)$, and satisfies $F(0)=0$, and (HF\ref{f1}) - (HF\ref{f2}). If $\rho_0\in{\cal P}_a(\Omega)$ is such that $\rho_0+\frac{1}{\rho_0}\in L^\infty(\Omega)$, and $t\mapsto u(t)$ is a nonnegative function in $C^2_c(\R)$, then
\begin{equation}
\label{weakconv6.0}
\lim_{h\downarrow 0}\int_{\Omega_\infty} \left\langle\,\rho^h\sigma^h, \nabla \left(F^\prime(\rho^h)\right)\,\right\rangle\,u(t) = \int_{\Omega_\infty} \left\langle\,\rho\sigma, \nabla \left(F^\prime(\rho)\right)\,\right\rangle\,u(t),
\end{equation}
where $\rho$ and $\sigma$ are defined in Lemma \ref{lemstrongconv1} and Lemma \ref{lemweakconv1}.\\
Therefore, $\left(\mbox{div}(\rho^h\sigma^h)\right)_h$ converges weakly to $\mbox{div}(\rho\sigma)$ for a subsequence in $\left[C^2_c(\R\times\R^d)\right]'$, and  
\begin{equation}
\label{eqn:new1}
\mbox{div}(\rho\sigma)=\mbox{div}\left(\rho\nabla c^\star\left[\nabla\left(F^\prime(\rho)\right)\right]\right).
\end{equation}
\end{theorem}
{\tt Proof:} Let $T>0$ be such that  $\mbox{spt}\,u\subset [-T,T]$, and assume that $\rho(t)=\rho_0$, for $t\leq 0$. Denote by $(\rho^h)_h$ the subsequence of $(\rho^h)_h$, such that
\begin{enumerate}
\item $(\rho^h)_{h\downarrow 0}$ converges to $\rho\,$ a.e.,
\item $\{\,\nabla\left(F^\prime(\rho^h)\right)\,\}_{h\downarrow 0}\,$ converges weakly to $\,\nabla\left(F^\prime(\rho)\right)\,$ in $L^{q^\star}(\Omega_T)$, and
\item $\{\,\sigma^h=\nabla c^\star\left[\nabla\left(F^\prime(\rho^h)\right)\right]\,\}_{h\downarrow 0}\,$ converges weakly to $\,\sigma,\,$ in $L^q(\Omega_T)$,
\end{enumerate}
as in Proposition \ref{propstrongconv} and Lemma \ref{lemweakconv1}. We first observe that
\begin{equation}
\label{theoweak1}
\lim_{h\downarrow 0} \int_{\Omega_T} \langle\sigma^h,\rho^h\nabla\left(F^\prime(\rho)\right)\rangle\,u(t) = \int_{\Omega_T}  \langle\sigma,\rho\nabla\left(F^\prime(\rho)\right)\rangle\,u(t),
\end{equation}
and
\begin{equation}
\label{theoweak2}
\lim_{h\downarrow 0} \int_{\Omega_T} \langle \rho^h\nabla c^\star\left[\,\nabla\left(F^\prime(\rho)\right)\,\right], \nabla\left(F^\prime(\rho^h)\right)-\nabla\left(F^\prime(\rho)\right)\rangle\,u(t) = 0.
\end{equation}
Indeed, since $(\rho^h)_h$ is bounded in $L^\infty(\Omega_T)$ (see (\ref{strongconv2.1})), and $\nabla\left(F^\prime(\rho)\right)\in L^{q^\star}(\Omega_T)$, (\roman{mycounter}) and the dominated convergence theorem imply that $\{\rho^h\nabla\left(F^\prime(\rho)\right)\,\}_{h\downarrow 0}$ converges to $\rho\nabla\left(F^\prime(\rho)\right)$ in $L^{q^\star}(\Omega_T)$. Then, we use (\roman{iii}) and the fact that $u\in C^2_c(\R)$, to conclude (\ref{theoweak1}).\\
Because of Proposition \ref{lemappend3}, the convexity of $c$, and the fact that $c(z)\geq \beta\,|\,z\,|^q$, we have that
\begin{eqnarray*} 
|\,\nabla c^\star(z)\,|^q \leq \frac{c\left(\nabla c^\star(z)\right)}{\beta} &=& \frac{1}{\beta}\,\left(\langle z,\nabla c^\star(z)\rangle - c^\star(z)\right)\\
&\leq& \frac{1}{\beta}\,\langle z,\nabla c^\star(z)\rangle \leq M(\beta,q)\,|\,z\,|^{q^\star}.
\end{eqnarray*}
We deduce that 
\begin{equation}
\label{eqn:new2}
\Big|\,\nabla c^\star\left[\,\nabla\left(F^\prime(\rho)\right)\,\right]\,\Big|^q \leq M(\beta, q)\,\Big|\,\nabla\left(F^\prime(\rho)\right)\,\Big|^{q^\star}, 
\end{equation}
which shows that $\nabla c^\star\left[\,\nabla\left(F^\prime(\rho)\right)\,\right]\in L^q(\Omega_T)$. Then, we use (\roman{mycounter}) and the dominated convergence theorem, to have that $\{\,\rho^h\nabla c^\star\left[\,\nabla\left(F^\prime(\rho)\right)\,\right]\,\}_{h\downarrow 0}$ converges to $\rho\nabla c^\star\left[\,\nabla\left(F^\prime(\rho)\right)\,\right]$ in $L^q(\Omega_T)$. We conclude  (\ref{theoweak2}), because of (\roman{ii}).\\
\noindent The proof of (\ref{weakconv6.0}) follows directly from the following three claims:\\

\noindent {\bf Claim 1.}
\[\int_{\Omega_T} \left\langle\, \rho\sigma, \nabla \left(F^\prime(\rho)\right)\,\right\rangle\,u(t) \leq\, \liminf_{h\downarrow 0}\; \int_{\Omega_T}\left\langle\, \rho^h\sigma^h, \nabla \left(F^\prime(\rho^h)\right)\,\right\rangle\,u(t). \]
{\tt Proof:} Because $c^\star$ is convex, and $u$ and $\rho^h$ are nonnegative, we have that
\[\int_{\Omega_T} \rho^h\left\langle\nabla c^\star\left[\nabla\left(F^\prime(\rho^h)\right)\right]-\nabla c^\star\left[\nabla\left(F^\prime(\rho)\right)\right], \nabla\left(F^\prime(\rho^h)\right)-\nabla\left(F^\prime(\rho)\right)\,\right\rangle u(t) \geq 0,\]
and then,
\begin{eqnarray}
\label{weakconv6.1}
\lefteqn{\liminf_{h\downarrow 0} \int_{\Omega_T}\langle\sigma^h, \rho^h\nabla \left(F^\prime(\rho)\right)\rangle u(t)}\nonumber\\
& & \leq \liminf_{h\downarrow 0} \int_{\Omega_T} \langle \rho^h\sigma^h, \nabla\left(F^\prime(\rho^h)\right)\rangle u(t)\nonumber\\
& & \quad + \limsup_{h\downarrow 0} \int_{\Omega_T} \langle \rho^h\nabla c^\star\left[\nabla\left(F^\prime(\rho)\right)\,\right], \nabla \left(F^\prime(\rho)\right)-\nabla\left(F^\prime(\rho^h)\right)\rangle u(t).
\end{eqnarray}
We combine (\ref{theoweak1}) - (\ref{weakconv6.1}), to conclude Claim 1.\\

\noindent {\bf Claim 2.}
\begin{eqnarray*}
\lefteqn{\limsup_{h\downarrow 0}\int_{\Omega_T} \left\langle \rho^h\sigma^h, \nabla \left(F^\prime(\rho^h)\right)\right\rangle\,u(t)}\\
& & \leq \int_\Omega \left[\,\rho_0F^\prime(\rho_0) - F^\star\left(F^\prime(\rho_0)\right)\,\right]\,u(0) \nonumber \\ 
& & \quad +\;\int_{\Omega_T} \left[\,\rho(t,x)F^\prime\left(\rho(t,x)\right)-F^\star\left(F^\prime\left(\rho(t,x)\right)\right)\,\right]\,u^\prime(t). 
\end{eqnarray*}
{\tt Proof:} First, we observe that
\begin{equation}
\label{theoweak4}
\lim_{h\downarrow 0} \int_{\Omega_T} F(\rho^h)\partial^h_t\,u(t) = \int_{\Omega_T} F(\rho)u^\prime(t).\end{equation}
Indeed, it is clear that
\begin{eqnarray}
\label{theoweak5}
\lefteqn{\Big|\,\int_{\Omega_T} F(\rho^h)\partial^h_tu(t)-F(\rho)u^\prime(t)\,\Big|} \nonumber\\
& & \leq \int_{\Omega_T} |\,F(\rho^h)-F(\rho)\,|\,|u^\prime(t)\,|+\int_{\Omega_T} |\,F(\rho^h)\,|\,|\partial^h_t\,u(t)-u^\prime(t)\,|.
\end{eqnarray}
Because of (\ref{bba}), and the continuity of $F$, we have that $\left(F(\rho^h)\right)_h$ is bounded in $L^\infty(\Omega_T)$. We let $h$ go to $0$ in (\ref{theoweak5}), and we use (\roman{mycounter}), the fact that $u\in C^2_c(\R)$, and the Lebesgue dominated convergence theorem, to conclude (\ref{theoweak4}).\\
\noindent Lemma \ref{lemweakconv2} gives that 
\begin{eqnarray*}
\lefteqn{\limsup_{h\downarrow 0} \int_{\Omega_T} \left\langle \rho^h\sigma^h, \nabla \left(F^\prime(\rho^h)\right)\right\rangle\,u(t)}\\
& &\leq \liminf_{h\downarrow 0} \frac{1}{h}\int_0^h\int_\Omega F(\rho_0)u(t) + \limsup_{h\downarrow 0} \int_{\Omega_T} F(\rho^h)\partial^h_t\,u(t),
\end{eqnarray*}
and by (\ref{theoweak4}) and the continuity of $u$, we deduce that 
\begin{equation}
\label{theoweak4.0}
\limsup_{h\downarrow 0} \int_{\Omega_T} \left\langle \rho^h\sigma^h, \nabla \left(F^\prime(\rho^h)\right)\right\rangle\,u(t) \leq \int_\Omega F(\rho_0)u(0) + \int_{\Omega_T} F\left(\rho(t,x)\right)u^\prime(t).
\end{equation}
Since $F\in C^1\left((0,\infty)\right)$ is strictly convex and satisfies $F(0)=0$ and (HF\ref{f1}), we have that
\begin{equation}
\label{theoweak5.0}
F^\star\left(F^\prime(a)\right)=aF^\prime(a)-F(a),\quad \forall\,a>0.
\end{equation}
% COMMENTS
We substitute (\ref{theoweak5.0}) into (\ref{theoweak4.0}) for $a=\rho(t,x)$ and $a=\rho_0(x)$, to conclude Claim 2.\\

\noindent{\bf Claim 3.} 
\begin{eqnarray*}
\lefteqn{\int_\Omega \left[\,\rho_0 F^\prime(\rho_0)-F^\star\left(F^\prime(\rho_0)\right)\,\right]\,u(0)}\\
& & \quad +\int_{\Omega_T} \left[\,\rho(t,x)F^\prime\left(\rho(t,x)\right)-F^\star\left(F^\prime\left(\rho(t,x)\right)\right)\,\right]\,u^\prime(t) \\ 
& & \leq \int_{\Omega_T} \langle \rho\sigma, \nabla \left(F^\prime(\rho)\right)\rangle\,u(t).
\end{eqnarray*}
{\tt Proof:} Set $\xi(t,x):= F^\prime\left(\rho(t,x)\right)u(t)$ for $(t,x)\in \R\times\Omega$. Because of (\roman{mycounter}), (\roman{ii}), (\ref{bba}), and the fact that $F\in C^2\left((0,\infty)\right)$, we have that $F^\prime(\rho)\in L^\infty(\Omega_T)$ and $\nabla\left(F^\prime(\rho)\right)\in L^{q^\star}(\Omega_T)$. We approximate $F^\prime(\rho)$ by $C^\infty(\Omega_T)$-functions in $W^{1,q^\star}(\Omega_T)$, and we use (\ref{exibc25}) with the backward derivative  $\partial^{-h}_t\,\xi(t,x):= \frac{\xi(t,x)-\xi(t-h,x)}{h}$, and Proposition \ref{proplimsecond}, to obtain that
\[\int_{\Omega_T} (\rho_0-\rho^h)\partial^{-h}_t\,\xi + \int_{\Omega_T} \langle \sigma^h, \rho^h\nabla \left(F^\prime(\rho)\right)\rangle\,u(t) = 0(h^{\epsilon(q)}),\]
where $\epsilon(q)=\min(1,q-1)$. We let $h$ go to $0$ in the subsequent equality, and we use (\ref{theoweak1}), to conclude that
\begin{equation}
\label{exibc44}
\lim_{h\downarrow 0} \int_{\Omega_T}(\rho_0-\rho^h)\,\partial^{-h}_t\,\xi +\int_{\Omega_T}\langle \sigma, \rho\nabla\left(F^\prime(\rho)\right)\rangle\,u(t)=0.
\end{equation}
Since $\mbox{spt}\,u\subset [-T,T]$, we have that
\[\int_{\Omega_T}\rho_0\,\partial^{-h}_t\xi = -\frac{1}{h}\int_{-h}^0\int_\Omega \rho_0(x)\xi(t,x),\]
 and then,
\begin{equation}
\label{weakconv6.6}
\lim_{h\downarrow 0} \int_{\Omega_T} \rho_0\,\partial^{-h}_t\,\xi = -\int_\Omega \rho_0(x)\,\xi(0,x) = -\int_\Omega \rho_0 F^\prime(\rho_0) u(0).
\end{equation}
We combine (\ref{exibc44}), (\ref{weakconv6.6}) and (\roman{mycounter}), to have that
\begin{equation}
\label{weakconv6.7}
\int_{\Omega_T} \langle\sigma, \rho\nabla\left(F^\prime(\rho)\right)\rangle\, u(t) = \lim_{h\downarrow 0} \int_{\Omega_T} \rho(t,x)\partial^{-h}_t \xi(t,x) + \int_\Omega \rho_0 F^\prime(\rho_0) u(0).
\end{equation}
By direct computations, we obtain that
\begin{eqnarray*}
 \rho(t,x)\partial^{-h}_t\xi(t,x) &=& \rho(t,x)\,F^\prime\left(\rho(t,x)\right)\partial^{-h}_t\,u(t)\\
& &\quad +\, \frac{1}{h}\,\rho(t,x)\,u(t-h)\left[\,F^\prime\left(\rho(t,x)\right)-F^\prime\left(\rho(t-h,x)\right)\,\right].
\end{eqnarray*}
Since $F\in C^1\left((0,\infty)\right)$ is strictly convex, and satisfies $F(0)=0$ and (HF\ref{f1}),  we have that
\[\left(F^\prime(b)-F^\prime(a)\right)b \geq F^\star\left(F^\prime(b)\right)-F^\star\left(F^\prime(a)\right), \quad \forall\,a,b>0,\]
% COMMENTS
and then, we deduce that
\begin{eqnarray*}
\rho(t,x)\,\partial^{-h}_t\xi(t,x) &\geq& \rho(t,x)\,F^\prime\left(\rho(t,x)\right)\,\partial^{-h}_t u(t)\\
& & + \frac{1}{h}\,u(t-h)\left[\,F^\star\left(F^\prime\left(\rho(t,x)\right)\right) - F^\star\left(F^\prime\left(\rho(t-h,x)\right)\right)\,\right].
\end{eqnarray*}
We integrate both sides of the subsequent inequality over $\Omega_T$, and we use that $u=0$ on $(T-h,T)$ for $h$ small enough, and $\rho(t,x)=\rho_0(x)$ for $t\in(-h,0)$, to obtain that
\begin{eqnarray*}
\int_{\Omega_T}  \rho(t,x)\partial^{-h}_t\xi(t,x) &\geq& \int_{\Omega_T} \left[\,\rho(t,x)F^\prime\left(\rho(t,x)\right) - F^\star\left(F^\prime\left(\rho(t,x)\right)\right)\,\right]\partial^{-h}_t u(t)\\
& & \quad - \frac{1}{h}\int_0^hu(t-h)\int_\Omega F^\star\left(F^\prime(\rho_0)\right). 
\end{eqnarray*}
We let $h$ go to $0$ in the above inequality, to deduce that
\begin{eqnarray}
\label{weakconv6.8}
\lim_{h\downarrow 0} \int_{\Omega_T} \rho(t,x)\partial^{-h}_t\xi(t,x) &\geq&  \int_{\Omega_T} \left[\,\rho(t,x)F^\prime\left(\rho(t,x)\right) -  F^\star\left(F^\prime\left(\rho(t,x)\right)\right)\,\right]u^\prime(t)\nonumber\\
& &- \int_\Omega F^\star\left(F^\prime(\rho_0)\right)u(0).
\end{eqnarray}
We combine (\ref{weakconv6.7}) and (\ref{weakconv6.8}), to conclude Claim 3.\\
In the end, we show that $\sigma = \nabla c^\star\left[\,\nabla\left(F^\prime(\rho)\right)\,\right]$, which combined with Lemma \ref{lemweakconv1} completes the proof of Theorem \ref{theoweakconv}. Indeed, let $\epsilon>0,\;\psi\in C^\infty(\Omega)$, and set $\omega_\epsilon(t,x):= F^\prime\left(\rho(t,x)\right) - \epsilon\psi(x)$. It is clear that $\nabla \omega_\epsilon\in L^{q^\star}(\Omega_T)$, and 
\[ \Big|\,\nabla c^\star\left(\nabla \omega_\epsilon\right)\,\Big|^q\leq M(\beta,q)\,|\,\nabla \omega_\epsilon\,|^{q^\star},\]
as in the proof of (\ref{theoweak2}). We deduce that $\nabla c^\star\left(\nabla \omega_\epsilon\right)\in L^q(\Omega_T)$. We use that $c^\star$ is convex, and $\rho^h$ and $u$ are nonnegative, to have that 
\[\int_{\Omega_T} \rho^h\,\langle\nabla c^\star\,\left[\,\nabla\left(F^\prime(\rho^h)\right)\,\right]-\nabla c^\star\left(\nabla w_\epsilon\right),\, \nabla\left(F^\prime(\rho^h)\right)-\nabla w_\epsilon\rangle\,u(t) \,\geq\, 0. \] 
We let $h$ go to $0$ in the above inequality, to obtain that
\begin{eqnarray}
\label{theoweak11}
\limsup_{h\downarrow 0}\int_{\Omega_T} \langle \rho^h\sigma^h, \nabla\left(F^\prime(\rho^h)\right)\rangle\,u(t) - \liminf_{h\downarrow 0}\int_{\Omega_T} \langle \sigma^h, \rho^h\nabla w_\epsilon\rangle\,u(t)\nonumber\\
\qquad- \liminf_{h\downarrow 0} \int_{\Omega_T}\langle \rho^h \nabla c^\star(\nabla w_\epsilon),\nabla\left(F^\prime(\rho^h)\right)-\nabla w_\epsilon \rangle\,u(t) \geq 0.
\end{eqnarray}
As in the proof of (\ref{theoweak1}) and (\ref{theoweak2}), we have that
\begin{equation}
\label{theoweak12.0}
\liminf_{h\downarrow 0} \int_{\Omega_T} \langle\sigma^h, \rho^h\nabla w_\epsilon\rangle\,u(t) = \int_{\Omega_T} \langle\sigma, \rho\nabla w_\epsilon\rangle\,u(t),
\end{equation}
and
\begin{eqnarray}
\label{theoweak12}
\lefteqn{\liminf_{h\downarrow 0} \int_{\Omega_T} \langle \rho^h \nabla c^\star(\nabla w_\epsilon),\nabla\left(F^\prime(\rho^h)\right)-\nabla w_\epsilon \rangle\,u(t)}\nonumber\\
& & = \int_{\Omega_T} \langle \rho\nabla c^\star(\nabla w_\epsilon),\nabla\left(F^\prime(\rho)\right)-\nabla w_\epsilon \rangle\,u(t).
\end{eqnarray}
We combine (\ref{weakconv6.0}) and (\ref{theoweak11}) - (\ref{theoweak12}), to have that
\[\int_{\Omega_T} \langle \rho\sigma-\rho\nabla c^\star(\nabla w_\epsilon), \nabla\left(F^\prime(\rho)\right)-\nabla w_\epsilon\rangle\,u(t) \geq 0.\] 
We divide the subsequent inequality by $\epsilon$, and we let $\epsilon$ go to $0$, to obtain that
\[ \int_{\Omega_T}\langle \rho\sigma - \rho\nabla c^\star\,\left[\,\nabla\left(F^\prime(\rho)\right)\,\right],\nabla\psi(x)\,u(t)\rangle \geq 0. \] 
Choosing $-\psi$ in place of $\psi$, we get that
\[ \int_{\Omega_T}\langle \rho\sigma - \rho\,\nabla c^\star\,\left[\,\nabla\left(F^\prime(\rho)\right)\,\right], \,\nabla\psi(x)\,u(t)\rangle = 0. \] 
And since $\psi$ and $u\geq 0$ are arbitrary test functions, we deduce (\ref{eqn:new1}). This completes the proof of the theorem  \hfill$\Box$

\subsection{Existence and uniqueness of solutions}
Here, we state and prove the theorem of existence and uniqueness for (\ref{intro8}).

\begin{theorem}
\label{theoexistbound}
{\it (Case $\,V=0\,$).}\\
Assume that $c:\R^d\rightarrow [0,\infty)$ is strictly convex, of class $C^1$, and satisfies $c(0)=0$ and (HC\ref{c4}), and $F:[0,\infty)\rightarrow \R$ is strictly convex, of class $C^2\left((0,\infty)\right)$, and satisfies $F(0)=0$ and (HF\ref{f1}) - (HF\ref{f2}). If $\rho_0\in{\cal P}_a(\Omega)$ is such that $\rho_0+\frac{1}{\rho_0}\in L^\infty(\Omega)$, and $V=0$, then, (\ref{intro8}) has a unique weak solution  $\rho:[0,\infty)\times\Omega \rightarrow [0,\infty)$, in the sense that
\begin{enumerate}
\item $\rho+\frac{1}{\rho}\in L^\infty\left((0,\infty);L^\infty(\Omega)\right),\; \nabla\left(F^\prime(\rho)\right)\in L^{q^\star}(\Omega_T)$ for $0<T<\infty$, and
\item for $\xi\in C^2_c(\R\times\R^d)$,
\begin{equation}
\label{theobound0}
\int_{\Omega_\infty} \Big\{-\rho\frac{\partial\,\xi}{\partial\,t} + \left\langle \rho\,\nabla c^\star\,\left[\,\nabla\left(F^\prime(\rho)\right)\,\right], \nabla \xi\right\rangle \Big\}=\int_\Omega \rho_0(x)\,\xi(0,x)\,\mbox{d}x.
\end{equation}
\end{enumerate}  
\end{theorem}
{\tt Proof:} Proposition \ref{propstrongconv} gives that $(\rho^h)_h$ converges to $\rho$ a.e. for a subsequence, and since $\rho^h\geq 0$ for all $h$, we deduce that $\rho\geq 0$. We combine (\ref{bba}) and Proposition \ref{propstrongconv}, to have that $\rho+\frac{1}{\rho}\in L^\infty\left((0,\infty);L^\infty(\Omega)\right)$. We use that $\nabla\left(F^\prime(\rho)\right)\in L^{q^\star}(\Omega_T)$ (Lemma \ref{lemweakconv1}) to conclude $(\roman{mycounter})$.\\
\noindent Recall that (\ref{eqn:new2}) gives that $\nabla c^\star\left[\,\nabla\left(F^\prime(\rho)\right)\,\right]\in L^q(\Omega_T)$ for $0<T<\infty$, and (\ref{bba}) and Proposition \ref{propstrongconv} imply that $\rho\in L^\infty(\Omega_T)$. We deduce that $\rho\nabla c^\star\left[\,\nabla\left(F^\prime(\rho)\right)\,\right]\in L^q(\Omega_T)$. Now, fix $0<T<\infty$, and let $\xi\in C^2_c(\R\times\R^d)$ be such that $\mbox{spt}\,\xi(.,x)\subset [-T,T]$ for $x\in\Omega$. Because of  Proposition \ref{propinterpolation} and Proposition \ref{proplimsecond}, we have that
\begin{equation}
\label{theobound1}
\lim_{h\downarrow 0}\int_{\Omega_T} \Big\{(\rho_0-\rho^h)\,\partial^h_t\xi + \langle\,\rho^h \nabla c^\star\left[\,\nabla\left(F^{\prime}(\rho^h)\right)\,\right],\nabla \xi\,\rangle\Big\} = 0.
\end{equation}
Lemma \ref{lemstrongconv1} gives that $(\rho^h)_h$ converges weakly to $\rho$ in $L^1(\Omega_T)$ for a subsequence, and then, we have that
\begin{equation}
\label{theobound3}
\lim_{h\downarrow 0}\int_{\Omega_T}(\rho_0-\rho^h)\,\partial^h_t\xi = \int_{\Omega_T}(\rho_0-\rho)\,\frac{\partial \xi}{\partial t}= -\left[\,\int_{\Omega_T}\rho\,\frac{\partial \xi}{\partial t} + \int_\Omega \rho_0(x)\xi(0,x)\,\right]. 
\end{equation}
From Theorem \ref{theoweakconv}, we have $\left(\mbox{div}\{\rho^h\nabla c^\star\left(\nabla\left(F^\prime(\rho^h)\right)\right)\}\right)_h$ converges weakly to\\
 $\mbox{div}\{\rho\nabla c^\star\left(\nabla\left(F^\prime(\rho)\right)\right)\}$ in $\left[C^2_c(\R\times\R^d)\right]'$, for a subsequence, then we deduce that
\begin{equation}
\label{eqn:new3}
 \lim_{h\downarrow 0} \int_{\Omega_T}\langle\rho^h\nabla c^\star\left[\nabla\left(F^\prime(\rho)\right)\right],\nabla\xi\rangle = \int_{\Omega_T}\langle \rho\nabla c^\star\left[\nabla\left(F^\prime(\rho)\right)\right],\nabla\xi\rangle.
\end{equation}
We combine (\ref{theobound1}) - (\ref{eqn:new3}), and we use the fact that $\mbox{spt}\,\xi(.,x)\subset [-T,T]$, to conclude (\ref{theobound0}).\\
Here, we prove uniqueness of solutions to (\ref{intro8}) when $\frac{\partial \rho}{\partial t}\in L^1\left((0,T)\times\Omega\right)$, for $0<T<\infty$. Using the arguments in \cite{otto:contraction}, it is easy to extend the proof to the general case. In fact, assumption (\ref{intro2}) imposed in \cite{luckhaus:quasilinear} would not be required here. The convexity of $c^\star$, that is, $\langle \nabla c^\star(z_1)-\nabla c^\star(z_2), z_1-z_2\rangle \geq 0$, for $z_1,\,z_2\in\R^d$, suffices to extend the proof. \\
Let $T>0$, and assume that $\rho_1$ and $\rho_2$ are weak solutions of (\ref{intro8}) with the same initial data, such that $N\leq \rho_j\leq M$ a.e., and $\frac{\partial \rho_j}{\partial t}\in L^1(\Omega_T),\;j=1,2$. Since $\nabla\left(F^\prime(\rho_j)\right)\in L^{q^\star}(\Omega_T)$, and
\[ \Big|\,\nabla c^\star\left[\,\nabla\left(F^\prime(\rho_j)\right)\,\right]\,\Big|^q\leq M(\beta,q)\,\Big|\,\nabla\left(F^\prime(\rho_j)\right)\,\Big|^{q^\star},\]
we have that $\nabla c^\star\left[\,\nabla\left(F^\prime(\rho_j)\right)\,\right]\in L^q(\Omega_T)$. For $\delta>0$, we define
\[ \Omega_T\ni (t,x)\mapsto \xi_\delta\,(t,x) := \varphi_\delta\left(F^\prime\left(\rho_1(t,x)\right)-F^\prime\left(\rho_2(t,x)\right)\right),\]
where 
\[\varphi_\delta(\tau):=\left\{\begin{array}{lcl}
0 &\mbox{if}& \tau\leq 0\\
\frac{\tau}{\delta} &\mbox{if}& 0\leq \tau\leq \delta\\
1 &\mbox{if}& \tau\geq \delta.
\end{array}\right.\]
Using a smooth approximation of $\xi_\delta$ as a test function in the differential equations satisfied by $\rho_1$ and $\rho_2$, and passing to the limit, we have that 
\[\int_{\Omega_T} \xi_\delta\,\partial_t(\rho_1-\rho_2) = -\int_{\Omega_T} \langle \rho_1\nabla c^\star\left[\,\nabla\left(F^\prime(\rho_1)\right)\,\right]-\rho_2\nabla c^\star\left[\,\nabla\left(F^\prime(\rho_2)\right)\,\right], \nabla\xi_\delta\,\rangle,\]
which reads as
\begin{eqnarray*}
\lefteqn{\int_{\Omega_T} \xi_\delta\,\partial_t(\rho_1-\rho_2)}\\
& & = -\frac{1}{\delta}\int_{\Omega_{T,\delta}^{(1,2)}} \rho_1\langle\nabla c^\star\left[\,\nabla\left(F^\prime(\rho_1)\right)\,\right]-\nabla c^\star\left[\,\nabla\left(F^\prime(\rho_2)\right)\,\right], \nabla\left(F^\prime(\rho_1)-F^\prime(\rho_2)\right)\rangle\\
& &\quad -\frac{1}{\delta} \int_{\Omega_{T,\delta}^{(1,2)}} (\rho_1-\rho_2)\langle \nabla c^\star\left[\,\nabla\left(F^\prime(\rho_2)\right)\,\right], \nabla\left(F^\prime(\rho_1)-F^\prime(\rho_2)\right)\rangle,
\end{eqnarray*}
where $\Omega_{T,\delta}^{(1,2)}:=\Omega_T\cap [\,0<F^\prime(\rho_1)-F^\prime(\rho_2)<\delta\,]$.
Because $c^\star$ is convex, the first term on the right hand side of the above equality is nonpositive. And since $F\in C^1\left((0,\infty)\right)$ is strictly convex and satisfies (HF\ref{f1}), and $N\leq \rho_1,\,\rho_2\leq M$ a.e., we have, a.e. on $\Omega_{T,\delta}^{(1,2)}$, that
\[ |\,\rho_1-\rho_2\,| = \Big|\left[(F^\star)^\prime\circ F^\prime\right](\rho_1) - \left[(F^\star)^\prime\circ F^\prime(\rho_2)\right]\Big| \leq\,\delta\,\sup_{\tau\in \left[\,F^\prime(N),F^\prime(M)\,\right]}\, (F^\star)^{\prime\prime}(\tau).\]
We deduce that
\begin{eqnarray*}
\lefteqn{\int_{\Omega_T} \xi_\delta\,\partial_t(\rho_1-\rho_2)}\\
& & \leq \sup_{\tau\in \left[\,F^\prime(N),F^\prime(M)\,\right]}\,(F^\star)^{\prime\prime}(\tau) \int_{\Omega_{T,\delta}^{(1,2)}} \Big|\,\langle \nabla c^\star\left[\,\nabla\left(F^\prime(\rho_2)\right)\,\right], \nabla\left(F^\prime(\rho_1)-F^\prime(\rho_2)\right)\rangle\,\Big|.
\end{eqnarray*}
We let $\delta$ go to $0$ in the subsequent inequality, and we use that $\varphi_\delta\rightarrow \1_{[0,\infty)},\;|\,\Omega_{T,\delta}^{(1,2)}\,|\rightarrow 0$, and $\left[\,F^\prime(\rho_1)-F^\prime(\rho_2)\geq 0\,\right]=\left[\,\rho_1-\rho_2\geq 0\,\right]$, to have that 
\[\int_{\Omega_T} \partial_t\left[\,(\rho_1-\rho_2)^+\,\right] \leq 0,\]
which reads as  
\[ \int_\Omega\left[\,\rho_1(T)-\rho_2(T)\,\right]^+ \leq \int_\Omega \left[\,\rho_1(0)-\rho_2(0)\,\right]^+ = 0,\]
for $0<T<\infty$. Interchanging $\rho_1$ and $\rho_2$ in the above argument, we conclude that $\rho_1=\rho_2$ \hfill$\Box$\\

\begin{theorem}
\label{theoexistbound1}
{\it (General case).}\\
Assume that $\,V:\bar{\Omega}\rightarrow [0,\infty)\,$ is convex, of class $C^1,\,$ and $\,c:\R^d\rightarrow [0,\infty)\,$ is strictly convex, of class $C^1$, and satisfies $\,c(0)=0\,$ and (HC\ref{c4}). Assume that $\,F:[0,\infty)\rightarrow \R\,$ is strictly convex, of class $C^2\left((0,\infty)\right),\,$ and satisfies $\,F(0)=0\,$ and (HF\ref{f1}) - (HF\ref{f2}). If $\,\rho_0\in{\cal P}_a(\Omega)\,$ is such that $\,\rho_0+\frac{1}{\rho_0}\in L^\infty(\Omega),\,$ then, (\ref{intro8}) has a unique weak solution  $\,\rho:[0,\infty)\times\Omega \rightarrow [0,\infty),\,$ in the sense that
\begin{enumerate}
\item $\rho+\frac{1}{\rho}\in L^\infty\left((0,\infty);L^\infty(\Omega)\right),\;\nabla\left(F^\prime(\rho)\right)\in L^{q^\star}(\Omega_T)$ for $0<T<\infty$, and
\item for $\,\xi\in C^2_c(\R\times\R^d)$,
\begin{equation}
\label{theobound1.0}
\int_{\Omega_\infty} \Big\{-\rho\frac{\partial\,\xi}{\partial\,t} + \left\langle \rho\,\nabla c^\star\,\left[\,\nabla\left(F^\prime(\rho)+V\right)\,\right], \nabla \xi\right\rangle \Big\} =\int_\Omega \rho_0(x)\,\xi(0,x)\,\mbox{d}x.
\end{equation}
\end{enumerate}
\end{theorem}
{\tt Proof:} The proof of the uniqueness of solutions is similar to that of Theorem \ref{theoexistbound}. Here, we only prove existence of solutions to (\ref{intro8}). Let $\xi\in C^2_c(\R\times\R^d)$ be such that $\mbox{spt}\,\xi(\cdot,x)\subset [-T,T]$ for $x\in\Omega$ and for some $0<T<\infty$. Following the arguments in the previous sections, where the minimization problem $(P)$ (\ref{convana2}) is replaced by 
\[ (P^V):\quad \inf\,\{\,hW^h_c(\rho_0,\rho)+E(\rho):\quad \rho\in{\cal P}_a(\Omega)\,\},\]
and
\[ E(\rho):=E_i(\rho)+\int_\Omega \rho\,V\,\mbox{d}x,\]
we have, as in Proposition \ref{propinterpolation}, that
\begin{eqnarray}
\label{existV1}
\Big|\int _{\Omega_T} (\rho_0-\rho^h)\,\partial^h_t \xi\,\mbox{d}x\,\mbox{d}t + \int_{\Omega_T} \langle\,\rho^h \nabla c^\star\left[\,\nabla\left(F^{\prime}(\rho^h)+V\right)\,\right],\nabla \xi \,\rangle \,\mbox{d}x\,\mbox{d}t\,\Big| \nonumber\\ 
\leq \frac{1}{2}\,\sup_{[0,T]\times\bar{\Omega}}\,\Big|D^2 \xi(t,x)\Big|\; \sum _{i=1}^{T/h} \int _{\Omega \times \Omega}| x-y|^2\, d\gamma_i^h (x,y),
\end{eqnarray}
and, as in Proposition \ref{proplimsecond}, that
\begin{equation}
\label{existV2}
\sum_{i=1}^{T/h} \int_{\Omega\times\Omega}|\,x-y\,|^2\,\mbox{d}\gamma_i^h(x,y) \leq  M(\Omega,T,F,\rho_0,q,\beta)\,h^{\epsilon(q)}.
\end{equation}
We let $h$ go to $0$ in (\ref{existV1}), and we use (\ref{existV2}), to deduce that
\begin{equation}
\label{theobound1.1}
\lim_{h\downarrow 0}\int_{\Omega_T} \Big\{(\rho_0-\rho^h)\,\partial^h_t\xi + \langle\,\rho^h \nabla c^\star\left[\,\nabla\left(F^{\prime}(\rho^h)+V\right)\,\right],\nabla \xi\,\rangle\Big\} = 0.
\end{equation}
The following claim suffices to conclude Theorem \ref{theoexistbound1}.\\

\noindent{\bf Claim.}\hspace{2mm} For $\,0<T<\infty,\,$ the estimates
\begin{equation}
\label{theobound1.2}
\|\,\rho^h\,\|_{L^\infty\left((0,\infty); L^\infty(\Omega)\right)} \leq \|\,\rho_0\,\|_{L^\infty(\Omega)},
\end{equation}

\begin{equation}
\label{theobound1.3}
\int_{\Omega_T} \rho^h\Big|\,\nabla\left(F^\prime(\rho^h)\right)\,\Big|^{q^\star} \leq M(\Omega, T, F, \rho_0, V, q, \alpha),
\end{equation}
and the energy inequality in time-space 
\begin{eqnarray}
\label{theobound1.4}
\lefteqn{\int_{\Omega_\infty} \langle \rho^h\nabla\left(F^\prime(\rho^h)+V\right), \nabla c^\star\left[\,\nabla\left(F^\prime(\rho^h)+V\right)\,\right]\rangle\,u(t)}\\
& &\leq \frac{1}{h}\,\int_{\Omega_h}\left[\,F(\rho_0)+\rho_0V\,\right]\,u(t) + \int_{\Omega_\infty} \left[\,F(\rho^h)+\rho^hV\,\right]\,\partial^h_t\,u(t),\nonumber
\end{eqnarray}
hold, for nonnegative functions $u$ in $C^2_c(\R)$.\\

\noindent Indeed, because of (\ref{theobound1.2}), there exists $\rho:[0,\infty)\times\Omega \rightarrow [0,\infty)$, such that \\
$(iii).$\hspace{.2cm} $(\rho^h)_h\,$ converges to $\rho,\,$ weakly in $L^1(\Omega_T),\,$ for a subsequence.\\
As a consequence, 
\begin{equation}
\label{theobound1.5}
\lim_{h\downarrow 0} \int_{\Omega_T} (\rho_0-\rho^h)\,\partial^h_t\,\xi = \int_{\Omega_T} (\rho_0-\rho)\,\frac{\partial \xi}{\partial t}.
\end{equation}
Using (\ref{theobound1.2}) and (\ref{theobound1.3}), we deduce the space-compactness and the time-compactness of $(\rho^h)_h$ in $L^1(\Omega_T)$, as in the case where $V=0$. Hence,\\
$(iv).$\hspace{.2cm} $(\rho^h)_h\,$ converges strongly to $\rho$ in $L^1(\Omega_T),\,$ for a subsequence.\\
Then, we use $(iv)$, (\ref{theobound1.4}), and we follow the lines of the proof of Theorem \ref{theoweakconv}, where we use $F^\prime(\rho^h)+V$ in place of $F^\prime(\rho^h)$, and $F(\rho^h)+\rho^hV$ in place of $F(\rho^h)$, to conclude that\\
$(v).$\hspace{.2cm} $\left(\mbox{div}\{\rho^h\nabla c^\star\left[\nabla\left(F^\prime(\rho^h)+V\right)\right]\}\right)_h\,$ converges weakly to $\mbox{div}\{\rho\nabla c^\star\left[\nabla\left(F^\prime(\rho)+V\right)\right]\},\,$ in $\left[C^2_c(\R\times\R^d)\right]'$, for a subsequence. \\
Hence,
\begin{equation}
\label{theobound1.6}
\lim_{h\downarrow 0} \int_{\Omega_T} \langle \rho^h\nabla c^\star\left[\,\nabla\left(F^\prime(\rho^h)+V\right)\,\right],\nabla\xi\rangle = \int_{\Omega_T} \langle \rho\nabla c^\star\left[\,\nabla\left(F^\prime(\rho)+V\right)\,\right],\nabla\xi\rangle.
\end{equation}
We combine (\ref{theobound1.1}), (\ref{theobound1.5}), and (\ref{theobound1.6}), to conclude (\ref{theobound1.0}).\\
As in Theorem \ref{theoexistbound}, $(\roman{mycounter})$ follows directly from (\ref{theobound1.2}), (\ref{theobound1.3}), and the maximum/minimum principle of Proposition \ref{propbound} for $\nabla V\not =0$.\\

\noindent{\tt Proof of the Claim:} (\ref{theobound1.2}) is a direct consequence of the maximum principle of Proposition \ref{propbound} for $\nabla V\not= 0$.\\
As in the case $V=0$, we have, because of Proposition \ref{propproperty} and the maximum/minimum principle of Proposition \ref{propbound}, that $P(\rho^h_i)\in W^{1,\infty}(\Omega)$, and $\nabla\left(F^\prime(\rho^h_i)\right)\in L^\infty(\Omega)$. Then, choosing $G:=F$ in Theorem \ref{theoenergyrel}, the (internal) energy inequality (\ref{en1}) read as
\[\int_\Omega F(\rho^h_{i-1}) - \int_\Omega F(\rho^h_i)  \geq \int_\Omega \langle\nabla\left(F^\prime(\rho^h_i)\right), S^h_i(y)-y\rangle\, \rho^h_i(y)\,\mbox{d}y,\]
where $S^h_i$ is the $c_h$-optimal map that pushes $\rho^h_i$ forward to $\rho^h_{i-1}$. We use that $(S^h_i)_{\#}\rho^h_i=\rho^h_{i-1}$, and $V\in C^1(\Omega)$ is convex, to deduce the {\it potential energy inequality} 
\[\int_\Omega \rho^h_{i-1}\,V - \int_\Omega \rho^h_i\,V \geq \int_\Omega \langle\,\nabla V , S^h_i(y)-y\,\rangle\, \rho^h_i(y)\,\mbox{d}y.\]
We add both of the subsequent inequalities, and we use the Euler-Lagrange equation of $(P^V)$, that is,
\begin{equation}
\label{eulerV}
\frac{S^h_1(y)-y}{h} \,=\, \nabla c^\star \left[ \nabla (F^\prime (\rho^h_1(y))+V(y))\right], \quad \mbox{for a.e.}\;\;y\in\Omega,
\end{equation}
(where $S^h_1$ is the $c_h$-optimal map that pushes $\rho^h_1$ forward to $\rho_0$), to deduce the {\it free energy inequality}
\begin{equation}
\label{en}
E(\rho^h_{i-1}) - E(\rho^h_i)\geq h\int_{\Omega_T} \langle \nabla\left(F^\prime(\rho^h_i)+V\right), \nabla c^\star\left[\,\nabla\left(F^\prime(\rho^h_i)+V\right)\,\right]\rangle\,\rho^h_i,
\end{equation}
for $i\in\N$. We sum (\ref{en}) over $i$, and we use that $V$ and $\rho^h_{T/h}$ are nonnegative, and Jensen's inequality, to have that
\[h\int_{\Omega_T} \langle \nabla\left(F^\prime(\rho^h)+V\right), \nabla c^\star\left[\,\nabla\left(F^\prime(\rho^h)+V\right)\,\right]\rangle\,\rho^h\leq E(\rho_0) - |\,\Omega\,|\,F\left(\frac{1}{|\,\Omega\,|}\right).\]
We conclude, as in the proof of (\ref{strongconv2.3}), that 
\begin{equation}
\label{theobound1.7}
\int_{\Omega_T} \rho^h\,\Big|\,\nabla\left(F^\prime(\rho^h)+V\right)\,\Big|^{q^\star} \leq \overline{M}(\Omega, T, F, \rho_0, q, \alpha).
\end{equation}
On the other hand, because of (\ref{theobound1.2}) and the fact that $V\in C^1(\bar{\Omega})$, we have that 
\begin{equation}
\label{theobound1.8}
\Big\|\,(\rho^h)^{1/q^\star}\,\nabla V\,\Big\|_{L^{q^\star}(\Omega_T)} \leq \|\,\rho_0\,\|^{1/q^\star}_{L^\infty(\Omega)}\,\|\,\nabla V\,\|_{L^\infty(\Omega)}.
\end{equation}
We combine (\ref{theobound1.7}) and (\ref{theobound1.8}), to conclude (\ref{theobound1.3}).\\
The proof of (\ref{theobound1.4}) follows the lines of the proof of Lemma \ref{lemweakconv2} where we use the free energy inequality (\ref{en}) in place of the internal energy inequality (\ref{en1}) \hfill$\Box$

\begin{rem}
\label{remextension}
(Existence of solutions to (\ref{intro8}) for a wider class of $\rho_0$).\\
Here, we show how to extend our existence theorem \ref{theoexistbound1} to a wider class of initial probability densities $\rho_0:\;\rho_0\in L^\infty(\Omega)$ and $\frac{1}{\rho_0}\not\in L^\infty(\Omega)$, or $\rho_0\in L^p(\Omega),\;p\geq q$. For simplicity, we assume that $V=0$.
\end{rem}
{\bf Case 1}: $\rho_0\in L^\infty(\Omega)$, and $\frac{1}{\rho_0}\not\in L^\infty(\Omega)$.\\
Let $(\rho_{0,\delta})_\delta$ be a sequence in ${\cal P}_a(\Omega)$, such that 
\begin{equation}
\label{rem1}
\left\{\begin{array}{l}
\eta_\delta \leq \rho_{0,\delta}\leq \|\,\rho_0\,\|_{L^\infty(\Omega)}\,\mbox{a.e.},\;\mbox{where}\;\; 0<\eta_\delta\leq \delta\\ \\
E_i(\rho_{0,\delta})\leq E_i(\rho_0)\\ \\
\rho_{0,\delta}\rightarrow \rho_0\;\mbox{in}\; L^1(\Omega),\;\mbox{as}\;\delta\downarrow 0
\end{array}\right.
\end{equation}
(see \cite{agueh:thesis}, Proposition 1.4.2), and define the approximate solution $\rho^h_\delta$ to (\ref{intro8}), by
\[\rho^h_\delta:=\left\{\begin{array}{lcl}
\rho_{0,\delta} &\mbox{if}& t=0\\ \\
\rho^h_{i,\delta} &\mbox{if}& t\in ((i-1)h,ih],
\end{array}\right.\]
where $\rho^h_{i,\delta}$ is the unique minimizer of 
\[(P_{i,\delta}):\quad \inf\,\{\,hW^h_c\left(\rho^h_{i-1,\delta},\rho\right)+ E_i(\rho):\;\rho\in {\cal P}_a(\Omega)\}.\]
Since $\rho_{0,\delta}+\frac{1}{\rho_{0,\delta}}\in L^\infty(\Omega)$, we have as before, that 
\begin{equation}
\label{rem2}
\left\{\begin{array}{l}
\|\,\rho^h_\delta\,\|_{L^\infty\left((0,\infty);L^\infty(\Omega)\right)} \leq \|\,\rho_{0,\delta}\,\|_{L^\infty(\Omega)},\\ \\
\int_{\Omega_T}\rho^h_\delta\Big|\,\nabla\left(F^\prime(\rho^h_\delta)\right)\,\Big|^{q^\star} \leq \overline{M}(\Omega, T, F, \rho_{0,\delta}, q, \alpha),
\end{array}\right.
\end{equation}
and
\begin{equation}
\label{rem3}
\int_{\Omega_T} (\rho_0-\rho^h_\delta)\,\partial^h_t \xi + \int_{\Omega_T} \langle\,\rho^h_\delta \nabla c^\star\left[\,\nabla\left(F^{\prime}(\rho^h_\delta)\right)\,\right],\nabla \xi \,\rangle = 0(h^{\epsilon(q)}), 
\end{equation}
where $\epsilon(q):=\min(1,q-1),\,\xi\in C^2(\R\times\R^d)$, and
\[\overline{M}(\Omega, T, F, \rho_{0,\delta}, q, \alpha):= M(\alpha,q)\left(E_i(\rho_{0,\delta})-|\,\Omega\,|F\left(\frac{1}{|\,\Omega\,|}\right)+\alpha T|\,\Omega\,|\,\|\,\rho_{0,\delta}\,\|_{L^\infty(\Omega)}\right).\]
We introduce a convex function $H:[0,\infty)\rightarrow \R$, such that,
\begin{list}{({\bf HH\arabic{line}}) : }{\usecounter{line}}
\item $H\in C^1\left([0,\infty)\right)\cap C^2\left((0,\infty)\right)\;$ and $\;H^{\prime\prime}(x)=x^{1/q^\star}F^{\prime\prime}(x),\,\forall\,x>0$.  \label{h}
\end{list}
Combining (\ref{rem1}), (\ref{rem2}), and (HH\ref{h}), we have that 
\begin{equation}
\label{rem4}
\left\{\begin{array}{l}
\|\,\rho^h_\delta\,\|_{L^\infty\left((0,\infty);L^\infty(\Omega)\right)} \leq \|\,\rho_0\,\|_{L^\infty(\Omega)},\\ \\
\Big\|\,\nabla\left(H^\prime(\rho^h_\delta)\right)\,\Big\|^{q^\star}_{L^{q^\star}(\Omega_T)} \leq \overline{M}(\Omega, T, F, \rho_0, q, \alpha).
\end{array}\right.
\end{equation}
We deduce that there exists $\rho:[0,\infty)\times\Omega\rightarrow [0,\infty)$, such that $(\rho^h_\delta)_{h,\delta}$ converges to $\rho$ in $L^1(\Omega_T)$ for a subsequence, and $\{\tilde{\sigma}^h_\delta:=(\rho^h_{\delta})^{1/q}\nabla c^\star\left[\,\nabla\left(F^\prime(\rho^h_\delta)\right)\,\right]\}_{h,\delta}$ converges weakly to $\rho^{1/q}\nabla c^\star\left[\,\nabla\left(F^\prime(\rho)\right)\,\right]$ in $L^q(\Omega_T)$ for a subsequence, as $(h,\delta)$ goes to $(0,0)$. Then, we let $(h,\delta)$ go to $(0,0)$ in (\ref{rem3}), to conclude that
 $\rho\in L^\infty\left((0,\infty);L^\infty(\Omega)\right)$ is a weak solution of (\ref{intro8}), as in Theorem \ref{theoexistbound}; but here, we do not require that $\frac{1}{\rho}\in L^\infty\left((0,\infty);L^\infty(\Omega)\right)$ \hfill$\Box$\\

\noindent{\bf Case 2}: $\rho_0\in L^p(\Omega),\,p\geq q$, and $E_i(\rho_0)<\infty$.\\
Using Corollary 1.4.3 \cite{agueh:thesis}, we approximate $\rho_0$ by a sequence $(\rho_{0,\delta})_\delta$ in ${\cal P}_a(\Omega)$, such that
\[\left\{\begin{array}{l}
\eta_\delta \leq \rho_{0,\delta}\leq \epsilon_\delta\,\mbox{a.e.},\;\mbox{where}\;\; 0<\eta_\delta\leq \delta \;\mbox{and}\;\epsilon_\delta\geq \frac{1}{\delta}\\ \\
E_i(\rho_{0,\delta})\leq E_i(\rho_0)\\ \\
\rho_{0,\delta}\rightarrow \rho_0\;\mbox{in}\; L^p(\Omega)\;\mbox{as}\;\delta\downarrow 0,
\end{array}\right.\]
and we define $\rho^h_\delta$ as in case 1. Since $\rho_0\not\in L^\infty(\Omega)$, we cannot obtain (\ref{rem4}) from (\ref{rem2}), as in case 1. Here, we take advantage of the fact that $\rho_0\in L^p(\Omega)$ and $E_i(\rho_0)<\infty$, as follows:
\begin{enumerate}
\item We choose $G(x):=\frac{x^p}{p},\,x>0,$ in the (internal) energy inequality (\ref{en1}), and we observe that $\langle \nabla\rho^{p-1},\nabla c^\star\left[\nabla\left(F^\prime(\rho)\right)\right]\rangle \geq 0$ for $\rho\in {\cal P}_a(\Omega)$, to have that
\begin{equation}
\label{rem5}
\|\,\rho^h_\delta\,\|_{L^\infty\left((0,\infty);L^p(\Omega)\right)} \leq \|\,\rho_0\,\|_{L^p(\Omega)}.
\end{equation}
As a consequence, there exists a function $\rho: [0,\infty)\times\Omega\rightarrow [0,\infty)$ such that $(\rho^h_\delta)_{h,\delta}$ converges weakly to $\rho$ in $L^p(\Omega_T),\,0<T<\infty$.
\item Next, we choose $G:=F$ in (\ref{en1}), to control the spatial derivatives of $\rho^h_\delta$, as
\begin{eqnarray}
\label{rem6}
\lefteqn{\int_{\Omega_T} \rho^h_\delta\Big|\nabla\left(F^\prime(\rho^h_\delta)\right)\Big|^{q^\star}} \nonumber\\
& & \leq M(\alpha,q)\left(E_i(\rho_0)-|\,\Omega\,|F\left(\frac{1}{|\,\Omega\,|}\right)+\alpha|\,\Omega_T\,|^{1/p^\star}\|\,\rho_0\,\|_{L^p(\Omega)}\right).
\end{eqnarray}
\end{enumerate}
We combine (\ref{rem5}),(\ref{rem6}), and we use (HH\ref{h}), to deduce, as in the previous sections, that $(\rho^h_\delta)_{h,\delta}$ converges strongly to $\rho$ in $L^1(\Omega_T)$, for a subsequence. We conclude, as in case 1, that $\rho$  is a weak solution of (\ref{intro8}) in the sense that, $\rho\in L^\infty\left((0,\infty);L^p(\Omega)\right)$,\\
 $\rho\nabla c^\star\left[\nabla\left(F^\prime(\rho)\right)\right]\in L^1(\Omega_T)$ for $0<T<\infty$, and (\ref{theobound0}) holds \hfill$\Box$\\

\noindent Few examples of energy density functions satisfying (HF\ref{f1}) - (HF\ref{f2}) and (HH\ref{h}) are $F(x)=\sum_{i=1}^n A_iF_i(x)$, where $F_i(x)\in \{x\ln(x), \frac{x^m}{m-1}\}$ with $m>1$ or $\max\left(\frac{1}{q}, 1-\frac{1}{d}\right)\leq m<1$, and $A_i>0$. For examples, for the fast diffusion equation $\frac{\partial\rho}{\partial t}=\Delta\rho^m$, this corresponds to the range $1-\frac{1}{d}\leq m<1$ if $d\geq 2$, and $\frac{1}{2}\leq m<1$ if $d\leq 2$; for $\frac{\partial\rho}{\partial t}=\Delta_p\rho^n$, we have $n\geq\frac{d-(p-1)}{d(p-1)}$ if $d\geq p$, and $n\geq\frac{1}{p(p-1)}$ if $d\leq p$, and in particular, for the p-Laplacian $\frac{\partial\rho}{\partial t}=\Delta_p\rho$, we require $\frac{2d+1}{d+1}\leq p\leq d$ or $p\geq \max\left(d,\frac{1+\sqrt{5}}{2}\right)$.

\section{Appendix}
In Proposition \ref{propappendix1}, we collect results of previous authors used in this work, and in Proposition \ref{lemappend3}, we establish intermediate results needed in the previous sections. Proposition \ref{propappendix1} is due to Cordero \cite{cordero:properties}  and Otto \cite{otto:doubly} . For its proof, we refer to these references. A sketch of proof of this proposition can also be found in \cite{agueh:thesis}, sections 5.1 and 5.2.

\begin{prop}
\label{propappendix1}
Let $\rho_0, \rho_1\in {\cal P}_a(\Omega)$, and assume that $c:\R^d\rightarrow [0,\infty)$ is strictly convex, and satisfies $c, c^\star\in C^2(\R^d)$. Denote by $S$, the $c$-optimal map that pushes $\rho_1$ forward to $\rho_0$, and define the interpolant map $S_t$, and the interpolant measure $\mu_{1-t}$, by
\[ S_t:=(1-t)\mbox{id}+tS\quad \mbox{and}\quad \mu_{1-t}:=(S_t)_{\#}\rho_1,\]
for $t\in [0,1]$. Then,
\begin{enumerate}
\item $S_t$ is injective, and $\mu_{1-t}$ is absolutely continuous with respect to Lebesgue.\\
Moreover, there exists a subset $K$ of $\Omega$, of full measure for $\mu_1:=\rho_1(y)\mbox{d}y$, such that, for $y\in K$ and $t\in [0,1]$, 
\item $\nabla S(y)$ is diagonalizable with positive eigenvalues.
\item The pointwise Jacobian $\mbox{det}\,(\nabla S)$ satisfies
\[0 \not= \rho_1(y) = \rho_{1-t}\left(S_t(y)\right)\,\mbox{det}\,\left[\,(1-t)\,\mbox{id}\,+t\nabla S(y)\,\right],\]
where $\rho_{1-t}$ is the density function of $\mu_{1-t}$.\\

In addition, if $\rho_1>0$ a.e., then
\item the pointwise divergence $\mbox{div}\,(S)$ is integrable on $\Omega$, and
\[\int_\Omega \mbox{div}\,(S(y)-y)\,\xi(y)\,\mbox{d}y \leq -\int_\Omega \langle\,S(y)-y,\,\nabla \xi\,\rangle\,\mbox{d}y,\]
for $\xi\geq 0$ in $C^\infty_c(\R^d)$.
\end{enumerate}
\end{prop}

\noindent The following estimates will be needed in the previous sections.

\begin{prop}
\label{lemappend3} Assume that $c:\R^d\rightarrow [0,\infty)$ is strictly convex, of class $C^1$, and satisfies $c(0)=0$ and (HC\ref{c3}). Then
\begin{equation}
\label{appendn3}
\langle\,z,\nabla c^\star(z)\,\rangle\geq c^\star(z) \geq 0,\;\;\forall\,z\in\R^d.
\end{equation}
In addition, if $\,c(z)\geq \beta\,|\,z\,|^q$ for some $\beta>0$ and $q>1$, then
\begin{equation}
\label{appendn6}
\langle\,z,\nabla c^\star(z)\,\rangle \leq M(\beta, q)\,|\,z\,|^{q^\star},
\end{equation}
where $M(\beta, q)$ is a constant which only depends on $\beta$ and $q$.
\end{prop}
{\tt Proof:} Since $c$ is strictly convex, differentiable and satisfies (HC\ref{c3}), we have that $c^\star\in C^1(\R^d)$ is convex. Then,
\begin{equation}
\label{appendn7.0}
\langle z, \nabla c^\star\,(z)\rangle = c^\star(z)+c\left(\nabla c^\star(z)\right) \geq  c^\star(z).
\end{equation}
Because $c(0)=0$ and $0$ minimizes $c$, we have that $c^\star(0)=0$ and $0$ minimizes $c^\star$. We conclude that $c^\star(z)\geq 0$, which proves (\ref{appendn3}).\\
\noindent Now, assume that $c(z)\geq \beta\,|\,z\,|^q$. Since $c^\star\in C^1(\R^d)$ is convex and nonnegative, we have that 
\begin{equation}
\label{appendn8.0}
\langle\,z,\,\nabla c^\star(z)\,\rangle \leq c^\star(2z)-c^\star(z) \leq c^\star(2z).
\end{equation}
Moreover, because  $c(z)\geq \beta\,|\,z\,|^q$, we have that 
\begin{equation}
\label{append8.1}
c^\star(2z)\leq M(\beta, q)\,|\,z\,|^{q^\star}.
\end{equation}
We combine (\ref{appendn8.0}) and (\ref{append8.1}), to conclude (\ref{appendn6}) \hfill$\Box$\\

\noindent{\bf \large Acknowledgements.} It is a pleasure to express my profound gratitude to my Ph.D advisor Wilfrid Gangbo, and to C\'edric Villani, Andrzej \'Swi\c ech and Eric Carlen for their support and suggestions.

\end{document}